%% file: main.tex
\documentclass{article}

\usepackage{hyperref}
\usepackage{verbatim}
\usepackage{graphicx}
\usepackage{amsmath}
\usepackage{amssymb}
\usepackage{amsthm}
\usepackage{color}

\usepackage{pgf}
\usepackage{tikz}
\usetikzlibrary{arrows,positioning,backgrounds,patterns,decorations.pathreplacing}

\newtheorem{theorem}{Theorem}[section]
\newtheorem{lemma}[theorem]{Lemma}
\newtheorem{proposition}[theorem]{Proposition}

\theoremstyle{definition}
\newtheorem{definition}[theorem]{Definition}

\theoremstyle{definition}

\numberwithin{equation}{section}

\include{macros}

\begin{document}

\title{On representations of
real numbers and the computational complexity of converting between such representations}

\linespread{1.1}

\author{Amir M.~Ben-Amram%
\thanks{Qiryat Ono, Israel}
\and
Lars Kristiansen%
\thanks{Department of Mathematics, University of Oslo, Norway
and
    Department of Informatics, University of Oslo, Norway}
\and
Jakob Grue Simonsen%
\thanks{Department of Computer Science, University of Copenhagen (DIKU), Denmark}}

\maketitle
\date{}

%    Abstract is required.
\begin{abstract}
We study the computational complexity of converting one representation of real numbers
into another representation. Typical examples of representations are Cauchy sequences, base-10 
expansions, Dedekind cuts and continued fractions.
\end{abstract}


\input{introduction}

\input{preliminaries}

\input{weihrauch}

\input{cauchy}

\input{sumeapproxbaseexp}

\input{bolzano}

\input{bestapproxequiv}

\input{continued}

\bibliographystyle{plain}
\bibliography{larsjakob.bib}

\end{document}

%% file: macros.tex
 \newcommand{\xcb} {{\mathcal B}}

\newcommand{\xcc} {{\mathcal C}}
\newcommand{\xcd} {{\mathcal D}}

\newcommand{\redrel}{ \preceq_S }
\newcommand{\redrelstrict}{\prec_S}

\newcommand{\dga}{\mbox{{\bf a}}}
\newcommand{\dgb}{\mbox{{\bf b}}}

\newcommand{\xd}{\texttt{D}}

\newcommand{\primeset}[1]{{\mathrm p}{\mathrm r}{\mathrm i}{\mathrm m}(#1)}

\newcommand{\xor}{\oplus}

\newcommand{\integer}{\ensuremath{\mathbb Z}}
\newcommand{\rational}{\ensuremath{\mathbb Q}}
\newcommand{\nat}{\ensuremath{\mathbb N}}

\newcommand{\paramorac}[2]{#1^{#2}}

\newcommand{\ddiv}[2]{#1 \, \mathrm{div} \, #2}

\newcommand{\mmod}[2]{#1 \, \mathrm{mod} \, #2}

\newcommand{\polylog}[1]{\mathrm{polylog}(#1)}
\newcommand{\poly}[1]{\mathrm{poly}(#1)}

\newcommand{\fareypairtree}[0]{\mathcal{T}_{\mathrm{F}}}
\newcommand{\sternbrocot}[0]{\mathcal{T}_{\mathrm{SB}}}

\newcommand{\bitlen}[1]{\Vert #1\Vert}

%% file: introduction.tex
\section{Introduction} \label{kortstokk}

\subsection{Motivations.}
In a computational setting
real numbers can be represented by Cauchy sequences, base-10 expansions, Dedekind cuts, continued fractions and a number of other representations (we will consider quite a few of them in this paper). Our goal is to analyze the computational complexity of converting one representation into another.
Let us say that we have access to the Dedekind cut of the real number $\alpha$. How hard will it be to compute a Cauchy sequence for $\alpha$?
How hard will it be to compute the continued fraction of $\alpha$?
Or let us say that we have access to the continued fraction of $\alpha$, 
how hard will it then
be to compute the base-10 expansion of $\alpha$? Will there be an efficient
algorithm? Can it be done in polynomial time? Exponential time?
%How much space will be required? (we don't actually discuss space - ABA)

These are very natural questions to ask, but they are also
naive, and the way the questions are posed above, does
hardly make any sense at all. We will do our best to pose such questions in a mathematically satisfactory manner, and when it is possible, we will
derive reasonably tight upper bounds on the computational complexity of conversions.  We cover most of the classic and well-known representations,  and we will also consider some representations of newer date.

It might very well be the case that there is an algorithm for converting one representation
in to another, but still it will not be possible to derive any upper bounds
on the computational resources the algorithm requires. In such situations
we will  give an intuitive explanation of why this is so (full proofs are mostly omitted as the results
can be found in references \cite{ciejouren} \cite{ciejourto} \cite{apalilf}).

\subsection{What is a representation?}

Formally, a {\em representation} of the {\em irrational} numbers
will be a class of functions. Every function in the class will represent a particular irrational number, and each
irrational number will be represented by some function in the class. The class of all Dedekind cuts 
of irrational numbers will be a canonical representation to us: The Dedekind cut of an irrational $\alpha$ is the function $\alpha:\rational \longrightarrow \{0,1\}$ where
$$
\alpha(q)=  \begin{cases}
0   & \mbox{if $q < \alpha$}   \\
1   & \mbox{if $q > \alpha$.}
\end{cases}
$$
Each irrational number has a unique representation in this class, and we can identify an irrational number $\alpha$
with its Dedekind cut $\alpha:\rational \longrightarrow \{0,1\}$.

We  will take advantage of the uniqueness of the Dedekind cuts to define what a representation in general is.
We  refer to the functions in a representation $R$ as
$R$-representations. When $f$ is an $R$-representation of 
$\alpha$, we will require that it is
 possible to compute
the Dedekind cut of $\alpha$ in $f$, that is, we will
require that there exists an oracle Turing machine $M$
such that
$$  \alpha(q) = \Phi^f_M(q)        $$
where $\Phi^f_M$ is the function computed by $M$ with oracle $f$. We will also require that at least one
$R$-representation $f$
of  $\alpha$ can be computed in the Dedekind cut of $\alpha$, that is, we will
require that there exists an oracle Turing machine $N$
such that
$$  f(x) = \Phi^\alpha_N(x)        $$
where $\Phi^\alpha_N$ is the function computed by $N$ when
the oracle is the Dedekind cut of $\alpha$.
We are now ready to give our
 formal definition.

\begin{definition} \label{begravelseimorgen}
A class of functions $R$ is a {\em representation (of
the irrational numbers)}  if
there exist Turing machines $M$ and $N$ such that
\begin{enumerate}
    \item for every irrational $\alpha$ there exists $f\in R$
    such that
    $$  \alpha = \Phi^f_M  \mbox{ and } f= \Phi^\alpha_N $$
\item for every $g\in R$ there exist an irrational $\alpha$  such that
$$ \alpha = \Phi^g_M = \Phi^f_M \mbox{ where } f= \Phi^\alpha_N\;$$
 When $\alpha = \Phi^g_M$, we say that $g$ {\em represents} $\alpha$ and that $g$ is an {\em $R$-representation} of $\alpha$.
\end{enumerate}

We say that an oracle Turing machine $M$
{\em converts an $R_1$-representation into an $R_2$-representation} if there 
for any $f\in R_1$ representing $\alpha$ exists $g\in R_2$ representing $\alpha$ such that $g = \Phi_M^f$.
\qed
\end{definition}

Let us study a few examples  in order to see how this 
definition works. 
We define a Cauchy sequence for 
  $\alpha$ as a function 
$C:\nat^{+} \longrightarrow \rational$ with the property
$$ \vert C(n) - \alpha \vert < n^{-1}\; .
$$ \label{flystreik}
Let $\mathcal{C}$ be the class  of all Cauchy sequences for all irrational numbers. We will now argue that  $\mathcal{C}$ is a representation according to the definition above.

First we observe that we can compute a Cauchy sequence
$C$ for $\alpha$ if we have access to the Dedekind cut of $\alpha$. We can use the Dedekind cut to find an integer $a$
such that $a< \alpha < a+1$. Thereafter, we can use the Dedekind cut and the equations
$$ C(1)= a+ \frac{1}{2} \;\;\; \mbox{ and } \;\;\;
C(i+1)=  \begin{cases}
C(i)- 2^{-i-1}  & \mbox{if $C(i) > \alpha$}   \\
C(i)+ 2^{-i-1}  & \mbox{if $C(i)< \alpha$}
\end{cases}
$$
to compute $C(n)$ for  arbitrary $n$. This is one possible
way to compute a Cauchy sequence in a Dedekind cut. There are
for sure other ways. Other algorithms may yield
Cauchy sequences that converge faster, or  slower, than the 
ones computed by
the algorithm suggested above. Anyway, there will be
an oracle Turing machine $N$ such that we have 
$f = \Phi^\alpha_N$ where $f$ is some Cauchy sequence for $\alpha$.

Next we observe that can compute the Dedekind cut of an irrational $\alpha$ in any Cauchy sequence for $\alpha$.
In order to compute $\alpha(q)$, we search for the least $n$ such that
$\vert C(n)- q\vert > n^{-1}$. This search terminates as
$q$ is rational and $\alpha$ is irrational. If $q<C(n)$,
it will be the case that $\alpha(q)=0$ (we have $q<\alpha$),
otherwise, we have $q>C(n)$, and then it will be case that  
$\alpha(q)=1$ (we have $q>\alpha$). Thus there will be
an oracle Turing machine $M$ such that $\alpha= \Phi^f_M$
whenever $f$ is a Cauchy sequence for $\alpha$. 

This shows that $\mathcal{C}$, that is, the class of all Cauchy sequences for all irrational numbers, is a representation.
We have Turing machines $M$ and $N$ such that there
for any irrational $\alpha$ exists a Cauchy sequence $C$
such that
  $\alpha = \Phi^C_M$  and  $C = \Phi^\alpha_N$, and
thus, clause (1) of Definition \ref{begravelseimorgen} is
satisfied, moreover, clause (2) is also satisfied since we have 
$\alpha = \Phi^g_M = \Phi^C_M$
for any Cauchy sequence $g$ for $\alpha$.

Let  $a_0, a_1, a_2, \ldots $ be an infinite sequence of integers where
$a_1, a_2, a_3 \ldots $ are positive. 
The   {\em continued fraction} $[a_0; a_1, a_2, \ldots ]$ is defined by
 $$[ \ a_0; a_1, a_2,a_3 \ldots \ ]  \; = \;
a_0 + \frac{\displaystyle 1}{\displaystyle a_1 + \frac{\displaystyle 1}{\displaystyle a_2 + \frac{1}{a_3 + \ldots}}}
$$
It is well known that any irrational number
can be uniquely written as an infinite continued
fraction, and moreover, each infinite continued fraction equals an irrational number (rationals have finite continued fractions). There is a one-to-one correspondence between the infinite continued fractions and the irrational
numbers. 

Let $\mathcal{F}$ be the class of all infinite continued fractions
where each $[ \ a_0; a_1, a_2,a_3 \ldots \ ]$ is identified
with a function $f$ where $f(n)=a_n$. Then $\mathcal{F}$ will be a representation according to Definition
 \ref{begravelseimorgen}. If we have 
 access to the Dedekind cut of $\alpha$,  we can compute the continued fraction of $\alpha$, and if have access to the continued fraction of $\alpha$, we
can compute the Dedekind cut of $\alpha$. Thus, there exist Turing machines $M$ and $N$ such that
$$  \alpha = \Phi^f_M  \mbox{ and } f= \Phi^\alpha_N $$
whenever $f$ is the continued fraction of the irrational number $\alpha$. This shows that clause (1) of 
Definition \ref{begravelseimorgen} is fulfilled. Clause (2)
is trivially fulfilled because of the uniqueness of the continued fractions, that is, because there is a one-to-one correspondence between the infinite continued fractions and the Dedekind cuts of the irrational numbers.

We will not give a formal definition of a representation
of all real numbers, and it is essential that Definition \ref{begravelseimorgen}
is restricted to the irrational numbers.
 If we involve
the rationals, the definition will not serve its
purpose as we cannot
always uniformly convert one standard representation into another,
even if we are dealing with representations 
of a computable nature. E.g., the algorithm above converting 
a Cauchy sequence $C$ into a Dedekind cut, 
search for a number $n$ such that 
$\vert C(n)- q\vert > n^{-1}$. 
This search will not terminate when $C$
is a Cauchy sequence for the rational number $q$. Thus,
if $\alpha$ might be rational,
the algorithm does not yield a Turing machine 
$M$ such that $\Phi^C_M$ is the Dedekind cut of $\alpha$
whenever $C$ is a Cauchy sequence $\alpha$, moreover,
it can be proven that no such Turing machine $M$ exits
(see Mostowski \cite{most}).
We {\em cannot} uniformly convert  Cauchy sequences 
for {\em real} numbers into  Dedekind cuts,
but we {\em can} uniformly convert  Cauchy sequences for 
{\em irrational} numbers  into Dedekind cuts.

The purpose of Definition \ref{begravelseimorgen} is to 
capture what we  intuitively  consider as  computable
representations of the real numbers, and maybe somewhat paradoxically, we achieve  just that by restricting the
definition to the irrational numbers. The reader should
be aware that standard notions of representations of
reals and irrationals in computable analysis tend to be more general than our notion, see Weihrauch \cite{wei}.
E.g., a sequence 
of rationals $q_0, q_1, q_2, \ldots$ containing
all rational numbers less than  $\alpha$, and nothing but
the rationals less than  $\alpha$, 
 will not yield a
representation according to our definitions.
We cannot use such a sequence 
 to compute the Dedekind cut of $\alpha$ (even if $\alpha$ is irrational). Such sequences will typically  be considered as a representation in the literature of computable analysis.

\subsection{An ordering relation on representations. }
We have seen that there is an algorithm for computing
the Dedekind cut of an irrational $\alpha$ in an arbitrary
Cauchy sequence for $\alpha$. The algorithm searches for the
least natural number $n$  that fulfills  certain criteria.
It is easy to see that such an unbounded search
is necessary. We cannot
 convert a Cauchy 
sequence into a Dedekind cut if we are not allowed to
carry out unbounded search. Neither can we convert a base-10
expansion into a base-3 expansion if we are not allowed to
carry out unbounded search.
Suppose an oracle tells us that the base-10 expansion of
an irrational starts
with $0.66666$. That will not be enough for us to decide if
the base-3 expansion starts with $0.1$ or $0.2$. 
Thus, in order to determine which of these two options
we should pick,
we have to ask the oracle for the next digit of the decimal
expansion, but of course, the next digit might also be 6, and so might the
next after the next. The oracle may continue for an arbitrarily long time to tell us that the
next digit  is $6$. Since the number is irrational the oracle will eventually  yield a digit that allows to determine if the base-3 expansion  starts with
$0.1$ or $0.2$, but we need to carry out an unbounded  search 
to get that digit.

Algorithms (conversions, computations, etc.) that do not perform unbounded
search will be referred as {\em subrecursive} algorithms
(conversions, computations, etc.), in general, the word {\em subrecursive} signifies {\em absence of unbounded search}. This terminology might not be
standard, but it will be very convenient.

We will now define an ordering relation $\redrel$
over the representations. Intuitively, the relation 
$R_1 \redrel R_2$ will indicate that the representation
$R_2$ is more informative than the representation $R_1$. 
If $R_1 \redrel R_2$ holds, a Turing machine with oracle access to an $R_2$-representation
of $\alpha$ can
subrecursively (yes, that means without carrying out unbounded search) compute 
an $R_1$-representation of $\alpha$. 
Thus, if the relation $R_1 \redrel R_2$ does not hold, it will not make much sense to talk about the computational complexity of converting an $R_2$-representation into an $R_1$-representation as such a conversion requires unbounded search, and thus, there will be no upper bound on the running time of a Turing machine undertaking the  conversion. On the other hand, if the relation holds, it should make sense to analyze the computational complexity of the conversion.  So far we have just indicated our intention
with the relation; below we shall give the formal definition.
First, however, we need an auxiliary definition. We are going to
formulate our definition in terms of time-bounded computation, and there is need to specify what
functions we admit as time bounds. We admit \emph{time-constructible functions}, as defined next. This
is a standard choice in complexity theory.

\begin{definition} 
A function $t:\nat\longrightarrow \nat$ is a {\em time bound}
if (i) $n\leq t(n)$, (ii) $t$ is increasing and (iii) $t$ is time-constructible:
there is a single-tape Turing machine that, on input $1^n$, computes
$t(n)$ in $\Theta(t(n))$ steps.
\qed
\end{definition}

Clause (iii) in the definition is needed because there are functions whose computational complexity
is disproportionate to the size of their values. We exclude such functions as time bounds, avoiding
certain pitfalls in proofs. The class of functions we admit as time bounds includes all the functions familiar
from analysis of algorithms such as polynomials (with positive coefficients), exponentials, the
tower-functions etc.

\begin{definition} \label{snartsommerigjen}
Let $t$ be a time-bound and let $R$ be a representation.
Then, $O(t)_R$ denotes the class of all irrational $\alpha$
in the interval $(0,1)$ such that at least one $R$-representation of $\alpha$ is computable by a Turing machine running in time 
$O(t(n))$ (where $n$ is the length of the input).  

Let $R_1$ and $R_2$ be representations. The
 relation $R_1 \redrel R_2$ holds if there for any
 time-bound $t$ exists a time-bound $s$ such that 
$$
O(t)_{R_2} \; \subseteq \; O(s)_{R_1}\; .
$$
If the relation $R_1 \redrel R_2$ holds, we will say that the representation $R_1$ is {\em subrecursive} in the representation $R_2$.
\qed 
\end{definition}

We will now study a few examples and discuss how the
definition above works. 
Recall that $\xcc$ denotes the representation 
by Cauchy sequences (see page \pageref{flystreik}). Let  
$\xcd$ denote the representation by Dedekind cuts.
It turns out that we have $\xcc\redrel \xcd$ and 
$\xcd\not\redrel \xcc$. Let us see why.

Intuitively we have $\xcc\redrel \xcd$ because
a Cauchy sequence for an irrational $\alpha$ in the interval
$(0,1)$ can be subrecursively computed in the Dedekind cut
of $\alpha$. No unbounded search is required. 
We simply set
$C(1)$ equal to $1/2$, and then we use the Dedekind cut
and, e.g., the equations
$$ 
C(i+1)=  \begin{cases}
C(i)- 2^{-i-1}   & \mbox{if $C(i) > \alpha$}   \\
(C(i)+ 2^{-i-1}  & \mbox{if $C(i)< \alpha$}
\end{cases}
$$
to compute $C(n)$. Then $C:\nat^+ \longrightarrow \rational$
will be a Cauchy sequence for $\alpha$. Formally we have
$\xcc\redrel \xcd$ because there for every time-bound $t$
exists a time-bound $s$ such that 
$O(t)_{\xcd} \; \subseteq \; O(s)_{\xcc}$. A Turing machine
that uses the equations above to compute  $C(n)$ needs to
compute the Dedekind cut  $O(n)$ times, that is,
$O(2^{\bitlen{n}})$ times where $\bitlen{n}$ is the
length of the input. Furthermore, assuming numbers are  represented 
in binary form, the Turing machine only needs to compute the Dedekind
cut for inputs of  size $O(2^{\bitlen{n}})$. Hence, if the Dedekind
cut of $\alpha$ is computable in time $O(t(m))$ where $m$ is the length of the input and $t$ is a time-bound,
then a Turing machine can compute a Cauchy sequence for $\alpha$
in time $O( nt( k 2^{\bitlen{n})})$ for some constant $k$.
Thus, let  $s$ be the time-bound $s(\bitlen{n})=  nk't2^{\bitlen{n}}$, for appropriate $k'> k$
and we have 
$O(t)_{\xcd}  \subseteq  O(s)_{\xcc}$.

Note that the set $O(t)_{\xcd}$ only contains irrationals from the interval
$(0,1)$. This is important.  When we compute a Cauchy sequence $C$
for an irrational in this interval, we can simply let $C(1)=1/2$.
We cannot in general compute a Cauchy sequence for an irrational $\beta$
subrecursively in the Dedekind cut of $\beta$. 
In order to set the value of $C(1)$ 
we will need a rational $a$ such that $a<\beta< a+1$, and we cannot 
get hold of such an $a$
without resorting to unbounded search. Thus, if we do not restrict
$O(t)_{\xcd}$ to irrationals in the interval $(0,1)$,
the relation $\xcc\redrel \xcd$ will not hold. But we want it to hold.
We are interested in representations of the fractional part of irrational
numbers. It is natural to abstract away the integral part.

Now, let us discuss why we have $\xcd \not\redrel \xcc$.
As explained above, we have $\xcd \not\redrel \xcc$ because
we cannot avoid unbounded search when we compute a Dedekind cut in a
Cauchy sequence, but what does a formal proof look like? In general
it is much harder to prove that the relation $\redrel$ does
not hold than it is to prove that it holds. In order to 
prove $R_1 \redrel R_2$, we just have to come up with a subrecursive
algorithm for converting an $R_2$-representation into an $R_1$-representation. In order to 
prove $R_1 \not\redrel R_2$, we have, according to our definitions,
to prove that there is a time-bound $t$ such that we have
$O(t)_{R_2} \not\subseteq O(s)_{R_1}$  for any time-bound
 $s$. This might not be all that easy. This might call for involved
 diagonalization arguments. In some cases we might do with a growth
 argument, \label{growthpage} that is, we might be able to prove,  for some time-bound $t$,
that there for any time-bound $s$ exists $\alpha_s\in O(t)_{R_2}$
such that any Turing machine computing an $R_1$-representation of $\alpha_s$
will have to give very large outputs, that is, outputs whose length is
not of order $O(s)$ (and thus the Turing machine cannot run in time $O(s)$).
Growth arguments tend to be easier, or at least less tedious, than 
diagonalization arguments, and we will present a rather detailed 
proof based on a growth argument in
Section \ref{skalmotemartintap10}.
But we cannot prove $\xcd \not\redrel \xcc$ by a growth
argument as a Turing machine computing 
a Dedekind cut gives outputs of length 1. Let us see how we can prove
$\xcd \not\redrel \xcc$ by a diagonalization argument.

%Let us see how we can prove
%$\xcd \not\redrel \xcc$ by a diagonalization argument (a growth
%argument will not be available as a Turing machine computing 
%a Dedekind cut gives outputs of length 1).

Let $s$ be an arbitrary time-bound. We will,  by standard diagonalization techniques, construct
a Cauchy sequence for an irrational $\alpha$ in the interval $(0,1)$
such that $\alpha$ becomes different from each 
$\beta\in O(s)_\xcd$, and hence we will have 
$\alpha\not\in O(s)_\xcd$. Our construction can be carried
out by a Turing machine, that is, the Cauchy sequence for $\alpha$ can be computed by
a Turing machine. That
Turing machine will run in time $O(t)$ for some 
time-bound $t$,  and thus we
have $\alpha\in O(t)_\xcc$. Moreover, it will turn out
that $t$ does not depend
on $s$.  Hence, we have a time-bound  $t$ and for every time-bound $s$ 
there exists $\alpha$ such that $\alpha\in O(t)_\xcc$
and $\alpha\not\in O(s)_\xcd$. Hence, we have  $t$
such that  $O(t)_\xcc \not\subseteq O(s)_\xcd$
for every $s$, and thus, by Definition \ref{snartsommerigjen},
we can conclude that $\xcd\not\redrel \xcc$.

We will now give an algorithm for computing a  Cauchy sequence $C$ such that $\lim_n C(n) \not\in O(s)_\xcd$
where $s$ is an arbitrary time-bound.
We will need a standard enumeration  $\{ e \}_{e\in \nat}$  of the Turing machines, 
and we use $\{ e\}(x)$ to denote the execution of
the $e$'th Turing machine on input $x$. Furthermore, we 
need an increasing time-bound function $S$ that eventually dominates any
time-bound of order $O(s)$, that is, for any $s_0$ of order
$O(s)$ we have $s_0(m)< S(m)$ for all sufficiently large $m$.
Such an $S$ will always exist, and moreover can be chosen so that $S(n)\ge 2n$ for all $n$.
Note that by the definition of  time-bound functions, we have a Turing machine that 
given $n$, computes $S(n)$ in at most $a_S\cdot S(n)$ steps for some constant $a_S$.
We define the sequence
$d_0,d_1,d_2, \ldots$ by $d_0=2$ and $d_{i+1} = S(d_i)$.
Finally, we will need a standard computable bijection 
$\langle \cdot , \cdot\rangle: \nat \times \nat \longrightarrow \nat$. For any $i\in\nat$, 
our algorithm needs to
compute the
unique $j,e\in \nat$ such that $\langle j , e \rangle= i$.

The algorithm sets $C(1):= 1/2$. If $n>1$, the algorithm
checks if there exists $i$ such that $n = d_{3i+2}$.
If such an $i$ does not exist, the algorithm simply
set $C(n):=C(n-1)$;  if such $i$  exists
the algorithm finds the unique $j,e$ such that
$\langle j,e\rangle =i$ and sets
    \begin{itemize}
        \item $C(n):= C(n-1)$ if $\{ e \}(C(n-1))$ does not terminate within $n$ steps
        \item $C(n):= C(n-1) - 2^{-n}$ if $\{ e \}(C(n-1))$   terminates within $n$ steps and outputs 0
        \item $C(n):= C(n-1) + 2^{-n}$ if $\{ e \}(C(n-1))$   terminates within $n$ steps and outputs something else than 0.
    \end{itemize}

It is clear that the algorithm indeed computes a Cauchy sequence for a real number
in the interval $(0,1)$, and it is also pretty easy to see that
the length of the output $C(n)$ will be
bounded by a function of order $O(n)$ if we represent
numbers in binary form and code rationals in a reasonable way. We can w.l.o.g.~assume that $S$ 
eventually
dominates any function of order $O(n)$ (we can just pick an $S$ that increases fast enough).
Thus, for all sufficiently large $n$, we have 
\begin{align} \label{snartslss}
    \bitlen{C(n)} < S(n)\; .
\end{align}
Now, let $\alpha= \lim_{n} C(n)$ and let $\beta\in O(s)_\xcd$. We will prove that $\alpha\neq \beta$.

The Dedekind cut of $\beta$ can be computed by a
Turing machine $\{ e \}$ running in time $O(s)$.
Thus, $S(\bitlen{x})$ will be an upper bound on the number
of steps in the computation $\{ e \}( x  )$
when $x$ is large. Pick a sufficiently large $j$ and let 
$i = \langle j,e\rangle$ (we can make $d_{3i}$  as big as
we want by picking a big $j$). Our algorithm is designed such that we have 
\begin{align} \label{lyongadestenger}
    C(d_{3i}) = C(d_{3i+2}-1)\; .
\end{align}
Hence, by (\ref{lyongadestenger}) and (\ref{snartslss}),
the number of steps in the computation
$\{ e \}( C(d_{3i+2}-1)  )$ will be bounded by
$$S(\vert C(d_{3i+2}-1) \vert) 
= S(\vert C(d_{3i}) \vert) <  S(S(d_{3i})) = d_{3i+2}\; .
$$
Assume that the output of the computation 
$\{ e \}( C(d_{3i+2}-1)  )$ is 0. 
Then we have
$$
C(d_{3i+2}) = C(d_{3i+2}-1) + 2^{-d_{3i+2}}
$$
but as $\{ e \}$ computes
the Dedekind cut of $\beta$, we have $\beta > C(d_{3i+2}-1)$.
It follows that
$$
\alpha = \lim_n C(n) < C(d_{3i+2}-1) - 2^{-d_{3i+2}} + \sum_{n> d_{3i+2}} 2^{-n} < C(d_{3i+2}-1) < \beta\; .
$$
Assume that the output of the computation 
$\{ e \}( C(d_{3i+2}-1)  )$ is different from 0. 
Then we have
$$
C(d_{3i+2}) = C(d_{3i+2}-1) + 2^{-d_{3i+2}}
$$
but as $\{ e \}$ computes
the Dedekind cut of $\beta$, we have $\beta < C(d_{3i+2}-1)$.
It follows that
$$
\alpha = \lim_n C(n) > C(d_{3i+2}-1) + 2^{-d_{3i+2}} - \sum_{n> d_{3i+2}} 2^{-d_{n}} > C(d_{3i+2}-1) > \beta\; .
$$
This proves that $\alpha\neq \beta$.

The same argument can be used to prove that the limit $\alpha$ is irrational. Among the Turing
machines enumerated there is a machine $e$ that computes the Dedekind cut of any given rational $q$; this 
computation can be done in linear time under a reasonable encoding of rationals,
so we may assume that $e$ has running time in $O(s)$; and
we conclude that $\alpha \neq q$.

We have proved that our algorithm 
computes a Cauchy sequence for an irrational number
and that irrational number cannot be in the class $O(s)_{\mathcal{D}}$. Let us 
undertake a 
complexity analysis of a Turing machine $M$ executing
the algorithm:
\begin{itemize}
    \item The input to $M$ is a natural number $n$ (we will estimate an upper bound for $M$'s
    running time as a function of $n$).
    \item First $M$ will compute $C(n-1)$.
    \item Then $M$ will check if there exists $i$
    such that $n=d_{3i+2}$. Recall that $d_0=2$ and 
    $d_{i+1}= S(d_i)$ where $S$ is a time-bound
    function. 
    It is possible to check if such
    an $i$ exists in time $O(n)$. Briefly, $M$ computes $d_0,d_1,d_2\dots$
    until it either hits $j$ such that $n=d_j$, or $j$ such that $n<d_j$, or $j$ such that
    the computation of $d_j$ exceeds $a_S n$ time. In all cases, the computation of the last $d_j$
    is either completed or stopped after $O(n)$ steps. The computation of previous elements of the
    sequence also take $O(n)$ because the sequence, and therefore its computation time, grows at 
    least geometrically.
    \item If $M$ finds $i$ such that $n=d_{3i+2}$,
    then $M$ will compute $j,e$
    such that $\langle j,e\rangle = i$, check if the
    computation $\{ e \}(C(n-1))$ terminates within $n$ 
    steps, and  finally, compute the output.
    All this can be done in time $O(n^2)$ on a multi-tape Turing machine.
    %(who cares - ABA) and hence in time
    %$O(n^4)$ on a single-tape machine.
\end{itemize}
These considerations show that $M$ runs in time $O(n^{2})$, where $n\in \nat$ is
the input, and thus in time $O(2^{2\bitlen{n}})$ where
$\bitlen{n}$ is the length of the input (time complexity is always 
stated as a function of the input bit-length).
This allows us to conclude that $M$
computes a Cauchy sequence for an irrational in the
class $O(2^{2\bitlen{n}})_{\mathcal{C}}$. We also know
that this irrational is not in the class 
$O(s)_{\mathcal{D}}$, and recall that  $s$ was an arbitrary chosen time-bound. Hence, we have 
$O(2^{2\bitlen{n}})_{\mathcal{C}} \not\subseteq O(s)_{\mathcal{D}}$ for any
time-bound $s$. This proves that $\xcd \not\redrel \xcc$.

Our proof that $\xcd \not\redrel \xcc$ is meant to illustrate how our definitions works.
In the current paper we will in general not formally prove that one representation is
not subrecursive in another, but for the benefit of the reader we will to a certain extent 
provide informal explanations and intuitive arguments of why subrecursive conversions between
certain representations are impossible.

\begin{definition} \label{rodtihjornet}
Let $R_1$ and $R_2$ be representations. The relation 
$R_1 \equiv_S R_2$ holds when $R_1 \redrel R_2$ and
$R_2 \redrel R_1$. If the relation $R_1 \equiv_S R_2$ holds, we will say that the representation $R_1$ is
{\em subrecursively equivalent} to the representation $R_2$.

The relation 
$R_1 \redrelstrict R_2$ holds when $R_1 \redrel R_2$ and
$R_2 \not\redrel R_1$.
\qed 
\end{definition}

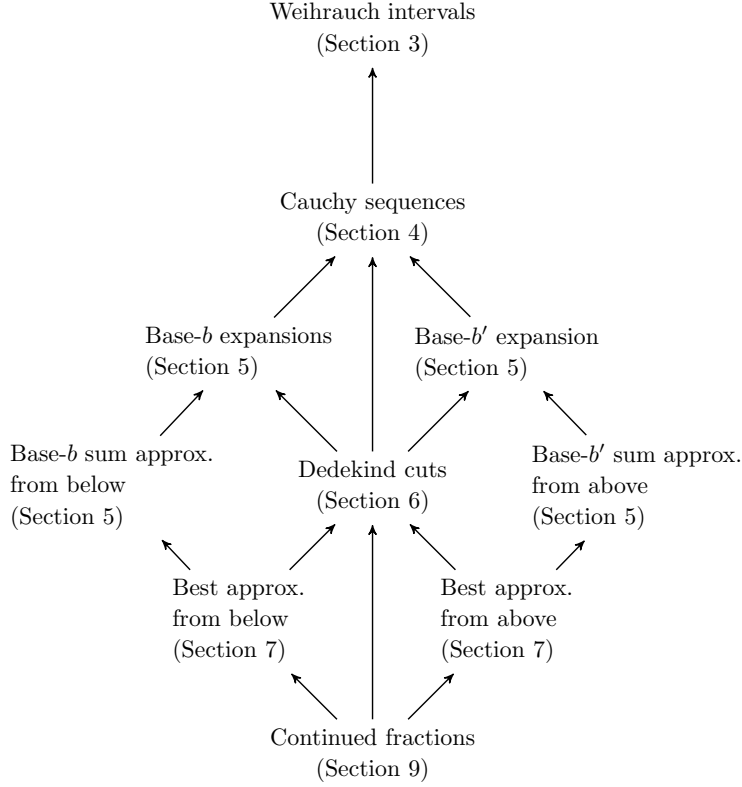
\begin{figure}
\begin{center}
\scalebox{0.9}{
\begin{tikzpicture}[>=stealth',shorten >=1pt,auto,node distance=2.8cm,
                    semithick]
  \tikzstyle{every state}=[fill=red,draw=none,text=white]

\node (Wei) [align = center] {Weihrauch intervals \\ (Section \ref{togangericannes})};

\node (Cau) [align=center, below of = Wei] {Cauchy sequences \\ (Section \ref{rosating})};

\node (Eb) [align = left, below left of = Cau] {Base-$b$ expansions\\ (Section \ref{tobiassnarteksamen})};

\node (Eb') [align=left, below right of = Cau] {Base-$b'$ expansion \\ (Section \ref{tobiassnarteksamen})};

\node (Ded) [align=center, below right of = Eb] {Dedekind cuts \\(Section \ref{skittenskje})};

\node (Sumbelow) [align= left, left = 1cm of Ded] {Base-$b$ sum approx.\\ from below \\
(Section \ref{tobiassnarteksamen})};

\node (Sumabove) [align=left, right = 1cm of Ded] {Base-$b'$ sum approx.\\ from above \\ (Section \ref{tobiassnarteksamen})};

\node (Bab) [align=left, below left of=Ded] {Best approx.\\ from below \\ (Section \ref{skitteskebest}) };

\node (Baa) [align=left, below right of=Ded] {Best approx.\\ from above \\ (Section \ref{skitteskebest})};

\node (Contfrac) [align=center, below right of = Bab]
{Continued fractions\\ (Section \ref{skittenskjecont})};

\path (Cau) edge[->] node {} (Wei);
\path (Ded) edge[->] node {} (Cau);
\path (Eb) edge[->] node {} (Cau);
\path (Eb') edge[->] node {} (Cau);
\path (Sumbelow) edge[->] node {} (Eb);
\path (Sumabove) edge[->] node {} (Eb');
\path (Contfrac) edge[->] node {} (Ded);
\path (Contfrac) edge[->] node {} (Bab);
\path (Contfrac) edge[->] node {} (Baa);
\path (Baa) edge[->] node {} (Sumabove);
\path (Bab) edge[->] node {} (Sumbelow);
\path (Bab) edge[->] node {} (Ded);
\path (Baa) edge[->] node {} (Ded);
\path (Ded) edge[->] node {} (Eb);
\path (Ded) edge[->] node {} (Eb');

\end{tikzpicture}
}
\end{center}
\caption{Overview of subrecursive degrees (equivalence classes) of representations.
\label{fig:all_encompassing}}
\end{figure}

\label{naaskaljegsnartgaaaatrene}
The equivalence relation $\equiv_S$ induces a degree structure on the 
representations. The directed graph in Figure \ref{fig:all_encompassing}
gives an overview of the relationship between some
natural degrees (equivalence classes). 
The nodes depict degrees of representations, and
each degree is labeled with one of the most well known representations in the degree.
 For two representations $R_1$ and $R_2$, there is a directed path 
 from a node labeled 
 $R_1$ to a node labeled $R_2$ if and only if $R_2 \redrelstrict R_1$. Thus, if there is a directed
 path from $R_1$ to $R_2$, we can subrecursively convert an $R_1$-representation into an $R_2$-representation,
 and if there is no directed path from $R_1$ to $R_2$, we cannot subrecursively convert an $R_1$-representation into an $R_2$-representation. Unfortunately we are not able
 to accurately  depict the complex relationship between the degrees of the base-$b$ expansions and the degrees of
 the base-$b$ sum approximations from below and above (for $b=2,3,4,\ldots)$, but our graph gives a
 rough idea of what this world looks like. See Section \ref{tobiassnarteksamen} for more
  on how these degrees relate to each other.

\subsection{Our goals and some references.}

We present a (degree) theory of
representations (of real numbers) which is based on Turing machines
and standard complexity theory.
This theory should be considered as 
a recast and an improvement of the theory developed in Kristiansen 
\cite{ciejouren} \cite{ciejourto} and Georgiev et al.~\cite{apalilf} which
is based on honest functions and subrecursive classes. The two approaches 
studying representations of reals and conversions between them are essentially
the same, even if, e.g., the reducibility relation $\redrel$ is never formally defined
in any other paper, and
it follows more or less straightforwardly from results proved in
\cite{ciejouren} \cite{ciejourto} \cite{apalilf} that the picture drawn
in Figure \ref{fig:all_encompassing} is correct.

In papers like \cite{ciejouren} \cite{ciejourto} \cite{apalilf}, and furthermore Georgiev \cite{longzero}
and Kristiansen \cite{bairelars}, the authors are just concerned with the existence or inexistence of a subrecursive conversion
from one representation to another. They never analyze the computational complexity of subrecursive
conversions, and they do not make any effort to find efficient conversions. In this paper we will
care about such matters, indeed, such matters will be our primary concern: We will impose
tight upper bounds on the running time of oracle Turing machines which convert 
one representation
into another. We will also give upper bounds on the number of oracle calls required and the
size of those calls.

In Section \ref{rosating} we study conversions between representations subrecursively equivalent
to the representation by Cauchy sequences. In Section \ref{tobiassnarteksamen} we
study conversions between representations subercursively equivalent to representations by
base-$b$ expansions and base-$b$ sum approximations. In Section \ref{skittenskje} we treat
 representations subrecursively equivalent to the representation by Dedekind cuts, 
and thereafter, in Section \ref{skitteskebest}, 
representations subrecursively equivalent to 
the representation by left/right best approximations.
Finally, in Section \ref{skittenskjecont}, we study conversions between representations subrecursively equivalent
to the representation by continued fractions. See Figure \ref{fig:all_encompassing}.

%% file: preliminaries.tex
\section{Preliminaries}
\label{prelims}

\subsection{Oracle Turing machines and complexity theory.}
We assume basic familiary with computability and computational complexity (standard textbooks are Sipser \cite{Sipser:book}, Du \& Ko \cite{DuKo:compcomp} and Arora \& Barak \cite{AroraBarak:compcomp}).

We will work with
Turing machines with oracle access to the representation being converted from. 
%All input, output and query tapes are assumed to have alphabet $\{0,1\}$ in addition to the blank symbol. All work tapes have alphabet $\{0,1,2\}$ in addition to the blank symbol.
Unless otherwise stated,  elements of $\mathbb{N}$ are assumed to be written on input, query, and output tapes in their binary representation, least-significant bit first.  
Pairs $(p,q)$ of integers
are assumed to be written using interleaved notation (i.e., the first bit of the binary representation of $p$
followed by the first bit of the binary representation of $q$, and so forth). 
Observe that the length of the representation of a pair $(p,q)$ is then $O(\log \max\{p,q\})$. Elements $p/q \in \mathbb{Q}$
are assumed to be represented by the representation of $(p,q)$.
 We denote the length of the binary representation of $x$ by $\bitlen{x}$.

Function-oracle machines are in standard use in complexity theory of functions on the set of real numbers (see, e.g., Ko \cite{kosbook}), and 
the next definition is a standard one.

\begin{definition}
A (parameterized) \emph{function-oracle Turing machine}  is a (multi-tape) Turing machine $M= (Q,q_0,F,\Sigma,\Gamma,\delta)$ with initial state $q_0 \in Q$, final states
$F \subseteq Q$, input and tape alphabets $\Sigma$ and $\Gamma$ (with $\Sigma \subseteq \Gamma$ and $\{\text{\textvisiblespace}\} \subseteq \Gamma \setminus \Sigma$), and partial transition function $\delta$ 
such that $M$ has a special \emph{query tape} and two distinct states
$q_q, q_a \in Q$ (the \emph{query} and \emph{answer} states).
%Assume a standard, computable bijection $b : (\Gamma \setminus \{\textvisiblespace\})^+ \longrightarrow \mathbb{N}$.

To be executed, $M$ is provided with a total function $f :(\Gamma\setminus\{\text{\textvisiblespace}\})^* \longrightarrow (\Gamma\setminus\{\text{\textvisiblespace}\})^*$ (the oracle) prior to execution on any input. We write $\paramorac{M}{f}$
for $M$ when $f$ has been fixed. We use $\Phi_\paramorac{M}{f}$ to denote the function
computed $\paramorac{M}{f}$. 

The transition relation of $\paramorac{M}{f}$ is defined as usual for Turing machines, except
for the query state $q_q$:
If $M$ enters state $q_q$, let $x$ be the word  currently on the query tape;
$M$ moves to state $q_a$ in a single step, and the contents
of the query tape are instantaneously changed to $f(x)$.
The query-tape head is reset to the origin, while other heads do not move.
The \emph{time- and space complexity} of a function-oracle machine is counted as for usual Turing machines, with the transition between $q_q$ and $q_a$ taking $\bitlen{f(x)}$ time steps.
%, and the space use
%of the query tape after the transition being $|f(x)|$. 
%
The \emph{input size} of a query is the number of non-blank symbols
on the query tape when $M$ enters state $q_q$.
\qed
\end{definition}

In other work on real number computation, there is a well-developed
notion of reducibility between representations that, roughly, requires
the representation to be written as an infinite string on one of the input tapes of a type-2 Turing machine \cite{DBLP:journals/tcs/KreitzW85,DBLP:journals/apal/WeihrauchK87,weirep,BRATTKA2002241}. In that setting, e.g., a function $f : \mathbb{Q} \cap [0,1] \longrightarrow \{0,\ldots,b-1\}$ is most naturally expressed by imposing a computable
ordering on its domain (e.g., rationals appear in non-decreasing order of their denominator), and the function values $f(q)$ appear encoded as bit strings in this order.
We strongly conjecture that our results carry over to the type-2 setting {\sl mutatis mutandis}.

\subsection{Some notation.}

We write
$f(n) = \poly{n}$ if $f : \mathbb{N} \longrightarrow \mathbb{N}$
is bounded above by a polynomial in $n$ with positive integer coefficients,
and $f(n) = \polylog{n}$ if $f$ is bounded above by a polynomial in $\log n$
with positive integer coefficients. 

We use the notation $f^{(n)}$ for the $n$th iterate of the function $f:\nat\to\nat$, that is,
$f^{(0)}(x) = x$ and $f^{(n+1)} = f\circ f^{(n)}$. Note the
parentheses in the superscript position, that distinguish this notation from
ordinary exponentiation. This notation is often used in conjunction with 
$\lambda$-notation, e.g. 
\begin{multline*}
    (\lambda x . g(x))^{(4)}(0) =  (\lambda x . g(x))(\lambda x . g(x))^{(3)}(0) \\
 =  g((\lambda x . g(x))^{(3)}(0)) = \dots = g(g(g(g(0))))\; .
\end{multline*}

\subsection{Farey sequences and the Stern-Brocot tree.} \label{brillelepper}

A {\em Farey sequence} is a strictly increasing sequence of fractions
between 0 and 1. The Farey sequence of {\em order} $k$, denoted $F_k$,
contains all fractions which when written in their lowest terms, have denominators 
less than or equal to $k$. Thus, e.g., $F_5$ is the sequence
$$
0/1 \; , \; 1/5 \; , \; 1/4 \; , \; 1/3 \; , \; 2/5 \; , \;  1/2 \; , \; 
3/5 \; , \; 2/3 \; , \;  3/4 \; , \; 4/5 \; , \; 
 1/1 \; .
$$
The ordered pair of two consecutive fractions in a Farey sequence is called a {\em Farey pair}.
Let $(a/b, c/d)$ be a Farey pair. The fraction $(a+c)/(b+d)$ is called the {\em mediant} of 
$a/b$ and $c/d$. 
The next theorem  was originally
proved by Cauchy \cite{Cauchy_exercises} in 1826.

\begin{theorem}\label{slankemeg}
Let $(a/b, c/d)$ be a Farey pair.
(i) We have $cb -ad = 1$ (or, equivalently $c/d- a/b = 1/(bd)$);
(ii) The mediant  $(a+c)/(b+d)$ is in its lowest terms and
lies strictly between $a/b$ and $c/d$, moreover, every other fraction lying strictly 
  between  $a/b$ and $c/d$ has denominator strictly greater than $b+d$.
\end{theorem}

E.g., $(1/3, 2/5)$ is a Farey pair as $1/3$ and $2/5$ are neighbors 
in the sequence $F_5$ (see above). The mediant of $1/3$ and $2/5$ is
$3/8$. Thus, $3/8$ lies in the open interval $(1/3, 2/5)$, and
any fraction in this open interval, with the exception of $3/8$, has
denominator strictly greater than $8$. For more on Farey pairs and Farey sequences,
see Hardy \& Wright \cite{hardy}.

We arrange the fractions strictly between 0 and 1 in a binary search 
tree $\fareypairtree$. 

\begin{definition}\label{def:farey}
The \emph{Farey pair tree} $\fareypairtree$ is the complete infinite binary tree where each node has an associated Farey pair $(a/b,c/d)$ defined by recursion on the position $\sigma \in \{0,1\}^*$ of a node in $\fareypairtree$ as follows: $\fareypairtree(\epsilon) = (0/1,1/1)$,
and if $\fareypairtree(\sigma) = (a/b,c/d)$, then $\fareypairtree(\sigma 0) = (a/b, (a+c)/(b+d))$ and $\fareypairtree(\sigma 1) = ((a+c)/(b+d),c/d)$. The \emph{depth} of a node in $\fareypairtree$ is the length of its position (with the depth of the root node being $0$).

%If $\fareypairtree(\sigma) = (a/b,c/d)$, we define %$\textrm{LEFT}(\fareypairtree(\sigma)) = a/b$
%and $\textrm{RIGHT}(\fareypairtree(\sigma)) = c/d)$. 
Abusing notation slightly, we do not distinguish
between the pair $\fareypairtree(\sigma) = (a/b,c/d)$ and the open interval $(a/b,c/d)$.

The (left) \emph{Stern-Brocot tree}\footnote{``Left'' because the Stern-Brocot tree originally concerns the interval $(0,2)$ and we are interested only in $(0,1)$ which corresponds to the left child of the Stern-Brocot tree.} $\sternbrocot$ is the infinite binary tree obtained from the Farey pair tree where each Farey pair
$(a/b,c/d)$ has been replaced by its mediant $(a+c)/(b+d)$.
\qed
\end{definition}

Thus, we have, for example
$$
\fareypairtree(0) = \left( \frac{0}{1}, \frac{1}{2} \right)  , \,  \fareypairtree(1) = \left( \frac{1}{2}, \frac{1}{1}\right)
 , \,  \fareypairtree(10) = \left( \frac{1}{2}, \frac{2}{3}\right)
 , \,  \fareypairtree(0000) = \left( \frac{0}{1}, \frac{1}{5}\right)\; .
$$
We will not use the Stern-Brocot tree directly, but we include it in the definition
for completeness.

Efficient computation of the elements of the Stern-Brocot tree
(and hence also the Farey pair tree) is possible, see Bates et al.~\cite{BATES20101020}; 
for our purposes, we simply need the next proposition.

\begin{proposition}\label{Hurwitzburdeslankesig}
There is a Turing machine
$M$ such that for any $\sigma \in \{0,1\}^*$, 
$\Phi_M(\sigma) = \fareypairtree(\sigma)$ and $M$ runs in time
$\poly{1 + \vert \sigma \vert}$.
\end{proposition}

We round off this section by  stating and proving a few  properties of
the Fairy pair tree.

\begin{proposition}\label{prop:viskalbeggeslankeos} If
$(a/b,c/d)$ is a Farey pair at depth $h$ in $\fareypairtree$, then
$a+b+c+d\geq h+3$.
\end{proposition}

\begin{proof}
For $h=0$, we have $0+1+1+1=3=h+3$. 

Let $h>0$. Assume the proposition for $h-1$ and
let $(a/b,c/d)$ be an arbitrary pair at depth $h-1$.
Then
 a pair at level $h$ is of the form (i) $(a/b, (a+c)/(b+d))$  or of the 
 form (ii) $((a+c)/(b+d), c/d)$. In case (i), we have
 \begin{multline*}
     a+(a+c)+b+(b+d) \; = \;  (a+b+c+d)+(a+b) \; \ge\;  (a+b+c+d)+1 \\ \; \ge \;
     (h+2)+1 = h+3\;.
 \end{multline*}
A symmetric argument will show that the proposition also holds in case (ii).
\end{proof}

\begin{proposition}\label{prop:alleskalslankesig}
Let $p/q \in \mathbb{Q} \cap [0,1]$ be a fraction in its lowest terms. Then,
$p/q$ is a fraction in a Farey pair at depth at most $p+q-2$ in $\fareypairtree$.
\end{proposition}

\begin{proof}
By construction, for any depth $n \geq 0$, the set of intervals $[a/b,c/d]$ occurring in $\fareypairtree$ at depth $n$ cover the unit interval, and each pair of intervals have
at most one point in common (which must be an end point). Hence, $p/q$ occurs
in some interval $[a/b,c/d]$ at any depth $n$, and by Theorem \ref{slankemeg}
we have $cb - ad = 1$. Assume for the sake of a contradiction that $n > p+q-2$ and 
$p/q \notin \{a/b,c/d\}$.
Since $p/q \notin \{a/b,c/d\}$, we have $a/b < p/q < c/d$, and
thus also $1 \leq pb - qa$ and $1 \leq qc - pd$.
Hence
\begin{multline*}
a + b + c + d \leq (a + b)(qc - pd) + (c + d)(pb - qa) \\
= p(cb- ad) + q(cb - ad)
= p + q\; .
\end{multline*}
This contradics  Proposition \ref{prop:viskalbeggeslankeos} 
which implies that $a + b + c + d \geq p+q+1$. Hence, $p/q$
must occur as an endpoint, and the first level at which $p/q$ appears as an endpoint must be at most $p+q-2$.
\end{proof}

\begin{lemma}\label{fareysizelemma}
Let $I = (a_n/b_n,c_n/d_n)$ be a Farey interval at depth $n$ in $\fareypairtree$,
and for $i=0,\dots,n$, let $(a_i/b_i,c_i/d_i)$,  denote the Farey pairs 
along the path from the root to $I$.  
Then the numbers $a_i,b_i,c_i,d_i$ are all bounded by $(n+1)(a_n+b_n)$.
\end{lemma}

\begin{proof}
We claim that for all $i$, $a_i+b_i+c_i+d_i \le (i+1)(a_n+b_n)$.
For $i=0$, we have $a_0+b_0+c_0+d_0 = 1\cdot(a_0+b_0)$.  Assume this holds for arbitary $i$. If the next node is a left
child, then
\begin{multline*}
 a_{i+1}+b_{i+1}+c_{i+1}+d_{i+1} = a_i+b_i+(a_i+c_i)+(b_i+d_i) \\
 \le a_i+b_i+(i+1)(a_n+b_n) 
 = (i+2)(a_{n}+b_{n})\; . 
\end{multline*}
If the next node is a right child, then
\begin{multline*}
 a_{i+1}+b_{i+1}+c_{i+1}+d_{i+1} = (a_i+c_i)+(b_i+d_i)+c_i+d_i \\
 \le (i+1)(a_n+b_n)+c_i+d_i 
\le  (i+1)(a_n+b_n) + a_{i+1}+b_{i+1} \\
\le (i+2)(a_n+b_n) \; .
\end{multline*}
\end{proof}

%% file: weihrauch.tex
\section{Weihrauch Intersections}                     
                     
\label{togangericannes}

\begin{definition}
A function $I: \nat \longrightarrow \rational \times \rational$ is a {\em Weihrauch intersection}
for the real number $\alpha$ if the left component of the pair $I(i)$ is strictly less that the right component of the pair $I(i)$ (for all $i\in \nat$) and
$$
\{ \ \alpha \ \} \; = \; \bigcap_{i=0}^{\infty} I_i^O
$$
where $I_i^O$ denotes the open interval given by the the pair $I(i)$. 
\end{definition}

\begin{theorem} \label{wi}
Any computable real number can be represented by a polynomial-time computable Weihrauch intersection.
\end{theorem}

\begin{proof}
A computable real number $\alpha$ has a computable Cauchy sequence  $C:\nat\longrightarrow \rational$
with the property $\vert C(n)- \alpha \vert < 2^{-n}$.
Let $M$ be a Turing machine computing $C$. 
We can w.l.o.g. assume that $\alpha\in (0,1)$. 

Compute $I(k)$ by the following algorithm: Find the greatest $n$ such that $n\leq k$
and $C(n)$ can be computed by $M$ in $k$ steps. Let $I(k) = (C(n)-2^{-n},C(n)+2^{-n})$.
Let $I(k) =(0,1)$ if no such $n$ exists
(it is possible to arrange this such that we have $I^O_{k+1} \subseteq I^O_{k}$).
\end{proof}

The  representation by Weihrauch intersections
is one of the main representations
in Weihrauch's seminal book \cite{wei}, and it is special among  the 
representations we consider in this paper: There exists  a time-bound $t$
such that every computable real has a Weihrauch
intersection computable by a Turing machine running in time $O(t)$ (by Theorem \ref{wi},
this will for sure be true for any $t$ that dominates all polynomials).
For every other representation $R$ considered in this paper,
there will for any time-bound $t$ exists a time-bound $s$ such that $O(t)_R$
is strictly included in $O(s)_R$. The degree of the representation by Weihrauch
intersections will be the zero degree of the degree structure described in Section \ref{kortstokk}.

Representation of reals by Weihrauch intersections are
also known as representation by {\em nested intervals}. In order to simplify our definition,
we have not required the intervals to be nested, but any
Weihrauch intersection can be easily  converted to a nested one.
A number of  subrecursively equivalent  representations can be found  in \cite{wei}, but 
they are all pretty similar from
our point of view, and we will not discuss any of them.

%% file: cauchy.tex
\section{Representations subrecursively equivalent to Cauchy sequences}
\label{rosating}

\subsection{Cauchy sequences.}
\begin{definition}
Let $\alpha \in (0,1)$ be an irrational number. 
Then $C:\nat^+ \longrightarrow \rational$ is a {\em Cauchy sequence} for $\alpha$ if
$\left\vert \alpha - C(n)\right\vert < {n}^{-1}$.
\qed
\end{definition}

\begin{lemma} \label{nikkedukke}
Let $C : \mathbb{N}^+ \longrightarrow \mathbb{Q}$ be a Cauchy sequence for an
irrational number $\alpha \in (0,1)$. Let $p/q = C(n)$.
There is a parameterized function-oracle Turing machine $M$
such that 
\begin{itemize}
    \item $\Phi_M^C :\mathbb{N} \longrightarrow \mathbb{Q} \times \mathbb{Q}$ is a Weihrauch intersection for
$\alpha$
\item  $M^C$ on input $n$ runs in time
$\polylog{ \max\{p,q\},n}$ and  uses exactly one oracle call of input size
$O(\log n)$.
\end{itemize}
\end{lemma}

\begin{proof}
Let 
$$
I(n)  =  (C(n)- n^{-1},C(n)+ n^{-1})\; .
$$
Then, $I$ is a Weihrauch intersection for $\alpha$ if $C$ is a Cauchy sequence for $\alpha$. Hence, only one oracle call to $C$ is needed,
and $I(n)$ can be obtained by basic arithmetic operations on (the binary representations of) $p,q$ and $n$. 
\end{proof}

Lemma \ref{nikkedukke} shows that a Cauchy sequence can be subrecursively converted into a Weihrauch intersection. {\sl Will it be possible to subrecursively convert a Weihrauch intersection into a Cauchy sequence?} In order to give a negative answer that question, we need a 
presumable  rather well known theorem.

\begin{theorem}\label{sologsommer}
For any time-bound  $t$
there exists a computable irrational number $\alpha$ such that  no Cauchy sequence for $\alpha$ can by computed by a Turing machine running in time $O(t)$.
\end{theorem}

It is not hard to see that the theorem holds. Let $t$ be a fairly fast-increasing time bound,
and let $A$ be any set of natural numbers
such that membership in $A$ can be decided by a Turing machine, but not by an $O(t)$-time Turing machine (the existence of such  a set can be shown by a standard diagonalization argument). 
Consider   the irrational number
$\alpha$  given by the base-$2$ expansion
$0.a_1 a_2 a_3 \dots$ where the two digits $a_{2i-1} a_{2i}$ are $11$ if $i\in A$; and $01$ otherwise. Now, $\alpha$ will obviously have a computable Cauchy sequence, but no 
$O(t)$-time Turing machine can compute such a Cauchy sequence. 
If a Cauchy sequence for $\alpha$ can be computed in time $O(t)$, 
then a  Turing machine
$M$ can decide if $m$ is in the set $A$ in time $O(t)$:
First $M$ computes $C(2^{2m})$. By assumption this can be done in time $O(t(\bitlen{2^{2m}})=O(t(2m+1))$.
Thereafter, $M$  determines the digits $0.\xd_1\xd_2 \ldots $ of the base-2 expansion of
$C(2^{2m})$. Observe that $C(2^{2m})$ lies sufficiently close to $\alpha$ to ensure that 
the digit $\xd_{2m-1}$ coincide with digit $a_{2m-1}$, and thus,  $m\in A$ iff $\xd_{2m-1}=1$.
Hence, $M$  can decide if $m$ is a member of
$A$ by computing  $\xd_{2m-1}$, and this can obviously be done in time $O(t)$.
Since no Turing machine can decide membership in $A$ in time $O(t)$, we can conclude that the theorem holds.

Theorem \ref{wi} states that any computable  real can be represented by
polynomial-time computable Weihrauch intersection. Thus, any real
can be represented by a  Weihrauch intersection computable in, let us say, time $O(2^n)$.
Now, $2^n$ is a fixed time-bound,
 and if it were possible 
to subrecursively convert a Weihrauch intersection into a Cauchy sequence,
the any computable irrational would be represented by a Cauchy sequence computable in
time $O(t)$ for some fixed time-bound $t$. By Theorem \ref{sologsommer}, such a $t$ does not exists, and we can conclude that Weihrauch intersections
cannot be subrecursively converted into Cauchy sequences.

\subsection{Definitions.}
The next definition gives some representations which are subrecursively equivalent 
to the representation by Cauchy sequences.

\begin{definition} \label{kmoCauchy}
Let $\alpha \in (0,1)$ be an irrational number. 
\begin{enumerate}
\item $C:\nat^+ \longrightarrow \rational$ is a {\em strictly increasing Cauchy sequence} for $\alpha$ if
(i) $C$ is a  Cauchy sequence for $\alpha$ and (ii) $C(n)< C(n+1)$.
\item Let $b\geq 2$ be a natural number. Then, $A:\nat^+ \longrightarrow \integer$ is a 
{\em converging base-$b$ sequence} for $\alpha$ if
$A(n)b^{-n}$ is a Cauchy sequence for $\alpha$.
\item  $D: \integer\times \nat^+ \longrightarrow \{0,1\}$ is a
{\em fuzzy (Dedekind) cut} for  $\alpha$ if
$$
D( p, q ) = 0\; \Rightarrow \;  \alpha < \frac{p+1}{q} \;\;\;\;\;\;\; \mbox{ and }  \;\;\;\;\;\;\;
D( p,q ) = 1\; \Rightarrow \;   \frac{p-1}{q} < \alpha \; .
$$
\item $S:\nat^+ \longrightarrow \{-1,0,1\}$ is a {\em signed  digit expansion} for $\alpha$ if
$$\alpha = \sum_{i=1}^{\infty} S(i)2^{-i}\; . $$
\end{enumerate}
\qed
\end{definition}

Converging base-2 sequences are used in Friedman and Ko \cite{KO1982323}
and also in the monograph Ko \cite{kosbook}.
Signed digit expansions  also seem to be well known.
The representation is discussed in Weihrauch's book \cite{wei} and
appears in several rather recent papers, eg.~Berger et al.~\cite{berger} and 
Bauer et al.~\cite{bauer}.
 The representations by strictly increasing
Cauchy sequences and fuzzy Dedekind cuts are discussed for the first time
in this paper.

\subsection{Cauchy sequences to fuzzy cuts.}

Let $C$ be a Cauchy sequence for $\alpha$.
We define the map $D: \integer\times \nat^+ \longrightarrow \{0,1\}$ by
$$D(p,q) \;\; = \; \; \begin{cases}
0   & \mbox{if $C(q)\leq pq^{-1}$}   \\
1   & \mbox{if $C(q)> pq^{-1}$}
\end{cases} 
$$

\begin{lemma}\label{lem:easy_fuzzy}
$D$ is a fuzzy cut for $\alpha$.
\end{lemma}

\begin{proof}
First we prove
\begin{align} \label{sandaker}
D(  p, q  ) = 0\; \Rightarrow \;  \alpha < \frac{p+1}{q}\; . 
\end{align}
Assume $D ( \ p, q \ ) = 0$.
If $\alpha < C(q)$, then we obviously have 
$$
\alpha \; < \;  C(q) \; \leq \;  \frac{p}{q} \;  < \;  \frac{p+1}{q}\; .
$$
Thus, (\ref{sandaker}) holds if $\alpha < C(q)$. Now, assume $\alpha > C(q)$.  By the definition  of a Cauchy sequence, we have
$\alpha - C(q) < q^{-1}$. Hence $\alpha   < C(q) + q^{-1}$. Furthermore, by the definition of $D$, we have
\[  \alpha   \; < \; C(q) + \frac{1}{q} \; 
\leq \; \frac{p}{q} + \frac{1}{q} \; = \; \frac{p+1}{q} \;  . \]
This proves (\ref{sandaker}). We also need to prove
\begin{align} \label{xsandakerveien}
D( p,q ) = 1\; \Rightarrow \;   \frac{p-1}{q} < \alpha \; . 
\end{align}
The proof of (\ref{xsandakerveien}) is symmetric to the proof (\ref{sandaker}). The lemma follows from 
(\ref{sandaker}) and (\ref{xsandakerveien}).
\end{proof}

\begin{lemma}
Let $C :\nat^+ \longrightarrow \rational$ be a Cauchy sequence for an irrational number $\alpha \in (0,1)$.
There is a parameterized function-oracle Turing Machine $M$ such that 
\begin{itemize}
    \item $\Phi^C_M : \integer\times \nat^+ \longrightarrow \{0,1\}$ is a fuzzy Dedekind cut for $\alpha$
  \item $M^C$ on input $(p,q)$ runs in time $\poly{\max{\bitlen p,\bitlen q,\bitlen{C(q)}}})$ 
  and uses a single oracle call of input size at most $O(\log q)$.
\end{itemize}
\end{lemma}

\begin{proof}
A single call to $C$ on input (the binary representation of) $q$ yields (the binary representation of)
$C(q)$, and by Lemma \ref{lem:easy_fuzzy}, a single comparison of $C(q)$ to $pq^{-1}$ yields
$D(p,q)$. A binary representation of the rational number $pq^{-1}$ can be computed in polynomial time
in the size of the representations of $p$ and $q$ (that is, in time $\polylog{\max\{p,q\}}$), 
and the final comparison of $C(q)$ and $pq^{-1}$ can be performed in time polynomial in $\max\{\bitlen{C(q)}, \bitlen{pq^{-1}}\}$.
\end{proof}

\subsection{Fuzzy cuts to signed digit expansions.}

Let $\alpha\in(0,1)$, and let $D$ be a fuzzy cut for $\alpha$.
We will use $D$ to define a signed digit expansion $S$.

For any map $S:\nat^+ \longrightarrow \{-1,0,1\}$, let 
$s_n = \sum_{i=1}^n S(i)2^{-i}$. Furthermore, 
let $a_n$ be the unique integer such that $a_n2^{-n}= s_n$, and let $M(r_1,r_2)$ denote the midpoint between the two
rationals $r_1$ and $r_2$, that is, $M(r_1,r_2)= r_1 + (r_2-r_1)/2$. Observe that
$$
M\left(s_n-\frac{1}{2^{n+1}} \, , \, s_n \right) = \frac{4a_n-1}{2^{n+2}} 
\;\;\; \mbox{ and } \;\;\; M\left(s_n \, , \, s_n+\frac{1}{2^{n+1}} \right) = \frac{4a_n+1}{2^{n+2}}
$$
Moreover, observe that $s_n= 2a_n/2^{n+1}$. 

We define the signed digit expansion $S$ by $S(1)=1$ and
$$S(n+1) \;\; = \; \; \begin{cases}
-1   & \mbox{if $D(2a_n,2^{n+1})= 0$ and  $D(4a_n-1,2^{n+2})= 0$ (Case 1)}   \\
0   & \mbox{if  $D(2a_n,2^{n+1})= 0$ and  $D(4a_n-1,2^{n+2})= 1$  (Case 2)} \\
1   & \mbox{if  ($D(2a_n,2^{n+1})= 1$ and  $D(4a_n+1,2^{n+2})= 1$ (Case 3)}  \\
0   & \mbox{if  $D(2a_n,2^{n+1})= 1$ and  $D(4a_n+1,2^{n+2})= 0$ (Case 4)} 
\end{cases} 
$$
The first conjunct in (Case 1) is not needed as it follows from the second conjunct. The same goes for the first conjunct in (Case 3).
The superfluous conjuncts  are included in order to make it easy to see that the four cases are mutually exclusive and exhaustive.

\begin{lemma}\label{lem:fuzzy_to_signed}
$$
\alpha \; = \; \sum_{i=1}^\infty S(i)2^{-i}\; . 
$$
\end{lemma}

\begin{proof}
We will prove
\begin{align}
s_n - 2^{-n} < \alpha < s_n + 2^{-n}  \label{amandabor}
\end{align}
by induction on $n$. The lemma follows straightforwardly from (\ref{amandabor}).
It is obvious that (\ref{amandabor}) holds when $n=1$  (as we have assumed $\alpha\in (0,1)$).

Assume by induction hypothesis that (\ref{amandabor}) holds. We need to prove that
\begin{align} \label{xamandabor}
\alpha \in (  s_{n+1} - 2^{-(n+1)} \, , \, s_{n+1} + 2^{-(n+1)}  )  
\end{align}

{\em (Case 1.)} In this case we have $s_{n+1}= s_n - 2^{-(n+1)}$.
Thus, in order to prove (\ref{xamandabor}), we need to prove $\alpha\in(s_n - 2^{-n} ,s_n)$.
We have $s_n - 2^{-n}< \alpha$ by the induction hypothesis (\ref{amandabor}). Moreover, since $D(4a_n-1,2^{n+2})= 0$,
we have 
$$
\alpha \; < \; \frac{4a_n-1+1}{2^{n+2}} \; = \; \frac{a_n}{2^n} \; = \; s_n\; .
$$
This proves that (\ref{xamandabor}) holds in (Case 1).

{\em (Case 2.)} In this case we have $s_{n+1}= s_n$.
Thus, in order to prove (\ref{xamandabor}), 
we need to prove $\alpha\in(s_n - 2^{-(n+1)} ,s_n +  2^{-(n+1)})$.
Since $D(2a_n,2^{n+1})= 0$, we have 
$$
\alpha \; < \; \frac{2a_n+1}{2^{n+1}}  \; = \; \frac{a_n}{2^n} + \frac{1}{2^{n+1}}\; = \; s_n  + \frac{1}{2^{n+1}}\; .
$$
Since $D(4a_n-1,2^{n+2})= 1$,
we have 
$$
\alpha \; > \; \frac{4a_n-1-1}{2^{n+2}} \; = \; \frac{a_n}{2^n} - 2^{-(n+1)} \; = \; s_n - 2^{-(n+1)}\; .
$$
This proves that (\ref{xamandabor}) holds in (Case 2).

(Case 3) is symmetric to (Case 1), and (Case 4) is symmetric to (Case 2).
\end{proof}

\begin{lemma}
Let
 $D : \integer\times \nat^+ \longrightarrow \{0,1\}$ be a fuzzy cut for
 of an irrational number $\alpha \in (0,1)$.
There is a parameterized function-oracle Turing Machine $M$ such that 
\begin{itemize}
    \item $\Phi^D_M : \nat^+ \longrightarrow \{-1,0,1\}$ is a signed digit expansion of $\alpha$
    \item $M^D$ on input $n$ runs in time $\poly{n}$ and uses $3(n-1)$ oracle calls of input size at most $O(n)$.
\end{itemize}
\end{lemma}

\begin{proof}
By Lemma \ref{lem:fuzzy_to_signed}, $S(1) = 1$, and computing $S(n)$ for $n > 1$
can be done by computing the integer $a_{n-1}$
satisfying $a_{n-1} s ^{-n-1} = s_n = \sum_{i=1}^{n-1} S(i)2^{-i}$,
and subsequently performing oracle calls returning
the values of
\begin{equation} \label{abcde}
      D(2a_{n-1},2^{n}) \; , \;  D(4a_{n-1}-1,2^{n+1}) \; \mbox{ and }
\; D(4a_{n-1}+1,2^{n+1})\; . 
\end{equation} 
Observe that the binary representations of the rational numbers $2^{-1}, \ldots,2^{-(n-1)}$ have length $O(n)$, and that
each representation can be computed in time $\poly{n}$. Hence,
if $S(1),\ldots,S(n-1)$ are known,
then $s_{n-1}$ can be computed in time
$\poly{n}$ using standard arithmetical operations. Furthermore, observe that the binary representation of $s_{n-1}$ (and thus $a_{n-1}$) has length $O(n)$. This implies that (i) the size of the oracle calls in (\ref{abcde}) 
is at most $O(n)$, and that (ii) $a_{n-1}$, $2 a_{n-1}$, and $4 a_{n-1}$ can  be computed
in time $\poly{n}$. 

Using the obvious recursive algorithm for $S(n)$ requires computing $S(1),\ldots$, $S(n-1)$, 
hence time $O(n \poly{n}) = \poly{n}$, and a total of $3(n-1)$ oracle calls,
each of size $O(n)$.
\end{proof}

\subsection{Signed digit expansions to Cauchy sequences.}

Let $S$ be a signed digit expansion of $\alpha\in(0,1)$.
Then we have
$\left\vert \alpha - \sum_{i=1}^n S(i)2^{-i} \right\vert \; < \; 2^{-n}$.
Let
$$C(n) \; = \;  \sum_{i=1}^{\lceil \log_2 n \rceil} S(i)2^{-i}$$
and we  have
$$\left\vert \alpha - C(n) \right\vert \; = \;
\left\vert \alpha - \sum_{i=1}^{\lceil \log_2 n \rceil} S(i)2^{-i} \right\vert
 \; < \; 2^{-\lceil \log_2 n \rceil} \; \leq \;  n^{-1}\; .$$
Hence, $C$ is a Cauchy sequence for $\alpha$.

\begin{lemma}
Let $S:\nat^+ \longrightarrow \{-1,0,1\}$ 
be a signed digit expansion for an irrational number $\alpha \in (0,1)$.
There is a parameterized function-oracle Turing Machine $M$ such that
\begin{itemize}
    \item $\Phi^S_M :\nat^+ \longrightarrow \rational$ is a Cauchy sequence for $\alpha$
    \item $M^S$ on input $n$ runs in time $\polylog{n}$ and uses $\lceil \log_2 n \rceil$ oracle calls of input size at most $O(\log n)$.
\end{itemize}
\end{lemma}

\begin{proof}
The result follows almost immediately from the text just prior
to the lemma. Observe that computation
of $C(n) = \sum_{i=1}^{\lceil \log_2 n \rceil} S(i)2^{-i}$ can be performed using $\poly{\log n} = \polylog{n}$ operations
on rationals whose representation has length at most
$O(\log n)$, hence in total time $\polylog{n}$.
\end{proof}

\subsection{From Cauchy sequences to strictly increasing Cauchy sequences.}

Let $C$ be a Cauchy sequence for some real number $\alpha$. 
Thus, for all $n$, we have
\begin{align}
 \alpha \in ( \ C(2^{n})-2^{-n} \, , \, C(2^{n})+2^{-n} \ ) \label{pencil}
\end{align}
We will use $C(n)_\ell$ to denote the left endpoint of the interval in (\ref{pencil}),
that is, $C(n)_\ell= C(2^{-n})-2^{-n}$, and
we define $\widehat{C}$ by 
$$\widehat{C}(n) \; = \; C(n+2)_\ell - 2^{-(n+1)}\; = \; C(2^{n+2}) - 2^{-(n+2)} - 2^{-(n+1)} \;.$$

\begin{lemma}  \label{junityvetyve} 
If $C$ is a Cauchy sequence for $\alpha$, then $\widehat{C}$ is a 
strictly increasing Cauchy sequence for $\alpha$.
\end{lemma}

\begin{proof}
Let $C$ be a Cauchy sequence for $\alpha$.
We will we prove
\begin{align} \label{kaniner}
    \left\vert \alpha - \widehat{C}(n)  \right\vert \; < \; 2^{-n}
\end{align}
and
\begin{align} \label{hjort}
  \widehat{C}(n)\; < \; \widehat{C}(n+1) \; .
\end{align}
and thus the lemma holds.

First we observe that
$\left\vert \alpha - C(n+2)_\ell \right\vert  <  2^{-(n+1)}$,
and thus, we have
\begin{multline*}
\left\vert \, \alpha - \widehat{C}(n) \, \right\vert \; = \; 
\left\vert \, \alpha - ( \ C(n+2)_\ell - 2^{-(n+1)} \ ) \, \right\vert \; = \; \\
\left\vert \, \alpha - C(n+2)_\ell  \, \right\vert + 2^{-(n+1)} \; < \; 
 2^{-(n+1)} +  2^{-(n+1)} \; = \;  2^{-n}\; .
\end{multline*}
The first equality holds by the definition of $\widehat{C}$, and 
the second equality holds since $C(n+2)_\ell$ lies below $\alpha$.
This proves that (\ref{kaniner}) holds.

Next we observe that 
\begin{align} \label{pinnsvin}
\left\vert C(2^{n}) - C(2^{n+1}) \right\vert \; < \; 2^{-n} +  2^{-(n+1)}\; .
\end{align}
holds for all $n$. Now, assume for the sake of a contradiction that (\ref{hjort}) does not hold, that is,
assume there exists $m$ such that
$\widehat{C}(m) \geq   \widehat{C}(m+1)$. 
Then we have
$$
\widehat{C}(m)- \widehat{C}(m+1) \geq 0\; .
$$
By the definition of $\widehat{C}$, we have 
$$
 C(m+2)_\ell - 2^{-(m+1)} - ( \ C(m+3)_\ell - 2^{-(m+2)} \ ) \geq 0\; .
$$
Hence
$$
 C(m+2)_\ell -  C(m+3)_\ell  \geq 2^{-(m+2)}\; .
$$
By the definition of  $C(\cdot)_\ell$, we have
$$
 C(2^{m+2}) - 2^{-(m+2)} - ( \  C(2^{m+3}) - 2^{-(m+3)} \ ) \geq 2^{-(m+2)}\; .
$$
Hence,  we have
$$
 C(2^{m+2})  - C(2^{m+3}) \geq 2^{-(m+2)} + 2^{-(m+2)} - 2^{-(m+3)} = 2^{-(m+2)} - 2^{-(m+3)} \; .
$$
This contradicts (\ref{pinnsvin}), and we have proved that (\ref{hjort}) holds.
\end{proof}

\begin{lemma}
Let
 $C : \nat^+ \longrightarrow \rational$ be a Cauchy sequence for an
 irrational $\alpha \in (0,1)$.
There is a parameterized function-oracle Turing Machine $M$ such that 
\begin{itemize}
    \item $\Phi^C_M :\nat^+ \longrightarrow \rational$ is a strictly increasing Cauchy sequence for $\alpha$
    \item $M^C$  on input $n$  runs in time $\poly{\max\{\bitlen{C(2^{n+2})},n\}}$ and uses a single oracle call of input size  $O(n)$.
\end{itemize}
\end{lemma}

\begin{proof}
Om input $n$, $M$ constructs the number $2^{n+2}$ (representable in $O(n)$ bits) and performs the oracle call $C(2^{n+2})$,
and then calculates and outputs the rational number $C(2^{n+2}) - 2^{-(n+2)} - 2^{-(n+1)} = \widehat{C}(n)$. This calculation
involves basic arithmetic on numbers representable 
in $\max\{\log \bitlen{C(2^{n+2})},O(n)\}$ bits,
hence in total time $\poly{\max\{\bitlen{C(2^{n+2})},n\}}$. By Lemma 
\ref{junityvetyve}, $\widehat{C}$ is a strictly increasing Cauchy sequence for $\alpha$.
\end{proof}

\subsection{From converging base-$b$ sequences to Cauchy sequences.}

 Let $A$ be a converging base-$b$ sequences for $\alpha$, and 
let $C(n)=A(n)b^{-n}$. Now, by our definitions, $C$ is a Cauchy sequence for $\alpha$. Thus the proof of the next lemma is straightforward.

\begin{lemma} Let $b\geq 2$, and let
$A:\nat^+ \longrightarrow \integer$ be a converging base-$b$ sequence for 
an irrational $\alpha \in (0,1)$.
There is a parameterized function-oracle Turing Machine $M$ such that
\begin{itemize}
    \item $\Phi^A_M :\nat^+ \longrightarrow \rational$ is a Cauchy sequence for $\alpha$
    \item $M^A$ on input $n$ runs in time $\poly{\max\{\log{A(n)},n\}}$ and uses a single oracle call of input size  $O(\log n)$.
\end{itemize}
\end{lemma}

\begin{proof}
On input $n \in \nat^+$, $M$ performs the oracle call $A(n)$
(where $n$ is representable in $O(\log n)$ bits), and then
computes the rational number $A(n)b^{-n}$. It is well-known that $b^{-n}$ can be
computed using $O(\log n)$ multiplications (where $b^{-1}$ is assumed to be hard-coded).
All the numbers involved occupy 
$O(\max\{\log A(n),n\})$ bits. Hence, $M$ runs in time
at most $\poly{\max\{\log{A(n)}, n\}}$.
\end{proof}

\subsection{From  Cauchy sequences to  converging base-$b$ sequences.}

Let $C$ be a Cauchy sequence for $\alpha$. We show how to compute $A(n)$ where $A$ is a 
converging base-$b$ sequence for $\alpha$.

First we define  the sequences $X_0,X_1,X_2,\ldots$ and $Y_0,Y_1,Y_2,\ldots$.
Let $p\in\integer$ and $q\in\nat^{+}$ be arbitrary. Let $X_0 = \ddiv{p}{q}$, let $Y_0 = \mmod{p}{q}$,
and let
$$  X_{n+1} = X_{n}\times b + (\ddiv{(Y_{n}\times b)}{q}) \;\;\; \mbox{ and } \;\;\;
Y_{n+1} = \mmod{(Y_{n}\times b)}{q}$$
where the operators
\begin{itemize}
    \item $\ddiv{x}{y}$ (integer division)
    \item $\mmod{x}{y}$ (the remainder of integer division)
\end{itemize}
have the property
\begin{align} \label{ikkesinoe}
    ( \ddiv{x}{y} )\times y + (\mmod{x}{y}) = x \; .
\end{align}

\begin{lemma} \label{hjemmeforste} For any $n$, we have
(i) $Y_n < q$, (ii) $p q^{-1}= X_n b^{-n} + Y_n b^{-n}q^{-1}$ 
and (iii) $X_n b^{-n}\leq pq^{-1} < (X_n +1) b^{-n}$.
\end{lemma}

\begin{proof}
It is obvious that (i) holds, and (iii) follows straightforwardly from (i) and (ii). We prove (ii)
by induction on $n$. It is obvious that (ii) holds if $n=0$. Furthermore, we have
\begin{multline*}
X_{n+1} b^{-(n+1)} + Y_{n+1} b^{-(n+1)}q^{-1} = \\
(X_{n} b + (\ddiv{(Y_{n} b)}{q}))  b^{-(n+1)} + (\mmod{(Y_{n} b)}{q}) b^{-(n+1)}q^{-1} = \\
[X_{n} bq + (\ddiv{(Y_{n} b)}{q})q + (\mmod{(Y_{n} b)}{q})] b^{-(n+1)}q^{-1} = \\
[X_{n} bq + Y_{n} b)] b^{-(n+1)}q^{-1} =
X_{n} b^{-n} + Y_{n} b^{-n}q^{-1} = p q^{-1} 
\end{multline*}
where the first equality holds by the definition of $X_{n+1}$ and $Y_{n+1}$, the third equality holds
by (\ref{ikkesinoe}) and the last equality holds by the induction hypothesis.
\end{proof}

We may now compute $A(n)$ by the following procedure:
\begin{itemize}
    \item let $p/q= C(2n)$ where $p\in\integer$ and $q\in\nat^{+}$
    \item let $k= \lceil \log_b 2n \rceil$
    \item compute $X_k$
    \item output $A(n) := X_k$.
\end{itemize}

By Lemma \ref{hjemmeforste} (iii), we have
$$
 A(n)b^{-k} \leq C(2n) < (A(n)+1)b^{-k} 
$$
and thus
$$
 \left\vert A(n)b^{-k} - C(2n) \right\vert < b^{-k} = b^{-\lceil \log_b 2n \rceil} \leq (2n)^{-1} \; .  
$$
Moreover, as $C$ is a Cauchy sequence for $\alpha$, we have 
$\left\vert \alpha - C(2n) \right\vert < (2n)^{-1}$, and
thus, we also have 
$\left\vert A(n)b^{-k} - \alpha \right\vert < n^{-1}$.
This proves that $A$ is a converging base-$b$ sequence for $\alpha$.

\begin{lemma}
Let $C : \nat^+ \longrightarrow \rational$ be a Cauchy sequence for an irrational $\alpha \in (0,1)$, and let $b \geq 2$. 
There is a parameterized function-oracle Turing Machine $M$ such that 
\begin{itemize}
    \item $\Phi^C_M :\nat^+ \longrightarrow \integer$ is converging base-$b$ expansion 
    for $\alpha$
    \item $M^C$ on input $n$ runs in time $\poly{(\log n) + \bitlen{C(2n)}}$ and
    uses
a single oracle call of input size  $O(\log n)$.
\end{itemize}
\end{lemma}

\begin{proof}
The algorithm given above  uses $C(2n)$, and hence an oracle call of input size $O(\log n)$.
Moreover the algorithm uses time at most $\polylog n$ to compute
$k = \lceil \log_b 2n \rceil$. The final for-loop consists of $k$ iterations
involving 5 arithmetical operations in each iteration; it is a straightforward
induction to see that each of these operations is applied to non-negative integer arguments of size at most 
$$\bitlen{C(2n)}+k(\bitlen{b}+1) \; = \; O((\log n) + (\log \bitlen{C(2n)})$$ 
bits. As each arithmetical operation is computable
in polynomial time in the size of the representation, the total
time use is $$O(k \, \poly{(\log n) + \bitlen{C(2n)}}) \; = 
\; \poly{(\log n) + \bitlen{C(2n)}}\; .$$
\end{proof}

\subsection{Summary.}

\label{oppsummeringcauchy}

Recall  that $O(t)_R$ denotes class of all irrational $\alpha$
in the interval $(0,1)$ such that at least one $R$-representation of $\alpha$ 
is computable by a Turing machine running in time 
$O(t(n))$ ($n$ is the length of the input, see Definition \ref{snartsommerigjen} at page \pageref{snartsommerigjen}).
When we combine the results on the complexity of conversions among representations 
in this section, we get
 the following theorem.

\begin{theorem}
Consider the representations by (1)  Cauchy sequences,
(2)  increasing Cauchy sequences, (3)  fuzzy cuts, (4)  signed-digit expansions
and (5) converging base-$b$ sequences, and let $R_1$ and $R_2$ be
any two of these five representations.
For an arbitrary time-bound $t$,  we have 
$$O(t)_{R_2} \; \subseteq \; O(\poly{t(2^{n})})_{R_1}\; .$$
\end{theorem}

%% file: sumeapproxbaseexp.tex
\section{Base-$b$ expansions and sum approximations}

\label{tobiassnarteksamen}

\subsection{The base-$b$ expansions.}

The representation of reals by base-$b$ expansions, or perhaps we should say base-10 expansions, is
very well known. We are talking about
the standard daily-life representation of reals.
  We will restrict our attention to reals 
 between 0 and 1.

\begin{definition} \label{defofbaseexpan}
A {\em base} is a natural number strictly greater than 1, and a {\em base-$b$ digit} is a natural number in the
set $\{0,1,\ldots, b-1\}$.

Let $b$ be a base, and let $\xd_1,\ldots, \xd_n$ be base-$b$ digits.  We
will use $(0.\xd_1\xd_2\ldots \xd_n)_b$ to denote the rational number 
$0 + \sum_{i=1}^n \xd_i b^{-i}$.

Let    $\xd_1,\xd_2,\ldots$ be an infinite sequence of base-$b$ digits. 
We say that  $(0.\xd_1\xd_2\ldots )_b$ is {\em the base-$b$ expansion of the real number $\alpha$} if we have
\begin{equation*}
(0.\xd_1\xd_2\ldots \xd_n)_b \;\; \leq \;\; \alpha \;\; < \;\;  (0.\xd_1\xd_2\ldots \xd_n)_b + b^{-n}
\end{equation*}
for all $n\geq 1$. Let $E^\alpha_b:\nat^+ \longrightarrow \{0,..,b-1\}$
be the function that yields the $i$th
digit of the base-$b$ expansion of $\alpha$, more precisely, let
 $E^\alpha_b(i) =\xd_i$ 
when $(0.\xd_1\xd_2\ldots)_b$ is the base-$b$ expansion of  $\alpha$.
We will say that $E^\alpha_b$  {\em is the the base-$b$ 
expansion of} $\alpha$.
\qed
\end{definition}

It is easy to see that a base-$b$ expansion can be subrecursively converted into a
Cauchy sequence: If  $(0.\xd_1\xd_2\ldots )_b$ is the base-$b$ expansion of  $\alpha$,
then
$$ (0.\xd_1)_b \; , \; (0.\xd_1\xd_2)_b \; , \; (0.\xd_1\xd_2\xd_3)_b \; , \;  \ldots  $$
will be the first elements of a Cauchy sequence for $\alpha$, and thus, we do
not need unbounded
search to compute a Cauchy sequence if we have access to the base-$b$ expansion.

\begin{lemma} \label{glemtelogsemplanlegging}
Let $E^\alpha_b:\nat\longrightarrow \{0,..,b-1\}$ be 
the base-$b$ expansion of an irrational
$\alpha \in (0,1)$.
There is a parameterized function-oracle Turing Machine $M$ such that 
\begin{itemize}
    \item $\Phi^{E^\alpha_b}_M :\nat \longrightarrow \rational$ is a Cauchy sequence for $\alpha$
    \item $M^{E^\alpha_b}$ on input $n$ runs in time $\poly{n}$ and uses
    $n$ oracle calls of input size at most $\log n$.
\end{itemize}
\end{lemma}

\begin{proof}
On input $n$, $M$ performs $n$ oracle queries to
$E^\alpha_b$, of size at most $\log n$ to obtain
the first $n$ digits $\xd_1, \ldots, \xd_n$ of the base-$b$
expansion of $\alpha$. Each digit requires $\log b$ space,
and computing the rational number
$p/q = \sum_{i=1}^n \xd_i b^{-i}$ can thus be done in time
$\poly{n}$.
\end{proof}

It turns out that we cannot subrecursively convert a Cauchy sequence into a base-$b$
representation. Neither can we in general subrecursively convert a base-$b$
expansion into a  base-$a$ expansion. Let us recall a definition from Kristiansen \cite{ciejourto}.

\begin{definition}
We will use $\primeset{b}$ denote the set of prime factors of the base $b$, that is, 
$\primeset{b}  =  \{ p\mid \mbox{$p$ is a prime and $p|b$}\}$.

Let $a$ and $b$ be bases such that $\primeset{a}\subseteq \primeset{b}$.
We will now define the {\em base transition factor} from  $a$ to $b$.
Let $b= p_1^{k_1} p_2^{k_2} \ldots p_n^{k_n}$, where $p_i$ is a prime and 
$k_i\in\nat^+$ (for $i=1,\ldots,n$),
be the prime factorization of $b$. Then, $a$ can be written of the form
$a= p_1^{j_1} p_2^{j_2} \ldots p_n^{j_n}$ where $j_i\in\nat$ (for $i=1,\ldots,n$).
The {\em base transition factor} from  $a$ to $b$ is the natural number $k$ such that
\begin{equation*}
 k \;\; = \;\; \max\{\, \left\lceil j_i / k_i \right\rceil \, \mid \,1\leq i \leq n \,\}\; .
\end{equation*}
\qed
\end{definition}

Note that the base transition factor from $a$ to $b$ is defined if and
only if $\primeset{a} \subseteq \primeset{b}$ (the definition does not make
sense when $\primeset{a} \not\subseteq \primeset{b}$). 
When we assume that the  base transition factor from $a$ to $b$ exists,  it is understood that we have $\primeset{a}\subseteq \primeset{b}$.

\begin{theorem}[The Base Transition Theorem]\label{thomas}
Let $k$ be the  base transition factor from  $a$ to $b$,
 and let 
$(0.\xd_1\xd_2\ldots)_a$  and  $(0.\dot{\xd}_1\dot{\xd}_2\ldots)_b$
be, respectively, the base-$a$ and base-$b$ expansion of the real number $\alpha$. Then,
for all $n\in\nat$, we have
\begin{multline}  
(0.\xd_1\ldots \xd_n)_a \;\; \leq \;\;  (0.\dot{\xd}_1\ldots \dot{\xd}_{kn})_b \;\; \leq \;\;  \alpha \;\; \\ < \;\;  
(0.\dot{\xd}_1\ldots \dot{\xd}_{kn})_b \; + \; b^{-kn} \;\;  \leq \;\; (0.\xd_1\ldots \xd_n)_a  \; + \;  a^{-n}\; . \tag{I}
\end{multline}
Moreover, for all $n,\ell\in\nat$, we have
\begin{equation}
(0.\dot{\xd}_1\ldots \dot{\xd}_{kn})_b \; < \; (0.\dot{\xd}_1\ldots \dot{\xd}_{\ell})_b \;\; 
\Rightarrow \;\; (0.\xd_1\ldots \xd_n)_a \; < \; (0.\xd_1\ldots \xd_{\ell m})_{a} \tag{II}
\end{equation}
where $m = \lceil \log_a b\rceil$.
\end{theorem}

A proof of the Base Transition Theorem can be found in \cite{ciejourto}.
Assume that the base transition factor $k$ from base $a$ to base $b$ exists.
Then, by clause (I) of the theorem, the $n$ first fractional digits $0.\xd_1\ldots \xd_n$
of the base-$a$
expansions of $\alpha$ will be determined by the $kn$ first fractional digits 
$0.\dot{\xd}_1\ldots \dot{\xd}_{kn}$
of the base-$b$ expansion of $\alpha$, and thus, we can subrecursively convert a base-$b$ expansion into a
base-$a$ expansion.

\begin{lemma}\label{lem:basesmusttransit}
Assume
that the  the base transition factor $k$ from base $a$
to base $b$ exists, and let
$E^\alpha_b:\nat^+ \longrightarrow  \{0, \ldots ,b-1\}$ be 
the base-$b$ expansion of an irrational
$\alpha \in (0,1)$.
There is a parameterized function-oracle Turing Machine $M$ such that 
\begin{itemize}
    \item $\Phi^{E^\alpha_b}_M :\nat^+ \longrightarrow \{0,\ldots ,a-1\} $ is the base-$a$ expansion  of $\alpha$
    \item $M^{E^\alpha_b}$ on input $n$ runs in time $\poly{n}$ and uses $kn$ oracle calls, each of input size at most $O(\log n)$.
\end{itemize}
\end{lemma}

\begin{proof}
 Note that the first $kn$ digits of the $b$-ary expansion of $\alpha$ can be computed using $kn$ oracle calls each of size at most $\log kn$ = $O(\log n)$.
Converting these $kn$ digits to a number of 
the form $q \cdot b^{-kn}$
can be done in time $\poly{kn} = \poly{n}$, and hence computing
$p = \lfloor a^n \cdot q/b^{kn} \rfloor$ can be done in time $\poly{kn} = \poly{n}$.
By Clause (I) of the Base Transition Theorem, $(0.\xd_1\ldots \xd_n)_a = p/a^n$, and as each base-$a$ digit $\xd_1,\ldots,\xd_n$ is in $\{0,\ldots,a-1\}$ and $p/a^n = \sum_{i=1} \xd_i a^{-i}$, a simple greedy  algorithm may compute the digits $\xd_1,\ldots,\xd_n$
of $(0.\xd_1\ldots \xd_n)_a$ in increasing order in time $\poly{n}$ when $p$ has been computed.
\end{proof}

If the base transition factor from $a$ to $b$ does not exist, then  we cannot \label{hvithoytaler}
compute $E^\alpha_a$ subrecursively in $E^\alpha_b$ even if we assume
that $\alpha$ is irrational. 
This is proved formally in \cite{ciejourto}, but intuitively it is not very hard to see why this is the case: 
Consider an irrational $\alpha$ which lies very close to the rational number $(0.1)_{10}$.
We have $(0.1)_{10} = (0.0(0011)^*)_2$, and let us say that
base-2 expansion of $\alpha$ starts with \label{biternegler}
$$
\alpha \; = \; 0.000110011001100110011001100110011001100110011\ldots
$$
Given the digits of the base-2 expansions displayed above, we cannot tell if the first fractional digit of the base-10 should be 0 or 1. 
We need more digits of the base-2 expansion to the determine the first digit of the base-10 expansion.  Now, $\alpha$ cannot
equal $(0.1)_{10}$ since we have assumed that $\alpha$ is irrational.
Thus, sooner or later we will find a digit in
the base-2 expansion  which allows to determine the first digit of the base-10 expansion,
but we  need unbounded search to find that digit. 

In general, if the base transition factor from $a$ to $b$ does not
exist, that is, if $\primeset{a}\not\subseteq\primeset{b}$, we cannot
subrecursively compute $E^\alpha_a$ in $E^\alpha_b$. It follows that
we cannot compute $E^\alpha_b$ subrecursively in a Cauchy sequence 
for $\alpha$ (for any base $b$). Assume for the sake of a contradictions
that we can, that is, assume that we can compute $E^\alpha_b$
subrecursively in an arbitrary Cauchy sequence for $\alpha$. Pick a base $b_0$ such that $\primeset{b}\not\subseteq\primeset{b_0}$. By Lemma
\ref{glemtelogsemplanlegging}, we can subrecursively compute a Cauchy sequence $C$ for $\alpha$ in $E^\alpha_{b_0}$. By our assumption we
can subrecursively compute $E^\alpha_b$ in $C$. Hence, we can subrecursively
compute $E^\alpha_{b}$ in $E^\alpha_{b_0}$ which is impossible 
as the base transition factor from $b$ to $b_0$ does not exist.

\subsection{Base-$b$ sum approximations.}
\label{shabatbaboker}

Base-$b$ sum approximations (from below and above) were
introduced by Kristiansen in \cite{ciejouren} and studied further
in \cite{ciejourto} and, with Georgiev and Stephan, in \cite{apalilf}.

\begin{definition} \label{amandaborhosmeg}
Let $(0.\xd_1\xd_2\ldots)_b$ be the base-$b$ expansion of the irrational 
$\alpha\in (0,1)$ (thus, we have $E^\alpha_b(n)= \xd_n$).

The  {\em  base-$b$ sum approximation from below} of $\alpha$  is the function $\hat{A}^{\alpha}_b:\nat\longrightarrow \rational$ defined
by  $\hat{A}^{\alpha}_b(0)= 0$ and 
$ \hat{A}^{\alpha}_b(n+1)   =   E^\alpha_b(m)   b^{-m} $
where $m$ is the least $m$ such that 
\begin{equation*}
\sum_{i=0}^n \hat{A}^\alpha_b(i) \;\; < \;\; (0.\xd_1\ldots \xd_m)_b
\end{equation*}
that is, $\hat{A}^{\alpha}_b(n)$ is the value represented by the $n$th non-zero digit of the base-$b$ expansion of $\alpha$.

Let $\overline{\xd}$ denote the {\em complement digit} of 
the base-$b$ digit $\xd$, that is, let $\overline{\xd}= (b-1)-\xd$ 
(observe that we have 
$(0.\xd_1 \xd_2 \xd_3\ldots)_b + (0.\overline{\xd}_1\overline{\xd}_2 \overline{\xd}_3\ldots)_b  = 1$ 
for any base $b$ and any base-$b$ expansion $(0.\xd_1 \xd_2\ldots)_b$).

The {\em  base-$b$ sum approximation from above} of $\alpha$  is the function $\check{A}^{\alpha}_b:\nat\longrightarrow \rational$ defined
by  $\check{A}^{\alpha}_b(0)= 0$ and 
$ \check{A}^{\alpha}_b(n+1)   =   \overline{E^\alpha_b(m)} b^{-m}$
where $m$ is the least $m$ such that 
\begin{equation*} 1 \; - \; \sum_{i=0}^n \check{A}^\alpha_b(n) 
\;\; >  \;\;  1 \; - \; (0.\overline{\xd}_1\ldots \overline{\xd}_m)_b \; .
\end{equation*}
\qed
\end{definition}

The functions $\hat{A}^{\alpha}_b$  and $\check{A}^{\alpha}_b$ are not defined if $\alpha$ is rational.
When we use the notation it is understood that $\alpha$ is irrational.
It is fairly straightforward to prove that
\begin{equation*}
\sum_{i=0}^\infty E^\alpha_b(i)b^{-i}  \;\; =  \;\;  \sum_{i=0}^\infty \hat{A}^\alpha_b(i)  \;\; =  \;\;       1 \; - \;       \sum_{i=0}^\infty \check{A}^\alpha_b(i) \; .
\end{equation*}
A detailed proof can be found in  \cite{ciejourto}.

We cannot subrecursively compute  $\hat{A}^\alpha_b$ 
in $\check{A}^\alpha_b$, and neither can we 
subrecursively compute $\check{A}^\alpha_b$ in $\hat{A}^\alpha_b$ (for any 
base $b$). The following growth argument (see page \pageref{growthpage})
 explains why we cannot subrecursively compute $\hat{A}^\alpha_2$ 
in $\check{A}^\alpha_2$: Let $t$ be a any time-bound, and let $f:\nat \longrightarrow \nat$
a monotone strictly increasing function that
grows faster than any function computable in time $O(t)$, moreover, let
the graph of $f$ be computable in
polynomial time, that is, the relation $f(x)=y$ can be decided in
 time polynomial in the size of the natural numbers  $x$ and $y$. It is straightforward to see that such an $f$ exists, and that we may w.l.o.g.~assume
that $2^x \leq f(x)$. Consider the
irrational number $\beta$ given by $\hat{A}^\beta_2(n)=2^{-f(n)}$.
We have
$$\beta = (0.0\ldots 010\ldots 010\ldots 010\ldots 010\ldots)_2$$ where the sequences 
$0\ldots 0$ of zeros are getting longer and longer.
Now
$$\beta = (0.0\ldots 010\ldots 010\ldots 010\ldots)_2
= 1 - (0.1\ldots 101\ldots 101\ldots 101\ldots)_2
$$
and thus $\check{A}^\beta_2(n)=2^{-g(n)}$ where $g:\nat \longrightarrow \nat$
is a slow growing function (we have $g(x)\leq 2x$), indeed, 
$g(x)$ is computable in polynomial time of the natural number $x$ as the graph of $f$ is computable in polynomial time. Since we can compute $g$ in polynomial time, we
can also compute
$\check{A}^\beta_2$  in polynomial time. Obviously, we cannot compute 
$\hat{A}^\beta_2$ in  time $O(t)$, if we could, then we could also
compute $f$ in  time $O(t)$, contrary to our assumption.

Hence we conclude that for any time-bound $t$ there exists
an irrational $\beta$ such that $\check{A}^\beta_2$ is computable in
polynomial time whereas $\hat{A}^\beta_2$ is not computable in time 
$O(t)$. This shows that we cannot subrecursively compute
$\hat{A}^\alpha_2$ in $\check{A}^\alpha_2$. The argument generalizes
easily to work for any base $b$, and hence, we cannot subrecursively compute
$\hat{A}^\alpha_b$ in $\check{A}^\alpha_b$.
A symmetric arugment
will show that we cannot subrecursively compute  $\check{A}^\alpha_b$ in
$\hat{A}^\alpha_b$. 
Detailed proofs can be found in  \cite{ciejourto}.

It is also proved in \cite{ciejourto} that the base transition factor
from $a$ to $b$ exists, if and only if, we can 
subrecursively compute $\hat{A}^\alpha_a$ in $\hat{A}^\alpha_b$,
if and only if, we can 
subrecursively compute $\check{A}^\alpha_a$ in $\check{A}^\alpha_b$.
We have already argued why
we cannot subrecursively compute $E^\alpha_a$ in $E^\alpha_b$
when the base transition factor from $a$ to $b$ does not exists (see page
\pageref{hvithoytaler}).
The very same argument 
should also give an intuitive explanation of 
why we cannot subrecursively compute
$\hat{A}^\alpha_a$ in $\hat{A}^\alpha_b$, or  $\check{A}^\alpha_a$ in $\check{A}^\alpha_b$, when the base transition factor from $a$ 
to  $b$ does not exist. In the next subsection we will analyze the
complexity of computing $\hat{A}^\alpha_a$ in $\hat{A}^\alpha_b$
when the needed base transition factor is available.

\subsection{From base-$b$ sum approximations to base-$a$ sum approximations.}

Assume that the base transition factor $k$ from base $a$ to base $b$ exists.
We will now give and explain an algorithm for computing $\hat{A}^\alpha_a$ ($\alpha$'s base-$a$ sum approximation from below)
using $\hat{A}^\alpha_b$ ($\alpha$'s base-$b$ sum approximation from below) as an oracle. Of course,
$\alpha$'s base-$a$ sum approximation from above   can be computed from $\alpha$'s base-$b$ sum approximation from above by a symmetric
algorithm.

%The algorithm works recursively, and we 
%ABA: I find the mention of recursion distracting, the algorithm is actually iterative
Assume
  the values of $\hat{A}^\alpha_a(0), \ldots , \hat{A}^\alpha_a(n)$ already are computed. Then, the algorithm
computes the value of $\hat{A}^\alpha_a(n+1)$ by carrying out the following instructions:
\begin{itemize}
    \item Step 1: Compute the rational number  $\sum_{i=0}^n \hat{A}^\alpha_a(i)$. The number will be of the form 
    $(0.\xd_1\ldots \xd_p)_a$
    for some $p$. Compute that $p$.
    \item Step 2: Ask the oracle $\hat{A}^\alpha_b$ for the value of $\hat{A}^\alpha_b(kp+1)$
    (where $k$ is the base transition factor). The oracle will yield
    a rational number of the form $\xd b^{-\ell}$ where $\ell\geq kp+1$ and $\xd$
    is a nonzero base $b$ digit. Compute $\ell$.
    \item Step 3: Use the oracle    $\hat{A}^\alpha_b$ to compute the rational number 
    $$R = \sum_{i=0}^{k\ell m} \hat{A}^\alpha_b(i)$$
    where $k$ is the base transition factor, $m = \lceil \log_a b\rceil$ and $\ell$
    is the value computed in Step 2.
    \item Step 4: Compute the least $i$ such that digit number $p+i$ in the base-$a$
    expansion of $R$ is nonzero (where $p$ is the value computed in Step 1 and $R$ is the value computed in Step 3).
     Give the output $\xd a^{-(p+i)}$ where $\xd$ is digit number $p+i$ in the base-$a$ expansion of $R$, that is, we have $\hat{A}^\alpha_a(n+1)= \xd a^{-(p+i)}$
\end{itemize}

We will now argue that the algorithm gives correct output.
Let
$(0.\xd_1\xd_2\ldots)_a$  and  $(0.\dot{\xd}_1\dot{\xd}_2\ldots)_b$
be, respectively, the base-$a$ and base-$b$ expansion of  $\alpha$.
In Step 1, the algorithm computes  $(0.\xd_1\ldots \xd_p)_a$.
According  to the definition of  $\hat{A}^\alpha_a$, the output should be
$\xd_{p+i}a^{-(p+i)}$ where $i$ is the least number such that $$(0.\xd_1\ldots \xd_p)_a < (0.\xd_1\ldots \xd_{p+i})_a \; .$$
In step 2, the algorithm computes $\ell$. By the definition of sum
approximations,
\begin{equation} \label{shabatbaerev}
(0.\dot{\xd}_1\ldots \dot{\xd}_{kp})_b  <  (0.\dot{\xd}_1\ldots \dot{\xd}_{\ell})_b \; .
\end{equation}
Hence, by clause (II) of  Base Transition Theorem, we have 
\begin{equation} \label{italiener}
    (0.\xd_1\ldots \xd_p)_a < (0.\xd_1\ldots \xd_{\ell m})_a 
\end{equation}
where $m=\lceil \log_a b \rceil$. This shows that the denominator of the
next term in the base-$a$ sum approximation is at most $a^{\ell m}$.
In Step 3, the algorithm computes $R$. By \eqref{italiener} and clause (I) of the Base Transition Theorem, we have
\begin{multline}
\label{erling}
(0.\xd_1\ldots \xd_p)_a < 
  (0.\xd_1\ldots \xd_{\ell m})_a \leq 
     (0.\dot{\xd}_1\ldots \dot{\xd}_{k \ell m})_b  = R < \alpha   \\
<  (0.\dot{\xd}_1\ldots \dot{\xd}_{k \ell m})_b + b^{-k\ell m} \le
  (0.\xd_1\ldots \xd_{\ell m})_a + a^{-\ell m}
     \; .  
\end{multline} 
Now, (\ref{erling}) implies that the first $\ell m$ fractional digits of the base-$a$ expansion of $\alpha$, that is $\xd_1\ldots \xd_{\ell m}$,
coincide with the first $\ell m$ fractional digits of the base-$a$ expansion of the rational number $R$. 
Thus, the algorithm computes a correct result in Step 4.

\begin{lemma}
Assume
that the  the base transition factor $k$ from base $a$
to base $b$ exists, and let $\hat{A}^\alpha_b : \nat \longrightarrow \rational$ is 
the base-$b$ sum approximation from below of an irrational 
$\alpha \in (0,1)$. Furthermore, assume
that $\hat{A}^\alpha_b$ is computable within a time bound $s$, and let
$g(n) = (\lambda x.\bitlen{m}+s(\bitlen{k}+x+1))^{(n)}(0)$ where $m=\lceil \log_a b \rceil$.
There is a parameterized function-oracle Turing Machine $M$ such that
\begin{itemize}
    \item $\Phi^{\hat{A}^\alpha_b}_M :\nat \longrightarrow \rational$ is the base-$a$ sum approximation from below of $\alpha$
    \item $M^{\hat{A}^\alpha_b}$ on input $n$ runs in time $\poly{2^{g(n)}+s(\bitlen{k}+g(n))}$
    and uses at most $k2^{g(n)}$
oracle calls, each of input size at most $\bitlen{k}+g(n)$.
\end{itemize}
\end{lemma}

\begin{proof}
By induction on $n$, we prove the time bound and also a bound on the exponent of $a$ in the denominator
of $\hat{A}^\alpha_a(n)$.

To compute $\hat{A}^\alpha_a(1)$, we call the oracle with input 1 and obtain
a position $\ell$; we then compute $R$ to obtain from it the first non-zero digit 
in base $a$. By assumption,  $\hat{A}^\alpha_b(1)$ is computable in time $s(\bitlen{1})$; this
implies that $s(\bitlen{1})$ is also a bound on the bit-length of the result of this computation, hence on that of $\ell$ (the result of the
computation is the binary representation of $c\cdot b^{-\ell}$ for some $c$).  The machine has to access values of $\hat{A}^\alpha_b(i)$ for
$i=1,\dots,k\ell m$, so it uses 
$k\ell m \le km 2^{s(\bitlen{1})} = k2^{g(1)}$ oracle calls.
Since the size of a result is bounded by its computation time, we have
that the size of the largest oracle answer is bounded by $s(\bitlen{k} + g(1))$.
It is routine to
verify that the execution time is dominated by the expression $k\ell m$, which is the length of the summation in Step~3, times the size of the largest number returned from an oracle call.
We have bounded $k\ell m$ by $k2^{g(1)}$, and the size of the last oracle call is bounded by
$s(\bitlen{k\ell m}) \le s(\bitlen{k}+g(1))$, justifying the bound on execution time for $n=1$.
We note that
the exponent in the denominator of  $\hat{A}^\alpha_a(1)$ is bounded by $\ell m$ 
\eqref{italiener}, hence by $2^{g(1)}$.

We turn to the induction step. Assume $n>1$, and assume inductively that the position $p$ is bounded by $2^{g(n-1)}$.
The value $\ell$ will be bounded by $2^{s(\bitlen{kp+1})}$ and the
largest input to an oracle call will be 
$$ k\ell m \; \le \; 
k m 2^{s(\bitlen{kp+1})} \; \le \; km 2^{s(\bitlen{k}+\bitlen{p}+1)}
\;  \le  \;  km 2^{s(\bitlen{k}+g(n-1)+1)}\; \le  k2^{g(n)}\;. $$
This also bounds the number of oracle
calls, since we can store results of previous queries and therefore never
query the oracle on the same input twice.
A bound on the bit-length of the result of the oracle call is 
$s(\bitlen{k} + g(n))$. 
The execution time is polynomial in the sum of this quantity and the number of calls,
yielding the bound $\poly{2^{g(n)}+s(\bitlen{k}+g(n))}$,
and the position of the last base-$a$ digit is bounded
by $\ell m \le 2^{g(n)}$.  
\end{proof}

We note that the lemma implies that if $\hat{A}^\alpha_b$ is computable by an 
$O(s)$-time Turing machine,  then $\hat{A}^\alpha_a$ is computable 
by an $O(s')$-time Turing machine where the time-bound $s'$ is primitive recursive in the 
time bound $s$.

\subsection{From base-$b$ sum approximations to base-$a$ expansions.}

This is very similar to converting
base-$b$ expansions to base-$a$ expansions.
As in Section~\ref{shabatbaboker}, we argue that if the base transition factor from $a$ to $b$ does not exist, then conversion from base-$b$ sum approximations to base-$a$ expansions cannot be done subrecursively. In the case that the transition factor exists, we have a straightforward
 conversion algorithm. 

\begin{lemma} \label{femblyanter}
Assume
that the  the base transition factor $k$ from base $a$
to base $b$ exists, and let $\hat{A}^\alpha_b : \nat \longrightarrow \rational$ be the base-$b$ sum approximation from below of an irrational $\alpha \in (0,1)$.
There is a parameterized function-oracle Turing machine $M$ such that
\begin{itemize}
    \item $\Phi_M^{\hat{A}^\alpha_b} : \nat \longrightarrow \{0,\ldots,a-1\}$
is the base-$a$ expansion of $\alpha$
\item $M^{\hat{A}^\alpha_b}$ on input $n$ runs in time $\poly{\bitlen{\hat{A}^\alpha_b(kn)}}$
and uses at most $kn$ oracle calls, each of input size $O(\log n)$.
\end{itemize}
\end{lemma}

\begin{proof}
As in the proof of Lemma \ref{lem:basesmusttransit}, computing the $n$th digit of $(0.\xd_1\xd_2\xd_3 \cdots)_a$ requires obtaining
the first $kn$ digits of $(0.\xd_1\xd_2\xd_3 \cdots)_b$.
By definition of $\hat{A}^\alpha_b$, if $\hat{A}^\alpha_b(n) = c \cdot b^{-m}$ for some $m \in \nat$ and $c \in \{0,\ldots,b-1\}$,
then the $n$th non-zero digit of $(0.\xd_1 \xd_2 \xd_3 \cdots)_b$ is $c$. Hence, to find the
first $kn$ digits of the base-$b$ expansion, $M$ may simply query $\hat{A}^\alpha_b(1),\hat{A}^\alpha_b(2),\ldots$ in order until the first $i$ 
is found for which $\hat{A}^\alpha_b(i) = c_i b^{-m_i}$ with $m_i \geq kn$. 
This is obviously a bounded search as such an $m_i$ will exist for some $i\leq kn$.
Once we have
$(0.\xd_1 \xd_2\dots \xd_{kn})_b$, we extract $(0.\xd_1 \xd_2 \xd_3 \dots \xd_n)_a$  by 
repeated multiplication by $b$ and division by $a$.
Clearly, this procedure uses at most $kn$ oracle calls, each of size at most $\log kn$.
The time complexity of the procedure is dominated by a polynomial in the size of the result of the last oracle call,
i.e., $\poly{\bitlen{\hat{A}^\alpha_b(kn)}}$ (note that the denominator of this number is at least $b^{kn}$).
\end{proof}

Lemma \ref{femblyanter} does of course also hold for sum approximations from above.

\subsection{Gray codes.}

The Grey code representation of real numbers was introduced by
Tsuiki \cite{DBLP:journals/tcs/Tsuiki02} and studied further 
in Berger et al.~\cite{berger}.
The representation is subrecursively
 equivalent to the representation by base-2 expansions. 

\begin{definition}
The function $G: \nat \longrightarrow \{0,1\}$ is the 
\emph{Gray code} of the irrational number
$\alpha$
 if $G(i) = 0$
if there is an even number $m$ with 
$$
m2^{-i} - 2^{-(i+1)} < \alpha < m 2^{-i} + 2^{-(i+1)}
$$
and $G(i) = 1$ if the same holds for an odd number $m$.
\qed
\end{definition}

Gray codes  are usually defined as maps $G: \nat \longrightarrow \{0,1,\bot\}$
where $G(i) = \bot$ if
$\alpha$ is of the form $m 2^{-i} - 2^{-(i+1)}$
for some $m$ and $i$, and hence is rational. As we are only considering irrationals, 
we do not need $\bot$.

\begin{lemma}[Tsuiki \cite{DBLP:journals/tcs/Tsuiki02}] \label{trustbust}
Let $E_2^\alpha$ and $G^\alpha$, respectively,  be the base-2 expansion and the Grey code of the irrational  $\alpha\in (0,1)$. 
Then we have
\begin{align*}
 E_2^\alpha(1) &= G^\alpha(0) \\
 E_2^\alpha(n+1) &= E_2^\alpha(n)\xor G^\alpha(n)
\end{align*}
where $\xor$ denotes the XOR function.
\end{lemma}

\begin{lemma}  \label{xtrustbust}
Let $G : \nat \longrightarrow \{0,1\}$ is the
Gray code of an irrational  $\alpha \in (0,1)$.
There is a parameterized function-oracle Turing Machine $M$ such that 
\begin{itemize}
    \item $\Phi^G_M :\nat^+ \longrightarrow \{0,1\}$ is the base-2 expansion of $\alpha$
    \item $M^G$ on input $n$  runs in time $O(n\log n)$ and uses $n$ oracle calls,
    each of input size  $O(\log n)$.
\end{itemize}
\end{lemma}

\begin{proof}
Straightforward from the conversion algorithm expressed by Lemma \ref{trustbust}.
\end{proof}

\begin{lemma}
Let  $E : \nat^+ \longrightarrow \{0,1\}$ be the base-2 expansion of
an irrational  $\alpha \in (0,1)$. There is a parameterized function-oracle Turing Machine $M$ such that 
\begin{itemize}
    \item $\Phi^{E}_M :\nat \longrightarrow \{0,1\}$ is the Gray code of $\alpha$
    \item $M^E$ on input $n$ runs in time $O(\log n)$ and uses $2$ oracle calls, each
    of input size  $O(\log n)$.
\end{itemize}
\end{lemma}

\begin{proof}
A simple rewrite of the equations in Lemma \ref{xtrustbust} give:
\begin{align*}
 G^\alpha(0) &= E_2^\alpha(1) \\
 G^\alpha(n+1) &= E_2^\alpha(n+1)\xor E_2^\alpha(n)
\end{align*} 
The implied algorithm has the complexity stated in the lemma.
\end{proof}

\subsection{Summary.}
In this section, our main results concerned the representations by
base-$b$ expansions (Gray codes being equivalent to base-2 expansions),
and the representations by base-$b$ sum approximations. These 
representations form clusters in which not all representations allow
for subrecursive conversion from one to the other. Base-$b$ expansions
are convertible to base-$a$ expansions when the base transition factor
from $a$ to $b$ exists---and then the overhead of the conversion is
exponential in the bit-length of the input. Base-$b$ sum approximations
from below (above) are convertible to base-$a$ sum approximations
from below (above) when the base transition factor
from $a$ to $b$ exists---and then the overhead of the conversion involves iteration, which guaranteed that an $O(s)$-time computable
base-$b$ representation becomes an $O(t)$-time computable base-$a$
representation with $t$ primitive recursive in $s$.

\newpage

%% file: bolzano.tex
\section{Representations Subrecursively Equivalent to Dedekind Cuts}

\label{skittenskje}

\subsection{Dedekind cuts.}

\begin{definition} \label{fiskebullyduckdolly}
Let $\alpha \in (0,1)$ be an irrational number.  The  {\em Dedekind
cut} of  $\alpha$ is the function $D^\alpha:\rational \longrightarrow \{0,1\}$ 
given by  $D^\alpha(q)=0$ iff  $q< \alpha$.
\qed
\end{definition}

The representation of irrational numbers by Dedekind cuts is not
subrecursively equivalent to the representation by base-$b$ expansions,
 neither is it subrecursively equivalent to the
representation by base-$b$ sum approximations from below or above, for any base $b$. It is fairly easy to see that we can subrecursively compute
base-$b$ expansion $E^\alpha_b$ in $D^\alpha$, see Section
\ref{bilboars} below, but it is not possible to subrecursively compute
  $D^\alpha$ in $E^\alpha_b$ for any fixed base $b$. An intuitive 
  explanation of why this is impossible is very similar to our
explanation   of why we cannot subrecurively compute $E_{10}^\alpha$
in $E_{2}^\alpha$ at page \pageref{biternegler}:
Let $\alpha$ be an irrational which 
base-2 expansion starts with 
$$
\alpha \; = \; 0.000110011001100110011001100110011001100110011\ldots \; .
$$
The period $0011$ may be repeated arbitrarily  many times, and thus,
we will need unbounded search to determine if $\alpha$ lies
above or below $10^{-1}$, that is, we need unbounded search to
compute $D^\alpha(10^{-1})$. This simple example should also serve as an 
intuitive explanation of
why we cannot subrecursively compute $D^\alpha$ in $\alpha$'s base-$b$
sum approximation from above or below, that is, $\hat{A}^{\alpha}_b$
or $\check{A}^{\alpha}_b$. Neither is it possible to subrecursively compute 
$\hat{A}^{\alpha}_b$,
or $\check{A}^{\alpha}_b$, in $D^\alpha$, but an intuitive explanation 
of why this is the case is not all that straightforward, and we refer
the interested reader to Section 7 and 8 of Kristiansen \cite{ciejourto}
for more
on the relationship between Dedekind cuts and base-$b$ sum approximations.

\subsection{Definitions.}
The next definition gives some representations being subrecursively equivalent to the
representation by Dedekind cuts.

\begin{definition} \label{bullyduckdolly}
Let $\alpha \in (0,1)$ be an irrational number. 
\begin{enumerate}
\item The 
{\em Beatty sequence} of  $\alpha$ is the
function $B^\alpha:\nat^+ \longrightarrow \nat$ given by
$$  \frac{B^\alpha(n)}{n} \; < \; \alpha \; < \; \frac{B^\alpha(n)+1}{n}\; .   $$
\item 
The {\em general base expansion} of  $\alpha$ is  the function 
$$E^\alpha:  (\nat \setminus \{0,1\} ) \times    \nat^+ \longrightarrow \nat$$
where $E^\alpha(b,n)=  E^\alpha_b(n)$ (recall that $E^\alpha_b$ is the base-$b$ expansion of $\alpha$, see Definition \ref{amandaborhosmeg}).
\item
The \emph{Hurwitz characteristic} of 
$\alpha$ is the function $H^\alpha : \mathbb{N} \longrightarrow \{0,1\}^*$
such that $H^\alpha(0),H^\alpha(1),H^\alpha(2),\ldots$ is a 
path in the Farey pair tree\footnote{Strictly speaking, the classic Hurwitz characteristic corresponds to a path through the full Stern-Brocot tree (not  the  Farey pair tree as we consider here), and hence the classic Hurwitz characteristic $H'$ of $\alpha \in (0,1)$ is the function defined by $H'(0) = 0$ and $H'(q) = 0 \cdot H(q-1)$ for $q > 0$. This does not change our results in any material way.} $\fareypairtree$, 
and moreover, for all $n\in \nat$, we have $\alpha \in \fareypairtree(H^\alpha(n))$.
\end{enumerate}
\qed
\end{definition}

Representations by Dedekind cuts \cite{Bertrand,Dedekind}, Beatty sequences \cite{beattyseq} and Hurwitz characteristics \cite{hurwitz} were known in the 19th century or earlier. Use of the Hurwitz characteristic to represent numbers rather than a stepping stone for other material is a much younger invention, see Lehman \cite{lehman}. Moreover, what is now known as Beatty sequences was used earlier by Bernard Bolzano \cite{bolzano}, whence this representation of reals could also be called \emph{Bolzano measures}.
The representation by  general base expansions might not
have been investigated before, but it is  very natural.

This section is based on the conference paper Kristiansen \& Simonsen \cite{dedekindcie}.

\subsection{Conversion between general base expansions and De\-dekind cuts.}
\label{bilboars}
We will compute $E^\alpha(b,n)$ by  computing the digits $\xd_1,\xd_2, \xd_3,\ldots$ of $\alpha$'s
base-$b$ expansion one by one. When we have determined the digits 
$\xd_1,\ldots, \xd_{n}$, we know that
$$
(0.\xd_1\ldots \xd_{n})_b < \alpha <  (0.\xd_1\ldots \xd_{n})_b + b^{-n}
$$
and then we can split the interval
$$
( \ (0.\xd_1\ldots \xd_{n})_b \, , \,  (0.\xd_1\ldots \xd_{n})_b + b^{-n} \ )
$$
into $b$ subintervals, each of length $b^{-n-1}$, and 
use the Dedekind cut of $\alpha$ to determine the digit $D_{n+1}$.

\begin{lemma}\label{lem:G_from_ded}
Let
$D: \mathbb{Q} \longrightarrow \{0,1\}$
be the Dedekind cut of ab irrational $\alpha \in (0,1)$.
There is a parameterized function-oracle Turing machine $M$ such that
\begin{itemize}
    \item $\Phi_\paramorac{M}{D}
: (\mathbb{N} \setminus \{0,1\}) \times\mathbb{N}^+ \longrightarrow \nat$ is the general base expansion of $\alpha$
\item $M^D$ on input $(b,n)$ runs in time $O(b^2\poly{n})$ and uses at most $n \log b$ 
oracle calls, each  of input size at most $O(n\log b)$.
\end{itemize}
\end{lemma}

\begin{proof}
$M$ constructs the sequence $E^{\alpha}_b(1), E^{\alpha}_b(2), \ldots, E^{\alpha}_b(n)$
inductively by maintaining an open interval $I_i = (v_i,w_i)$ with rational endpoints
$v_i,w_i \in \mathbb{Q}$ for each $i \in \{0,\ldots,n-1\}$ such that (i) $\alpha \in I_i$,
(ii) $v_i$ is a multiple of $b^{-i}$, and (iii) $w_i - v_i = b^{-i}$. 
Initially,
$I_0 = (0,1)$. For each interval $I_i$, $M$ splits $I_i = (v_i,w_i)$ into $b$ equal-sized
intervals 
$$
(v_i,v_i + b^{-(i+1)}), \ldots, (v_i+(b-1)b^{-(i+1)},v_i+b^{-i}) = (v_i+(b-1)b^{-(i+1)},w_i)\; .
$$
Observe that, for any interval $(r_1,r_2)$, if $D(r_1) = D(r_2) = 0$, then
$\alpha > r_2$, and if $D(r_1) = D(r_2) = 1$, then $\alpha < r_1$
(and the case $D(r_1) = 1 \land D(r_2) = 0$ is not possible).
Thus,
$M$ can use $D$ to perform binary search on (the endpoints of) the above set of intervals to find
the  interval 
\begin{align} \label{hvorforikke}
    ( \ v_i + jb^{-(i+1)} \, , \,  v_i + (j+1)b^{-(i+1)} \ )
\end{align}
that contains $\alpha$ (observe that, for this interval, $D(v_i + jb^{-(i+1)}) = 0$
and $D(v_i + (j+1)b^{-(i+1)}) = 1$). We then set $(v_{i+1},w_{i+1})$
to equal the interval (\ref{hvorforikke}).
By construction, we have $E^\alpha(b,i+1) = j$.

Clearly, in each step $i$, there are at most $\log b$ oracle calls to $D$,
and the construction of each of the $b$ intervals and writing on the query tape can be performed
in time polynomial in the binary representation of the numbers involved, hence in time 
$O(\polylog{b^{i}}) = O(\poly{i} \polylog{b})$. Hence,
the total time needed to produce $E(b,n)$ is at most $O(b n \poly{n} \polylog{b}) = O(b^2\poly{n})$
with at most $n \log b$ queries to $D$. 
In each oracle call, the rational numbers involved are all endpoints of intervals where the endpoints are sums of negative powers of $b$ and where the exponent of all powers are at most $n$. Hence, all oracle calls can be represented
by rational numbers using at most $O(n\log b)$ bits.
\end{proof}

Our algorithm for converting a general base expansion to a Dedekind cut is
based on the following observation: For any $n,m\in \nat$ such that $0<n/m<1$,
we have $D^\alpha(n/m) = 0$ iff $n/m<\alpha$ iff $n\leq E^\alpha(m,1)$.

\begin{lemma}\label{lem:ded_from_G}
Let $E:(\nat \setminus \{0,1\} ) \times    \nat^+ \longrightarrow \nat$ be
the general base expansion of an irrational
$\alpha \in (0,1)$. There is a parameterized function-oracle Turing machine $M$ such that
\begin{itemize}
    \item $\Phi_\paramorac{M}{E}
: \mathbb{Q} \longrightarrow \{0,1\}$ is the Dedekind cut of $\alpha$
\item $M^E$ on input $n/m$  runs in time $O(\log(\max\{n,m\}))$ and uses 
exactly one oracle call of input size at most $O(\log(m))$.
\end{itemize}
\end{lemma}

\begin{proof}
On input $n/m \in \mathbb{Q}$, $M$ first checks if $m = 1$, and outputs $0$ if $n \leq 0$
and $1$ if $n \geq 1$. Otherwise, $m > 1$, and
$M$ computes $E(m,1)$; by definition, this is an element of 
$\{0,\ldots,m-1\}$. Thereafter,
$M$ outputs $0$ if $n \leq E(m,1)$, and outputs $1$ otherwise.
$M$ needs to write the (representation of the) pair $(m,1)$ on the oracle tape and perform a single
comparison of numbers of magnitude at most $\max\{n,m\}$,
hence $M$ uses time $O(\log\max\{n,m\})$ for the comparison. $M$ uses exactly one oracle call to $E$ with the pair $(m,1)$, the representation
of which uses  at most $O(\log m)$ bits. 
\end{proof}

\subsection{Conversion between Beatty sequences and Dedekind cuts.}

It is easy to see how we can convert a Dedekind cut $D^\alpha$ into a 
Beatty sequence $B^\alpha$ as the value of $B^\alpha(n)$ is the 
natural number $m$ such that $m/n < \alpha < (m+1)/n$. We may use
$D^\alpha$ to search for that unique $m$.

\begin{lemma}\label{lem:ded_to_bol}
Let $D: \mathbb{Q} \longrightarrow \{0,1\}$
be the Dedekind cut of an irrational $\alpha \in (0,1)$.
There is a parameterized function-oracle Turing machine $M$ such that 
\begin{itemize}
    \item $\Phi_\paramorac{M}{D}
: \mathbb{N}^+ \longrightarrow \mathbb{N}$ is the Beatty sequence of $\alpha$
\item $M^D$ on input $n$ runs in time $O(\polylog{n})$ and uses
at most $\lceil \log n \rceil$ oracle calls, each   of input size at most $O(\log n)$.
\end{itemize}
\end{lemma}

\begin{proof}
On input $n$, $M$ finds the least $i \in \{1,\ldots,n\}$ such that
$D(i/n) = 1$. 
As $D(i/n) = 1$ and $j >i$
implies $D(j/n) = 1$, the least $i$ can be found by binary search,
halving the search range in each step. This can be done by maintaining two integers $l$ and $u$
ranging in $\{1,\ldots,n\}$, and requires a maximum of $\log n$ halving steps. In each halving step, $M$ finds the midpoint $m$ between $l$ and $u$, writes its binary representation on the query tape,
 queries $D$, and records the answer. Then, $l$ and $u$ are updated using basic binary arithmetic operations on integers, represented by at most $O(\log n)$ bits; 
 if $D(m/n) = 1$, then $u := m$, and if $D(m/n) = 0$, then $l := m$.
Clearly, in each step, the arithmetic and update operations can be performed in time polynomial in the size of the representation of the integers,
hence in time $\polylog{n}$.
As $(i-1)/n < \alpha < i/n$, we have $B(n) = i-1$, and $\paramorac{M}{D}$ thus returns $i-1$.
\end{proof}

In order to see that our algorithm for converting a Beatty sequence $B^\alpha$
into a Dedekind cut $D^\alpha$ is correct, observe that
we have $D^\alpha(n/m)=0$ iff $n/m< \alpha$ iff 
$n\leq B^\alpha(m)$, for any $m,n\in \nat$,

\begin{lemma}\label{lem:bol_to_ded}
Let $B: \mathbb{N}^+ \longrightarrow \mathbb{N}$
be the Beatty sequence of an irrational  $\alpha \in (0,1)$.
There is a parameterized function-oracle Turing machine $M$  such that 
\begin{itemize}
    \item $\Phi_\paramorac{M}{B}
: \mathbb{Q} \longrightarrow \{0,1\}$ is the Dedekind cut of $\alpha$
\item $M^B$ on input $n/m$ runs in time $O(\log(\max\{n,m\}))$ and uses exactly one oracle call of input size $O(\log m)$.
\end{itemize}
\end{lemma}

\begin{proof}
On input  $n/m$,  $M$  perform the oracle call $B(m)$,
resulting in an integer $B(m)$ (where $B(m) \in \{0,1\ldots,m-1\}$). 
If $n\leq B(m)$, then $M$ outputs 0, otherwise, $M$ outputs 1.
The comparison
$n \leq B(m)$ can be performed bitwise using the binary
representations of $n$ and $B(m)$ which is clearly linear
in $\log(\max\{n, B(m)\}) \leq \log\max\{n, m\}$.
Writing $m$ on the oracle tape clearly also takes time linear
in $\log m$. 
\end{proof}

\subsection{Conversion between Hurwitz characteristics and De\-dekind cuts.}

Let us first discuss how we can compute a Dedekind cut $D^\alpha$
using a Hurwitz characteristic $H^\alpha$ as an oracle. 
   
Assume  $0<n/m<1$
where $n$ and $m$ are relatively prime natural numbers. 
Let $(a/b, c/d)= \fareypairtree(H^\alpha(n+m))$. 
By Proposition \ref{prop:alleskalslankesig}, any reduced fraction $n/m$ occurs as one of the fractions in a Farey pair in $\fareypairtree$ at depth at most $n+m-2$, and thus exactly one of
(i) $n/m \leq a/b$ and (ii) $c/d \leq n/m$ must hold. As $a/b<\alpha < c/d$,
we have $D^\alpha(n/m)=0$ iff $n/m \leq a/b$.

\begin{lemma}\label{lem:hur_to_ded}
Let $H: \mathbb{N} \longrightarrow \{0,1\}^*$
be the Hurwitz characteristic of a irrational  $\alpha \in (0,1)$.
There is a parameterized function-oracle Turing machine $M$  such that 
\begin{itemize}
    \item $\Phi_\paramorac{M}{H}
: \mathbb{Q} \longrightarrow \{0,1\}$ is the Dedekind cut of $\alpha$
\item $M^H$ on input $n/m$ runs in time $\poly{\max\{n,m\}}$ and uses exactly one oracle call of input size at most $O(\log\max\{n,m\})$.
\end{itemize}
\end{lemma}

\begin{proof}
If $n/m\leq 0$, then $M$ outputs 0. If $n/m\geq 1$, then $M$ outputs 1.
Let $0< n/m < 1$. We assume that $n/m$ is in its lowest terms.
Then $M$ computes $H(n+m)$  using $\polylog{\max\{n,m\}}$ operations to compute the binary representation
of $n+m$, and then performing a single oracle call; note that the result of the oracle $H(n+m)$ is a bit string of length exactly $n+m \le \poly{\max\{n,m\}}$.
$M$ then computes $\fareypairtree(H(n+m))$ 
to obtain a Farey pair $(a/b,c/d)$ such that
$a/b < \alpha < c/d$. If $n/m\leq a/b$, then $M$ outputs 0, otherwise, $M$ outputs 1.

By Proposition \ref{Hurwitzburdeslankesig}, $M$ can find $(a/b,c/d)$ 
in time $$\poly{1 + \vert H(n + m) \vert} = \poly{\max\{n,m\}}$$ 
and whether $n/m\leq a/b$ holds can be tested
in time $O(\log\max\{a,b,n,m\})$. 
It is an easy induction on the depth $d$ to see that a numerator or denominator in any fraction occurring in a Farey pair at depth $d$ in $\fareypairtree$ is at most
$2^d$. Hence, $\max\{a,b,n,m\} \leq 2^{n+m}$, and the test can thus be performed in time
$O(n+m) = O(\max\{n,m\})$. Thus, $M$ needs a total time of $\poly{\max\{n,m\}}$.
\end{proof}

Our algorithm for converting a Dedekind $D^\alpha$ cut 
into a Hurwitz characteristic
$H^\alpha$   is not very surprising. The value of $H^\alpha(n)$ 
is path of length $n$ in the Farey pair tree $\fareypairtree$ 
where every interval along the path
contains $\alpha$. We can easily compute such a path when we have access 
to $D^\alpha$.

\begin{lemma}\label{lem:ded_to_hur}
Let
$D: \mathbb{Q} \longrightarrow \{0,1\}$
be the Dedekind cut of an irrational  $\alpha \in (0,1)$.
There is a parameterized function-oracle Turing machine $M$  such that
\begin{itemize}
    \item $\Phi_\paramorac{M}{D}
: \mathbb{N} \longrightarrow \{0,1\}^*$ is the Hurwitz characteristic of $\alpha$
\item $M^D$ on input $n$ runs in time $\poly{n}$ and uses exactly $n$ oracle calls, each
of input size at most $O(n)$.
\end{itemize}
\end{lemma}

\begin{proof}
On input $n$, $M$ constructs a path of length $n$ in $\fareypairtree$. 
$M$ does this
by starting at $i=0$ and incrementing $i$, maintaining a \emph{current} Farey pair $(a_i/b_i,c_i/d_i)$ such that $\alpha \in (a_i/b_i;c_i/d_i)$ for $i = 0,\ldots,n$ as the mediant of $(a_i/b_i,c_i/d_i)$
gives rise to the two children $p_L = (a_i/b_i,(a_i+c_i)/(b_i+d_i))$ and $p_R = ((a_i+c_i)/(b_i+d_i),c_i/d_i)$ of $(a_i/b_i,c_i/d_i)$ in $\fareypairtree$. Because $\alpha$ is irrational, it must be in exactly one of the open intervals
$(a_i/b_i,(a_i+c_i)/(b_i+d_i))$ and $((a_i+c_i)/(b_i+d_i),c_i/d_i)$, and thus $(a_{i+1}/b_{i+1},c_{i+1}/d_{i+1})$ must be either $p_L$ or $p_R$.
Clearly, $\alpha \in (a_i/b_i,(a_i+c_i)/(b_i+d_i))$ if{f} $D(a_i+c_i/b_i+d_i) = 1$ if{f}
the $i$th bit of $H(n)$ is $0$.
Hence, $M$ starts with $(a_0/b_0,c_0/d_0) = (0/1,1/1)$, and constructs the $n$ intervals
$(a_i/b_i,c_i/d_i)$ for $i = 1,\ldots,n$ by computing the mediant and querying $D$ in each step. Observe
that the query in step $i$ is the (binary representation of the) mediant of a Farey pair at depth
$i-1$, thus its denominator is bounded above by $2^{i}$ and its binary representation uses at most $O(\log 2^i) = O(i)$ bits.
%mediant and the denominator of the mediant
%of any Farey pair at depth $d$ is bounded above by $2^{i}$, the binary representation at step $i$ uses at most $O(\log 2^i) = O(i)$ bits.

As the numerators and denominators at depth $i$ in $\fareypairtree$ are of size at most
$2^i$ (hence representable by $i$ bits), computing the mediant at step $i$ can be done in time at most $O(i) = O(n)$ by two standard schoolbook additions, and the step $i$ contains exactly one query to $D$.
Hence, the total time needed for $M$ to construct $H(n)$ is at most
$O(n \poly{n}) = \poly{n}$, with
exactly $n$ oracle calls, each of size at most $O(\log 2^n) = O(n)$. 
\end{proof}

\subsection{Summary.}
\label{dekindoppsum}

We can now give a summary of our results on the complexity of conversions
among representations subrecursively equivalent to the representation by
Dedekind cuts.
%, along the lines
%we gave a summary  at the end of
%Section \ref{rosating} of our corresponding results
%with respect to representations subrecursively equivalent to the representation 
%by Cauchy sequences.

\begin{theorem}
Consider the representations by (1)  Dedekind cuts,
(2)  general base expansions, (3)  Hurwitz characteristics and (4)  Beatty sequences,  and let $R_1$ and $R_2$ be
any two of these four representations. Then, for  an arbitrary time-bound $t$,
we have 
$$O(t)_{R_2} \subseteq O(\poly{t(2^{n})})_{R_1} \; .$$
\end{theorem}

%% file: bestapproxequiv.tex
\section{Representations Equivalent to Best Approximations}

\label{skitteskebest}

\subsection{Best approximations.}

The representation of real numbers by left  (or right) best approximations might not
not be very well known, but it is a natural representation which
is intuitively easy to understand, and the next definition should not require any
explanations.

\begin{definition}
Let $\alpha$ be an irrational number in the interval $(0,1)$, and
let $r= a/b$ where $a,b$ are relatively prime natural numbers.

The rational  $r$ is a {\em left best approximant} of  $\alpha$ if 
we have $c/d\leq a/b<\alpha$ or $\alpha< c/d$ for any natural numbers $c,d$ where $d\leq b$.
The rational  $r$ is a {\em right best approximant} of  $\alpha$ if 
we have $\alpha < a/b \leq c/d$ or $c/d<\alpha$ for any natural numbers $c,d$ where $d\leq b$.

A {\em left best approximation} of $\alpha$ is a sequence of rationals $\{r_i\}_{i\in \nat}$
such that
$$
0 \; = \; r_0 \; < \; r_1 \; < \; r_2 \;   < \; \ldots 
$$
and each $r_i$ is a left best approximant to $\alpha$.
A {\em right best approximation} of $\alpha$ is a sequence of rationals $\{r_i\}_{i\in \nat}$
such that
$$
1 \; = \; r_0 \; > \; r_1 \; > \; r_2 \;   > \; \ldots 
$$
and each $r_i$ is a right best approximant to $\alpha$.
\qed
\end{definition}

In this section we will study a number of representations subrecursively equivalent
to the representation by left best approximations and a number of representations equivalent
to the representation by right best approximations. These two equivalence classes are incomparable to each other, that is, a representation
in one of the classes cannot be subrecursively converted to 
a representation in the other class.
We will explain why towards the end of this section.

\begin{definition}
A left best approximation   $\{r_i\}_{i\in \nat}$ of $\alpha$ is {\em complete} if  every  left best approximant of $\alpha$
occurs in the sequence $\{r_i\}_{i\in \nat}$.
A right best approximation   $\{r_i\}_{i\in \nat}$ of $\alpha$ is {\em complete} if  every  right best approximant of $\alpha$
occurs in the sequence. % $\{r_i\}_{i\in \nat}$.
\qed
\end{definition}

There is a connection between complete best approximations and paths in the Farey pair tree $\fareypairtree$. Let $\alpha$ be an irrational number, and let $\sigma_0 \sigma_1 \sigma_2 \dots\in \{0,1\}^\omega$
be the unique path in $\fareypairtree$ such that  
\begin{equation}     \label{sovetlenge}
\alpha\in\fareypairtree(\sigma_0 \ldots \sigma_n) = (a_n/b_n , c_n/d_n)
 \end{equation}
holds for any $n\in\nat$. 
Then, a  fraction $p/q$ is a left best approximant to $\alpha$ if and only
if $p/q$ occurs in the sequence $a_0/b_0, a_1/b_1, a_2/b_2, \ldots$.

It is obvious that every $a_i/b_i$ in the sequence is a left best approximant
to $\alpha$ as any  fraction in the interval $(a_i/b_i, c_i/d_i)$ has denominator
strictly greater than $b_i$ (see Theorem \ref{slankemeg} ). To see that every left best approximant to $\alpha$
indeed occurs in the sequence, assume for the sake of contradiction that a left
best approximant $p/q$ is not there. Then we have
$$ 
\frac{a_i}{b_i} \; < \; \frac{p}{q} \; < \;  \frac{a_{i+1}}{b_{i+1}}  
\; = \; \frac{a_{i}+ c_{i}}{b_{i} + d_{i}} 
$$
for some $i$. Since $p/q$ is a best approximant we must have $q <  b_{i+1} = b_{i} + d_{i}$.
But  we also have $p/q \in (a_i/b_i, c_i/d_i)$, and
that contradicts Theorem \ref{slankemeg} which states that any fraction in the interval
$(a_i/b_i, c_i/d_i)$ has denominator greater than or equal to $b_{i} + d_{i}$.

By the same token, $p/q$ is a right best approximant to $\alpha$ if and only
if $p/q$ occurs in the sequence $c_0/d_0, c_1/d_1, c_2/d_2, \ldots$. Hence, we have
the next lemma.

\begin{lemma} \label{farey=approximants}
Let $\alpha$ be an irrational number such that (\ref{sovetlenge}) holds.
(i) The 
sequence $\{ a_i/b_i\}_{i\in\nat}$ contains all the left best approximants of $\alpha$
and  nothing but left best approximants of $\alpha$.
(ii) The 
sequence $\{ c_i/d_i\}_{i\in\nat}$ contains all the right best approximants of $\alpha$
and  nothing  but right best approximants of $\alpha$.
\end{lemma}

We can use $\fareypairtree$ to subrecursively compute the  complete left best approximation 
$\{a_i/b_i \}_{i\in \nat}$ of $\alpha$ from an arbitrary  left best approximation 
$\{\hat{a}_i/\hat{b}_i \}_{i\in \nat}$  of $\alpha$.
Observe that we have $a_i / b_i \leq \hat{a}_i/\hat{b}_i  < \alpha$ for any $i\in \nat$.
Thus, there will be at least $n$ left best approximants of $\alpha$ that are smaller than
or equal to $\hat{a}_n/\hat{b}_n$.
Hence we can find  $a_n/b_n$ by constructing a path $\sigma_0\ldots \sigma_{m}$ in $\fareypairtree$
such that $\hat{a}_n/\hat{b}_n$ is the left endpoint of
$\fareypairtree(\sigma_0\ldots \sigma_{m})$.
By Lemma \ref{farey=approximants}, every left best approximant to $\alpha$ 
smaller than or equal to 
$\hat{a}_n/\hat{b}_n$ will occur along the path $\sigma_0\ldots \sigma_{m}$. We have argued
that there will be at least $n$ of them, and thus we can pick the $n$th one. 
By Lemma \ref{prop:alleskalslankesig}, we have $m\leq \hat{a}_n + \hat{b}_n-2$, and thus, unbounded search 
is not required.
A symmetric algorithm will compute the  complete right best approximation 
from an arbitrary right best approximation.

\begin{lemma}
 Let $L : \nat \longrightarrow \rational$ be a 
left best approximation of an irrational
$\alpha \in (0,1)$. Assume  $L$ is computable by an $O(s)$-time Turing machine.
There is a parameterized function-oracle Turing Machine $M$ such that
\begin{itemize}
    \item $\Phi^{L}_M :\nat \longrightarrow \rational$ is the complete
left best approximation of $\alpha$
\item $M^L$ on input $n$ runs in time $\poly{2^{s(\bitlen{n})}}$ and uses
  exactly one oracle call
of input size $\bitlen{n}$.
\end{itemize}
\end{lemma}

\begin{proof}
By the algorithm shown above, noting that we trace in the Farey tree a path of length at most $\hat{a}_n + \hat{b}_n-2$, and $\bitlen{a_n/b_n}\le s(\bitlen{n})$. 
The arithmetic operations performed by the algorithm take polynomial time in the
length of the operands.
\end{proof}

Given our discussion above,
it is not very hard to see that the representation by Dedekind cuts is subrecursive 
in the representation by left best approximations and also in  the representation by right best approximations: If $\{a_i/b_i \}_{i\in \nat}$ is a left best approximation of $\alpha$,
then we have $b_n > n$, and thus also
$m/n < \alpha$ iff $m/n < a_n/b_n$.
If $\{c_i/d_i \}_{i\in \nat}$ is a right best approximation of $\alpha$, then we have $d_n> n$,
and thus also $\alpha < m/n$ iff $c_n/d_n < m/n$.

\begin{lemma} \label{ineedthislabel}
Let $L : \nat \longrightarrow \rational$ be a
left best approximation of an irrational
$\alpha \in (0,1)$. Assume  $L$ is computable by an $O(s)$-time Turing machine.
There is a parameterized function-oracle Turing Machine $M$ such that 
\begin{itemize}
    \item $\Phi^{L}_M :\rational  \longrightarrow \{0, 1 \}$ 
is the Dedekind cut of $\alpha$
\item $M^L$ on input $p/q$ runs in time $\poly{s(\bitlen{q})}$ and uses 
exactly one oracle call of  input size $\bitlen{q}$.
\end{itemize}
\end{lemma}

\begin{proof}
We ask the oracle for the $q$th best approximant $a_q/b_q$. Then we have
$p/q \le a_q/b_q < \alpha$ or $\alpha <  p/q$. Thus, using a single
oracle call and a comparison of rationals, we can decide whether $p/q<\alpha$.
\end{proof}

\subsection{Definitions.}

We will now  define and explain a few representations 
subrecursively equivalent to left, or right, best approximations.

The base-$b$ sum approximation of from below (above) of $\alpha$, denoted 
$\hat{A}^{\alpha}_b$ ($\check{A}^{\alpha}_b$), is defined and discussed in 
Section \ref{tobiassnarteksamen}, see Definition \ref{amandaborhosmeg}.
The general sum approximation from below (above) of $\alpha$ encompasses
the base-$b$ sum approximation from below (above) of $\alpha$ for any base $b$.
The formal definition follows.

\begin{definition}
The {\em general sum approximation from below} of $\alpha$ is the function 
$\hat{G}^{\alpha}:(\nat \setminus \{0,1 \})\times\nat\longrightarrow \rational$ given 
by  $\hat{G}^{\alpha}(b,n)=\hat{A}^{\alpha}_b(n)$. 
The {\em general sum approximation from above} of $\alpha$ is the function $\check{G}^{\alpha}:(\nat \setminus \{0,1 \})\times\nat\longrightarrow \rational$ given
by  $\check{G}^{\alpha}(b,n) = \check{A}^{\alpha}_b(n)$. 
\qed
\end{definition}

What we will call a {\em Baire sequence} is an infinite sequence of natural numbers.
Such a sequence  $a_0,a_1, a_2,\ldots $  represents an irrational number $\alpha$ in the interval $(0,1)$.
We split the interval $(0,1)$ into infinitely many open subintervals with rational endpoints.
Specifically, we use the splitting
$$
(\ 0/1 \; , \; 1/2 \ )\; ( \ 1/2 \; , \; 2/3 \ ) \; ( \ 2/3 \; , \;3/4 \ ) 
\; \ldots \; ( \ n/(n+1) \; , \; (n+1)/(n+2)  \ ) \ldots      \; .
$$
The first number of the sequence $a_0$ tells us in which of these intervals we find $\alpha$.
Thus if $a_0=17$, we find $\alpha$ in the interval $(17/18,18/19)$. Then we split the interval $(17/18,18/19)$ 
in a similar way.  The second number of the sequence $a_1$
tells us in which of these intervals we find $\alpha$, and thus we proceed.  

In general, in order to split the interval $(q,r)$, we need a strictly increasing
sequence of rationals $s_0, s_1, s_2 \ldots$
such that $s_0=q$ and $\lim_i s_i = r$. We will use the splitting $s_i = (a + ic)/ (b + id)$
where $a,b$ are (the unique) relatively prime natural numbers such that $q= a/b$
and  $c,d$ are (the unique) relatively prime natural numbers such that $r= c/d$  (let $0= 0/1$ and $1= 1/1$). 
This particular splitting ensures that every interval induced by a Baire
sequence can be found in the Farey pair tree $\fareypairtree$.

We will say that the Baire sequences explained above are {\em standard}.
 The standard Baire sequence of the irrational number $\alpha$  will lexicographically precede 
standard Baire sequence of the irrational number $\beta$ iff $\alpha<\beta$.
We will also work with what we  call {\em dual} Baire sequences.
The dual sequence of $\alpha$ will lexicographically precede  the dual sequence of $\beta$ iff $\alpha > \beta$.
We get the dual sequences by using decreasing sequences of rationals to split intervals, e.g., the interval $(0,1)$ will
be split into the intervals  
$$ (\ 1/1 \; , \; 1/2 \ )\;   (\ 1/2 \; , \; 1/3 \ )  \; ( \ 1/3 \; , \; 1/4 \ )
\; \ldots \; ( \ 1/n \; , \; 1/(n+1)  \ ) \;
\ldots    \;  \;       \; .
$$

\begin{definition} \label{juleto}
Let $f:\nat \longrightarrow \nat$ be any function, and let $n\in\nat$. We define the interval $I^n_f$ by
 $I^0_f = (0/1, 1/1)$ and 
$$
I^{n+1}_f \;\; = \;\; \left( \ \frac{a+ f(n)c}{b+ f(n)d}\; , \;  \frac{a+ f(n)c + c}{b+ f(n)d + d}  \ \right)
$$
if $I^n_f = (a/b,c/d)$.
We define the interval $J^n_f$ by $J^0_f = (0/1, 1/1)$ and
$$
J^{n+1}_f \;\; = \;\; \left( \ \frac{a+ f(n)a+ c}{b+ f(n)b+ d}\; , \;  \frac{f(n)a+ c}{f(n)b+ d}  \ \right)
$$
if $J^n_f = (a/b,c/d)$.
The  function $B:\nat \longrightarrow \nat$ is the {\em standard Baire representation} of the irrational number $\alpha\in (0,1)$
if we have $\alpha\in I^n_B$ for every $n$.
The  function $A:\nat \longrightarrow \nat$ is the {\em dual Baire representation} of the irrational number $\alpha\in (0,1)$
if we have $\alpha\in J^n_A$ for every $n$.
\qed
\end{definition}

 Unit fractions, that is, fractions with nominator 1,
 were studied in the ancient Egypt, see e.g.~\cite{RMP}, and are also known as 
 Egyptian fractions.
In the literature,  an {\em Egyptian fraction expansion} may refer to any representation of 
a number as a sum of  fractions with nominator 1. 
The definition we give below ensures that any irrational number in the interval $(0,1)$
has a unique 
Egyptian fraction expansion, see Cohen \cite{Cohen1973}.

\begin{definition} \label{korkzero}
The function
$E^\alpha: \nat^{+} \longrightarrow \nat$ is the \emph{Egyptian fraction expansion} for  $\alpha$ if 
$$
\alpha = \sum_{i = 1}^\infty \left( \prod_{j=1}^i E(j) \right)^{-1}
$$
and $E(i)\leq E(i+1)$ (for all $i\in \nat^{+}$).
\qed
\end{definition}

We have e.g.\
$$
\sqrt{2} - 1 = \frac{1}{3} + \frac{1}{3 \cdot 5} + \frac{1}{3 \cdot 5 \cdot 5}
+ \frac{1}{3 \cdot 5 \cdot 5 \cdot 16} + \frac{1}{3 \cdot 5 \cdot 5 \cdot 16 \cdot 18} + \ldots 
$$
and this is the unique Egyptian fraction expansion of $\sqrt{2} - 1$. 
Another possible representation of irrationals based on Egyptian fractions
 is related to left best approximation, see Beck et al.~\cite{Beck1969}. In fact, from 
Theorem~\ref{slankemeg} and the relation of left best approximations to the Farey pair tree, 
(Lemma~\ref{farey=approximants}) it is easy to deduce that the difference between consecutive
fractions in a complete left best approximation is a unit fraction. Thus a 
complete left best approximation of $\alpha$ induces a series of unit fractions that adds up to
$\alpha$.

\begin{definition}
A sequence $\{q_i\}$ of positive unit fractions is the Farey-Egyptian expansion of $\alpha\in (0,1)$
if $\sum_{i=1}^\infty q_i = \alpha$, and the sequence $(\sum_{i=1}^j q_i)_j$ is a complete left best approximation
of $\alpha$.
\qed
\end{definition}

For example
$$
\sqrt{2} - 1 = \frac{0}{1} + \frac{1}{3} + \frac{1}{15} + \frac{1}{85}
+ \frac{1}{493} + \ldots 
$$
is  Farey-Egyptian expansion, associated with the complete left best approximation  $0/1,\;1/3,\;2/5,\;7/17\dots$.

The representation by
Farey-Egyptian expansions is closely related to the representation by 
complete left best approximations, and it is easy to convert 
Farey-Egyptian expansion into a complete left best approximation, and vice versa.

General sum approximations (from above and below) were introduced in  Kristiansen \cite{ciejouren} and studied further in Georgiev et al.~\cite{apalilf}.
 Left and right best
approximations are studied in 
\cite{apalilf}, and standard and dual Baire sequences are studied in
Kristiansen \cite{bairelars}. The computational complexity of representations by
 Egyptian fractions expansions is studied for the first time in this paper.

\subsection{Conversion between general sum approximations and best approximations.}

Let $\{a_i/b_i\}_{i\in \nat}$  be a left best approximation of $\beta$.
We will give an algorithm for computing the general sum approximation from below of $\beta$, that is $\hat{G}^\beta$,
using $\{a_i/b_i\}_{i\in \nat}$ as an oracle.  The $i$th iteration of the algorithm generates
$\hat{G}^\beta(b,i)$.
Having computed
$\hat{G}^\beta(b,i)$ for all $i < n$, the algorithm
will also have computed the sum
$$
\frac{c}{d} \; = \; \sum_{i=1}^{n-1} \hat{G}^\beta(b,i) \;.
$$
It then computes $\hat{G}^\beta(b,n)$ by executing the following
instructions.
\begin{itemize}
\item Step 1: 
Ask the oracle for the value of $a_d/b_d$.
Let $c'/d' = a_d/b_d$. Note that $d$ is a power of $b$; and that we have $d'>d$,
and  $c/d < c'/d' < \beta$.
\item Step 2: Compute $m:= \lceil \log_b d'd \rceil$.
\emph{Comment\/}: We have $c'/d' -c/d > 0$, and thus also
$$
\frac{c'}{d'} - \frac{c}{d} = \frac{c'd - cd'}{d'd}
\;\;  \geq \;\; \frac{1}{d'd} \; .
$$
It follows that $c/d +  1/b^{m} \leq  c'/d' < \beta$.

\item Step 3: Ask the oracle for the value of $a_{b^m}/b_{b^m}$ and compute (as further explained below) the least $ k \le m$ and base-$b$ digit $\xd>0$ so that
\begin{equation}
    \label{larsdidthi}
\frac{c}{d} + \frac{\xd}{b^k} \;\; < \;\; \beta \;\; < \;\;  \frac{c}{d} + \frac{\xd}{b^k} + \frac{1}{b^k}\; .
\end{equation}
Give the output $\xd/b^k$.
\end{itemize}
\emph{Comments on Step 3}:   
We have $i< b_i$ for any left best approximant $a_i/b_i$. 
 Hence, for any  fraction $p/{b^k}$, with $k\le m$, we   have $p/{b^k} < \beta$ iff 
$p/{b^k} \le a_{b^m}/b_{b^m}$, and \eqref{larsdidthi} is equivalent, for such $k$, to
\begin{equation} \label{eq-k-and-D}
\frac{c}{d} + \frac{\xd}{b^k} \;\; \leq \;\; \frac{a_{b^m}}{b_{b^m}} \;\; < \;\;  \frac{c}{d} + \frac{\xd}{b^k} + \frac{1}{b^k}\; .
\end{equation}
The least $k$ that satisfies (\ref{eq-k-and-D}) for some $\xd$ is the least $k$ that satisfies
\begin{equation*}
\frac{c}{d} + \frac{1}{b^k} \;\; \leq \;\; \frac{a_{b^m}}{b_{b^m}}
\end{equation*}
In order to find that $k$, we rewrite the  inequality as
$$
\frac{1}{b^k} \;\; \leq \;\; \frac{a_{b^m}}{b_{b^m}} - \frac{c}{d} 
$$
which again can be rewritten as
$${b^k} \;\; \ge \;\; \frac{d b_{b^m}}{d a_{b^m} - c{b_{b^m}}} \; .$$
Hence, the desired $k$ is
$$
k = \left\lceil \log_b \frac{d\cdot b_{b^m}}{d\cdot a_{b^m} - c\cdot b_{b^m}} \right\rceil \;.
$$
Having computed $k$, we look for a value of $\xd$ such that \eqref{eq-k-and-D} holds, and
it should be clear that 
$$
\xd = \left\lfloor \left(  \frac{a_{b^m}}{b_{b^m}} - \frac{c}{d} \right)/ b^k \right\rfloor \;. 
$$

The correctness of the algorithm follows straightforwardly from the comments on Steps 1--3 and the definition
of a general sum approximation from below.
A right best approximation can be converted into a general sum approximation from above by a symmetric algorithm.

\begin{lemma}
Let $L : \nat \longrightarrow \rational$ is 
left best approximation of an irrational
$\alpha \in (0,1)$. Assume $L$ is computable within a time bound $s$,
and let $f(b,n) = (\lambda x.(\bitlen{b} + 2s(x)))^{(n)}(1)$.
There is a parameterized function-oracle Turing Machine $M$ such that 
\begin{itemize}
    \item $\Phi^{L}_M :(\nat\setminus \{0, 1 \} )\times \nat \longrightarrow \rational$ is the general sum approximation from below of $\alpha$
    \item $M^L$ on input $b,n$ runs in time $\poly{s(f(b,n))}$ and uses 
    at most $2n$ oracle calls, each of input size at most $f(b,n)$.
\end{itemize}
\end{lemma}

\begin{proof}The proof is, of course, inductive. 
We claim that the bit-lengths of $d, b^m$ are
bounded by $f(b,n)$. Consequently, by assumption and the rule that the running time of
a machine bounds the size of its output, we also have that the bit-lengths of
the denominators $b_d,b_{b^m}$  are bounded
by $s(f(b,n))$. 
The induction step works as follows: 
The algorithm first sets $d' = b_d$.  We inductively assume the bound $s(f(b,n-1))$ for the bit-length of the denominators of each of 
$\hat{G}^\alpha(b,i)$, $i<n$. The common denominator $d$ is the denominator of the last term, and so its bit-length is bounded
by $s(f(b,n-1))$.  
Applying our assumption of the time bound $s$ for the
computation of $b_d$, 
we have $\bitlen{d'}\le s(\bitlen{d}) \le s(f(b,n-1))$.
We bound the size of the next denominator, $d_{b^m}$ as follows. First,
\begin{multline*}
    \bitlen{b^m} \; = \;  \bitlen{b^{\lceil log_b(dd')\rceil}} 
   \; \le \; \bitlen{bdd'} \\ \; \le \; \bitlen{b}+\bitlen{d}+\bitlen{d'} 
  \;  \le  \; \bitlen{b}+s(f(b,n-1))+s(f(b,n-1))  \; \le \; \\ \bitlen{b} + 2s(f(b,n-1)) 
  \; = \; f(b,n)
\end{multline*}
and once more, by the assumption of the time bound $s$, we have 
\[
\bitlen{d_{b^m}} \; \le \;  s(\bitlen{b^m}) 
\;  \le  \; s(f(b,n))\; .
\]
The execution time is polynomial in the size of the numbers manipulated, hence polynomial in $s(f(b,n))$.
\end{proof}

Next we give an algorithm for computing a complete 
left best approximation  $\{a_i/b_i\}_{i\in \nat}$ of $\beta$
which uses $\hat{G}^\beta$ (the general sum approximation from below of $\beta$) as an oracle.
Observe that the oracle makes it easy to compute the Dedekind cut $D^\beta$ of $\beta$: 
We have $D^\beta(c/d) = 0$ iff  $\hat{G}^\beta(d,1)\geq c/d$.

To compute $a_0/b_0$ is trivial since we have $a_0/b_0 = 0/1$ by convention. 
In order to compute 
$a_1/b_1$ the algorithm  ask the the oracle for the value of $\hat{G}^\beta(2,1)$. We
have $\hat{G}^\beta(2,1)= 1/M$ for some  $M$.
Then the  algorithm  uses the Dedekind cut of $\beta$ to search for $N\leq M$ such that 
$1/N < \beta < 1/(N-1)$ and set $a_1/b_1$ to $1/N$. (Such an $N$ will exist since $1/M<\beta$.
Note that $2/N \ge 1/(N-1) > \beta$, and hence, the algorithm computes $a_1/b_1$ correctly.)

When $n\geq 1$, the algorithm computes $a_{n+1}/b_{n+1}$
by the following procedure:
\begin{itemize}
    \item  Let 
$$
\frac{c'}{d'} \;\; = \;\; \frac{a_{n}}{b_{n}} \; + \;  \hat{G}^\beta(b_n,2)\;
$$ 
with relatively prime $c'$ and $d'$.
\item Use the Dedekind cut of $\beta$ to search for $c$ and the smallest $d$ such that $b_n < d \leq d'$ and $a_n/b_n < c/d < \beta < (c+1)/d$.
Let $a_{n+1}/b_{n+1}= c/d$.
\end{itemize}
In order to see that the algorithm is correct, observe that $a_{n}/b_{n}= \hat{G}^\beta(b_n,1)$. Hence, we have $a_n/b_n < c'/d' < \beta < (c'+1)/d'$,
and $c'/d'$ will be a left best approximant unless there exists a fraction $c/d$  such that $d< d'$ and $a_n/b_n < c/d < \beta < (c+1)/d$.  We are looking for the least such $d$, therefore we find a complete best approximation.

A general sum approximation from above can be converted into a right best approximation by a symmetric algorithm.

\begin{lemma}
Let $G :(\nat\setminus \{0, 1 \} )\times \nat \longrightarrow \rational$
be the general sum approximation from below of an irrational $\beta \in (0,1)$.
Assume that $G$ is computable within a time bound $s$, and let 
$f(n) = (\lambda x.2s(x))^{(n)}(1)$.
There is a parameterized function-oracle Turing Machine $M$ such that 
\begin{itemize}
    \item $\Phi^{G}_M :\nat \longrightarrow \rational$ is the
complete left best approximation  of $\beta$
\item $M$ on input $n$ runs in time $\poly{2^{f(n)}}$ and uses  at most $O(2^{f(n)})$
oracle calls, each of input size at most $f(n)$.
\end{itemize}
\end{lemma}

\begin{proof}
First we explain the implementation of the search for $c,d$ in the inductive step. 
The condition we test is
$b_n < d \leq d'$ and $a_n/b_n \leq c/d < \beta < (c+1)/d$.  
The oracle can be used to replace the
second conjunct by $a_n/b_n \leq c/d = G(d,1)$.
It follows that to perform the search, we first ask for $G(b_n,2)$, then for
$G(d,1)$ for successive values of $d$, from $b_n$ upwards, 
and at most up to $G(b_n,2)$.  
We conclude that the largest number involved in computing
$a_{n+1}/b_{n+1}$ is $G(b_{n},2)$.

Now, consider the whole process of computing $a_0/b_0,\dots,a_n/b_n$:
in the process of computing $a_{i+1}/b_{i+1}$ we make a single oracle call
to bound the search, and then we search from $b_i+1$ up to $b_{i+1}$ using
a single oracle call for each test. It is easy to see that the total number
of oracle calls is $n + b_n$.  To get a bound in terms of the input to the
algorithm, we note that $d'$ in the induction step is at most the denominator of 
$G(b_{n-1},2)$.
By assumption, $G(b_{n-1},2)$ occupies at most $s(2\bitlen{b_{n-1}})$ bits since the representation of a fraction is at most twice the size of
its denominator. Hence
$$
  b_n < 2^{2s(\bitlen{b_{n-1}})}
$$ 
so $\bitlen{b_n} \le 2s(\bitlen{b_{n-1}})$.  
This gives the bound $f(n)$ on the bit-length of the largest number
involved in the computation. The value of this number,
bounded by $2^{f(n)}$, bounds the number of \emph{steps} in the computation, and hence the execution time up to a polynomial.
\end{proof}

\subsection{Conversion between Baire sequences and best approximations.}

\begin{lemma} \label{lysnisse}
We have
\begin{align} \tag{i}
I^{n+1}_f = \fareypairtree(1^{f(0)}0 1^{f(1)}0 \ldots 1^{f(n)}0)    
\end{align}
and
\begin{align} \tag{ii}
J^{n+1}_f =   \fareypairtree( 0^{f(0)}1 0^{f(1)}1 \ldots 0^{f(n)}1 ) \; .
\end{align}
\end{lemma}

\begin{proof} 
We prove (i). The proof of (ii) is symmetric. 
Let $\sigma = 1^{f(0)}0 1^{f(1)}0 \ldots 1^{f(n-1)}0$. Observe that we have  $\fareypairtree (\sigma) = (0/1,1/1) =  I^0_f$ when $\sigma$ is the empty
sequence. Assume  that $\fareypairtree (\sigma) = I^n_f = (a/b, c/d)$.
We need to prove that 
\begin{align} \label{eee}
\fareypairtree (\sigma 1^{f(n)}0) \; = \; I^{n+1}_f \; . 
\end{align}
Let $k= f(n)$.
We prove (\ref{eee}) by a secondary induction on $k$.

Assume $k=0$. By Definition \ref{def:farey}, we have 
$$\fareypairtree ( \sigma 1^{f(n)}0 )\; = \;  \fareypairtree ( \sigma 1^{0}0 ) \; =  \; \fareypairtree ( \sigma 0 )\; = \; ( \ a/b \, , \, (a+c)/(b+d) \ )\; .$$
By Definition \ref{juleto}, we have $$I^{n+1}_f \; = \; ( \ (a+ kc)/(b+kd) \, , \, (  a+ kc+c)/(b+kd+d) \ ) 
\; = \; ( \ a/b \, , \, (a+c)/(b+d) \ )\; .$$
Thus, (\ref{eee}) holds when $f(n)=0$. Now, assume by induction hypothesis that
\begin{align} \label{fff}
\fareypairtree ( \sigma 1^{k}0 ) \; = \; \left( \ \frac{a+kc}{b+kd} \, , \, \frac{a+kc+c}{b+kd+d} \ \right) \; . 
\end{align}
Observe that the right hand side of (\ref{fff}) is the definition of $I^{n+1}_f$ with $k$ for $f(n)$.
Now, by  (\ref{fff}) and Definition \ref{def:farey}, we have 
\begin{align} \label{femticl}
\fareypairtree  ( \sigma 1^{k} ) \; = \; \left( \ \frac{a+kc}{b+kd} \, , \, \frac{c}{d}      \ \right) \; . 
\end{align}
Furthermore, by  (\ref{femticl}) and Definition \ref{def:farey}, we have 
\begin{align} \label{hundrecl}
\fareypairtree ( \sigma 1^{k+1} ) \; = \; \left( \ \frac{a+kc + c}{b+kd + d} \, , \, \frac{c}{d} \ \right) 
\; = \; \left( \ \frac{a+(k+1)c}{b+(k+1)d} \, , \, \frac{c}{d} \ \right) 
\end{align}
and by   (\ref{hundrecl}) and Definition \ref{def:farey}, we have 
\begin{align} \label{ggg}
\fareypairtree ( \sigma 1^{k+1}0 ) \; = \; \left( \ \frac{a+(k+1)c}{b+(k+1)d} \, , \, \frac{a+(k+1)c + c}{b+(k+1)d +d} \ \right) \; . 
\end{align}
Observe that the right hand side of (\ref{ggg}) is the definition of $I^{n+1}_f$ with $k+1$ for $f(n)$.
This proves that (\ref{eee}) holds.
\end{proof}

Given the lemma above it is easy to see how we can convert a standard Baire sequence $B$ into 
a complete right best approximation. We use $B$ to compute an interval $I$ in the Farey pair tree
such that $I = \fareypairtree(1^{B(0)}0 1^{B(1)}0 \ldots 1^{B(n)}0)$. By Lemma \ref{lysnisse} (i), we have
$I = I^{n+1}_B$. It follows that the right endpoint of $I$
is $n$th approximant in the complete right best approximation of $\alpha$ (see Lemma \ref{farey=approximants}).

By Lemma \ref{lysnisse} (ii), we have a symmetric algorithm for converting a dual Baire sequence into a left best approximation.

\begin{lemma}
Let  $B : \nat \longrightarrow \nat$
be the standard Baire sequence of an irrational $\alpha \in (0,1)$.
There is a parameterized function-oracle Turing Machine $M$ such that
\begin{itemize}
    \item $\Phi^{B}_M :\nat \longrightarrow \rational$ is the complete right best approximation  of $\alpha$
    \item  $M^B$ On input $ n \in \nat$ runs in time $\poly{n+\sum_{i=0}^n B(i)}$ and
    uses $n$ oracle calls, each  of input size $O(\log n)$
\end{itemize}
moreover, if $B$ is computable within the time bound $s$, then $\Phi^{B}_M$ runs in
time $\poly{n 2^{s(\bitlen{n})}}$.
\end{lemma}

\begin{proof}
We have to compute $I = \fareypairtree(1^{B(0)}0 1^{B(1)}0 \ldots 1^{B(n)}0)$.  By Proposition~\ref{Hurwitzburdeslankesig} this is polynomial in $n+\sum_{i=0}^n B(i)$.
By the standard argument, this shows that if $B$ is computable within the time bound $s$, then $\Phi^{B}_M$ runs in
time $\poly{n 2^{s(\bitlen{n})}}$.
\end{proof}

Lemma \ref{lysnisse} (i) also yields an algorithm for converting complete right best approximations into standard
Baire sequences: Given the complete right best approximation $\{a_i/b_i\}_{i\in \nat}$ of $\alpha$, we can compute
a (unique) string of the form $1^{k_0}0 1^{k_1}0\ldots 1^{k_n}0$ such that the right endpoint of the interval 
$\fareypairtree(1^{k_0}0 1^{k_1}0\ldots 1^{k_i}0)$ equals  $\{a_{i+1}/b_{i+1}\}$ (for all $i\leq n$).
By Lemma \ref{lysnisse}, we have $B(i)= k_i$ where $B$ is the standard  Baire sequence of $\alpha$.
Lemma \ref{lysnisse} (ii)  yields an algorithm for converting complete left best approximations into dual
Baire sequences.

\begin{lemma}
Let $R :  \nat \longrightarrow \rational$
be the complete right best approximation of an irrational $\alpha \in (0,1)$.
There is a parameterized function-oracle Turing Machine $M$ such that
\begin{itemize}
    \item $\Phi^{R}_M :\nat \longrightarrow \nat$ is  the standard Baire sequence of $\alpha$
    \item $M^R$ on input $n$ runs in time $\poly{2^{\bitlen{R(n)}}}$ and uses $n$ oracle
 calls, each of input size $O(\log n)$
\end{itemize}
moreover, if $R$ is computable within the time bound $s$, then $\Phi^{R}_M$ runs in
time $\poly{2^{s(\bitlen{n})}}$.
\end{lemma}

\begin{proof}
The machine has to trace a path in $\fareypairtree$ from the root up to the first occurrence of $R(n)=a_n/b_n$.
By Lemma~\ref{prop:alleskalslankesig}, this happens at most at depth $a_n+b_n-2 < 2^{\bitlen{R(n)}}$.
The branch taken at each level is dictated by the corresponding best approximant, so the machine has to compute
$R(0)$ through $R(n)$. The work at each level is polynomial in the depth and the size of the numbers involved (see Lemma~\ref{fareysizelemma}), which are all
polynomial in $a_n+b_n$, hence in $2^{\bitlen{R(n)}}$.

By the standard argument, this shows that if $R$ is computable within the time bound $s$, then $\Phi^{R}_M$ runs in
time $\poly{2^{s(\bitlen{n})}}$.
\end{proof}

\subsection{From general sum approximation from below to Egyptian fraction expansions.}

Let $\hat{G}^{\alpha}:(\nat \setminus \{0,1 \})\times\nat\longrightarrow \rational$  be the
general sum approximation from below of
$\alpha \in (0,1)$. We show how to convert $\hat{G}^{\alpha}$ to the 
Egyptian fraction expansion $E^\alpha$. 

The algorithm works recursively. For the base case, assume $\hat{G}^\alpha(2,1) = a/b$
(where $b$ is some power of 2). Then we search for the least $2\le d\le b$ such that the denominator of $\hat{G}^\alpha(d,1)$ equals $d$. 
This means that $1/d < \alpha < 1/(d-1)$, and we set $E^\alpha(1) = d$.
In the general case,
first the algorithm computes  
$E^\alpha(i)$  for all $i\leq n$,  then the algorithm computes $E^\alpha(n+1)$ as follows:
\begin{itemize}
    \item Step 1: 
    Compute the sum of the first $n$ terms in the Egyptian fraction expansion
$$
\frac{c}{d} \; = \; \sum_{i=1}^n \left( \prod_{j=1}^i E^\alpha(i) \right)^{-1}  .
$$
\item Step 2:
Let $c'/d'  = \hat{G}^\alpha(d, 1) + \hat{G}^\alpha(d, 2)$. 
Now $d'$ is some power of $d$, and $c/d < c'/d' < \alpha$.
\item Step 3:
Search for the least $m$ such that $c/d + (dm)^{-1} < \alpha$.
This search can be performed by the Dedekind cut of $\alpha$, in turn simulated using
$\hat{G}^\alpha$ as already shown. The search is bounded since $E(n)\le m \le d'/d$.
Return $m$.
\end{itemize}
Except for the effort to bound the searches, our algorithm is the natural greedy algorithm
which is known to compute the Egyptian fraction expansion, see Cohen~\cite{Cohen1973}.

\begin{lemma} 
Let
$G:(\nat \setminus \{0,1 \})\times\nat\longrightarrow \rational$
be the general base approximation from below of an irrational number $\alpha \in (0,1)$.
Assume that $G$ is computable within a time bound $s$, and let
$f(n) = (\lambda x. {s(2x)})^{(n)}(2)$.
There is a parameterized function-oracle Turing machine $M$ such that
\begin{itemize}
    \item $\Phi^{G}_M :\nat^+ \longrightarrow \nat$ is the Egyptian fraction expansion
of $\alpha$
\item $M^G$ on input $n$ runs in time $\poly{2^{f(n)}}$ and uses at most $2^{f(n)+1}$
oracle calls, each of input size at most $2f(n)$.
\end{itemize}
\end{lemma}

\begin{proof}
The crucial quantity that determines the complexity of this algorithm is the size 
(bit-length) of the
largest denominator encountered, denoted $d'$ in the above algorithm. We bound 
$d'$ by $f(n)$ as follows.
In the general case, $d'$ is obtained by calling 
$G(d,2)$ where $d$ comes from the previous iteration, so we assume for induction that
$\bitlen{d} \le f(n-1)$. The input $d,2$ is represented in twice the bit-length
of $d$ (assuming that a pair of integers is represented by zipping two binary numbers, as suggested 
in Section~\ref{prelims}), so the size of the oracle input is bounded by
$2f(n-1)$ and the output of the oracle has size bounded by $s(2f(n-1)) = f(n)$.
In the base case, we query the oracle for $G(2,1)$: the size of the result
is bounded by $s(4) = (\lambda x. {s(2x)})(2)$.

Thus $2f(n)$ bounds the size of any oracle input throughout the algorithm. The number of oracle calls in the general case is bounded by $2+(E(n)-E(n-1)+1)$,
where the first call is to determine $d'$, and the rest are in the search for $E(n)$.  It easily follows that the \emph{total} number of calls is bounded by
$3n+E(n)$, which we bound by $2\cdot 2^{f(n)}$ since $E(n)$ is clearly bounded by
$d' < 2^{f(n)}$, and $3n\le 2^{f(n)}$ which is easily proved by induction.
\end{proof}

\subsection{From Egyptian fraction expansions to left best approximations.}

The next lemma indicates how we can subrecursively convert the representation by Egyptian fraction expansions
into Dedekind cuts.
Here we shall use the Egyptian fraction expansion of a rational number. A rational
number may have two expansions, one finite and one infinite (e.g., $1/9$
is also $1/10+ 1/100 + 1/1000+\dots$). The algorithm we give
uses the finite one. In the proof of the next lemma, we write the expansion as an infinite sequence anyway,
for uniformity of notation; assume the missing fractions to be zero, identified
with $1/\infty = 0$ (e.g., $1/9 +0+0+\dots$).

\begin{lemma} \label{lexicograph}
Let $x,y$ be two different real numbers (possibly rational) in the interval $(0,1)$, and let $E^x$, $E^y$ be their respective Egyptian fraction
expansions. Then, $x<y$ if and only if $E^x$
precedes $E^y$ in lexicographic order.
\end{lemma}

\begin{proof}
%If $E^x(1)\ne E^y(1)$, then, by the correctness of the greedy algorithm, 
%$E^x(1) < E^y(1)$ implies $x<y$ and vice versa.  In general, 
Suppose that
$E^x(i) = E^y(i)$ for all $i<n$, and $n$ is the position of the first difference.
Then regarding the $n$th element we have the
following cases: either one of the sequences terminates, in which cases it is clear
that the number expressed by the other sequences is larger; or both continue. In the
latter case, assume w.l.o.g.~that $E^x(n) < E^y(n)$. 
Denote $d = \prod_{j=1}^{n-1} E^x(j)= \prod_{j=1}^{n-1} E^y(j)$.
Then
\begin{align*}
y-x &= \sum_{i=n}^\infty \left( \prod_{j=1}^i E^y(j) \right)^{-1} -\ 
      \sum_{i=n}^\infty \left( \prod_{j=1}^i E^x(j) \right)^{-1}  \\
 &\ge \frac{1}{d E^y(n)} - \sum_{i=n}^\infty \left( \prod_{j=1}^i E^x(j) \right)^{-1}  \\
 &\ge \frac{1}{d E^y(n)} - \frac{1}{d}\sum_{i=n}^\infty \left( \prod_{j=n}^i E^x(j)) \right)^{-1}  \\
 &\ge \frac{1}{d E^y(n)} - \frac{1}{d}\sum_{i=n}^\infty (E^x(n))^{n-i-1}
 \ge \frac{1}{d E^y(n)} - \frac{1}{d}\cdot\frac{1}{E^x(n)-1} \ge 0 \;.
\end{align*}
We conclude that $x\le y$, but they are known to differ, so $x<y$.
\end{proof}

We obtain the following algorithm for computing $D^\alpha(q)$ in $E^\alpha$:
Compute the finite Egyptian fraction expansion of $q$ and compare it lexicographically to $E^\alpha$. For computing the expansion of $q$ we use the algorithm 
of \cite{Cohen1973}, which we review below. Importantly, from the algorithm
it is easy to see that the expansion $E^{p/q}$ of a rational number $p/q$
includes at most $p$ terms and the size of each $E^{p/q}(i)$ is bounded by $q$.

\medskip

\paragraph{\sl Expansion algorithm:} We define an auxiliary function $Div(q,p)$ that for positive integers $q,p$ returns a pair of non-negative integers 
$n\le q$ and $r < p$ such that $np$ is the smallest multiple of $p$ with $np\ge q$, and $r = np -q$.
Given a fraction $p/q$, we construct sequences $n_1, n_2,\dots$ and $r_0, r_1,\dots$ by setting $r_0 = p$ and
$(n_{i+1},r_{i+1}) = Div(q,r_i)$ until we reach $r_j = 0$.  The Egyptian fraction expansion is given by
$E(i) = n_i$ for $i=1,\dots,j$.

\begin{lemma}
Let  $E :\nat^+ \longrightarrow \nat$
is the Egyptian fraction expansion of an irrational number $\alpha \in (0,1)$.
There is a parameterized function-oracle Turing Machine $M$ such that
\begin{itemize}
    \item $\Phi^{E}_M :\rational \longrightarrow \{0,1\}$ is the 
Dedekind cut of $\alpha$
\item $M^E$ on input $p/q$ runs in time $p\cdot\poly{\bitlen{q},\bitlen{E(p)}}$ 
and uses at most $p$
oracle calls, each of input size at most $1+\log p$.
\end{itemize}
\end{lemma}

\begin{proof}
The computation of the expansion of $p/q$, according to \cite{Cohen1973},
makes at most $p$ iterations, and in each iteration $O(1)$ arithmetic operations
are performed on numbers bounded by $q$. For our purpose, we compare the
$i$th number in the expansion, $E^{p/q}(i)$, with $E(i)$.  Recalling that
$E$ is a non-decreasing series, the complexity bounds stated follow.
\end{proof}

Let $E^\alpha$  be the Egyptian fraction expansion of $\alpha$.
We will give an algorithm for computing a complete left best approximation 
$L^\alpha$ of $\alpha$,  
using $E^\alpha$ as an oracle. The algorithm works recursively. 
The base case is the trivial approximation $0/1$. In the general case, for $n \ge 1$,
the algorithm first computes  
$a_{n-1}/b_{n-1} = L^\alpha(n-1)$.  Then the algorithm computes $L^\alpha(n)$ as follows:
\begin{itemize}
    \item Step 1: Compute $m$ terms of the Egyptian Fraction expansion for
the least $m$ such that
$$
\sum_{i=1}^m \left( \prod_{j=1}^i E^\alpha(j) \right)^{-1} \; > \; \frac{a_{n-1}}{b_{n-1}} 
$$
and let
$d'  = \; \prod_{i=1}^{m} E^\alpha(i)$.
\item Step 2: Using an implementation of the Dedekind cut $D^\alpha$ as described
above, search for the least $d > b_{n-1}$ such that for some $c<d$, the fraction $c/d$ is
a better left approximation to $\alpha$ than $a_{n-1}/b_{n-1}$. 
Choose the largest such $c$. Return $c/d$
as $L^\alpha(n)$. Note that $d$ will be at most $d'$.
\end{itemize}
Step~1 may seem an unbounded search. But in fact, since the
Egyptian fraction expansion of $a_{n-1}/b_{n-1}$ comprises at most $a_{n-1}$ terms,
we know, by Lemma~\ref{lexicograph}, that $m\le a_{n-1}+1$.

\begin{lemma}
Let  $E :\nat^+ \longrightarrow \nat$
be the Egyptian fraction expansion of an irrational number $\alpha \in (0,1)$.
Assume $E$ is computable within a time bound $s$, and
let $f(n) = (\lambda x. 2^{x\log x + s(\bitlen{x})})^{(n)}(1)$.
There is a parameterized function-oracle Turing Machine $M$ such that 
\begin{itemize}
    \item $\Phi^{E}_M :\nat \longrightarrow \rational$ is the complete left best 
approximation of $\alpha$
\item 
$M^E$ on input $n$ runs in time $\poly{f(n)}$ and  uses at most $f(n)$
oracle calls, each of input size at most $\bitlen{f(n)}$.
\end{itemize}
\end{lemma}

\begin{proof}
We will see that the crucial quantity in the analysis of this algorithm is $b_n$.
We claim that $b_n \le f(n)$, which follows by induction from 
$$
b_n \le 2^{b_{n-1} \log {b_{n-1}} + s(\bitlen{b_{n-1}})} \;.
$$
To justify this, note that 
\begin{equation} \label{bound-d'}
 b_n \le d' \le \prod_{i=1}^m E(i) \;.
\end{equation}
Now $m \le a_{n-1}+1 \le b_{n-1}$, as argued before the lemma. Moreover, the first $m-1$ terms $E(1),\dots,E(m-1)$
coincide with the Egyptian fraction expansion of $a_{n-1}/b_{n-1}$, which imply that their size is at most
$b_{n-1}$, hence their product at most $b_{n-1}^{b_{n-1}}$. To bound $E(m)$, we recall that
$E(m)$ is computable in time $s(\bitlen{m})$, and therefore bounded by $2^{s(\bitlen{m})} \le 2^{s(b_{n-1})}$.
We deduce from \eqref{bound-d'} that
$$
b_n \le 2^{b_{n-1} \log b_{n-1} + s(\bitlen{b_{n-1}})}
$$
This completes the induction.

Now, to bound the oracle input size we claim that $b_n$ bounds the largest value passed to the oracle
throughout the algorithm; this is easy enough to verify. The bound on the number of oracle calls also follows, since
the Turing machine can record answers from the oracle, and hence, we can avoid querying the same input twice.

To bound the execution time, we bound the time for the last iteration;  the total time is polynomial in this quantity,
since the bound on the last iteration is the largest.
In the last iteration we use the Egyptian fraction expansion of $a_{n-1}/b_{n-1}$ (known from the previous iteration)
and compare it lexicographically with the series $E(i)$ to determine $m$; this takes $O(m)=O(b_{n-1})$
arithmetic operations on numbers
bounded by $b_n$. We have proved that $b_n \le f(n)$, and the time 
for the computation is polynomial in its
bit-length, times $O(b_{n-1})$.  Thus the time required by Step~1 will be $\poly{f(n-1)\log f(n)} = \poly{f(n)}$.

Step~2 is implemented as follows: for $d = b_{n-1}+1, b_{n-1}+2, \dots$ we compute 
$c = \lceil (a_{n-1}d) /b_{n-1} \rceil$, as this gives the smallest 
nominator such that $c/d > a_{n-1}/b_{n-1}$.
We test if $c/d < \alpha$ using the Dedekind cut. The whole process is $\poly{f(n)}$ since $d\le f(n)$.
When we find the first such $c,d$, we search
for the largest $c$ such that $c/d$ is still below $\alpha$. This involves less than $d$ applications of the Dedekind cut;
again we remain within $\poly{f(n)}$ time.
\end{proof}

\subsection{Summary.}

We will now give a summary of  this section
along the lines we have given summaries of the section of Section \ref{rosating} 
(page \pageref{oppsummeringcauchy}) and Section \ref{skittenskje} (page \pageref{dekindoppsum}).

\begin{theorem} \label{langrenn}
Consider the representations by (1)  right best approximations, 
(2) complete  right best approximations
(3) Baire sequences and (4)  general sum approximations from above, and let $R_1$ and $R_2$ be
any two of these four representations. Then, for  an arbitrary time-bound $t$,
there exists a time-bound $s$ primitive recursive in $t$ such that
$O(t)_{R_2} \subseteq O(s)_{R_1}$. 
\end{theorem}
A comment meant for the readers familiar with the Grzegorcyk hierarchy:
Our results are  a bit stronger than what the theorem above asserts.
One can easily check that the time-bound 
$s$, for any $i\geq 3$, will be in the Grzegorcyk class $\mathcal{E}_{i+1}$ if the time-bound $t$ is
in the Grzegorcyk class $\mathcal{E}_{i}$. The same goes for the next theorem if we leave out
the representation by Egyptian fraction expansions.

\begin{theorem} \label{xlangrenn}
Consider the representations by (1)  left best approximations, 
(2) complete  left best approximations
(3) dual Baire sequences, (4)   general sum approximation from below and (5) Egyptian 
fraction expansions, and let $R_1$ and $R_2$ be
any two of these five representations. Then, for  an arbitrary time-bound $t$,
there exists a time-bound $s$ primitive recursive in $t$ such that
$O(t)_{R_2} \subseteq O(s)_{R_1}$. 
\end{theorem}

\section{A Little Bit on the Degrees of Representations}
\label{skalmotemartintap10}

Recall the definition of the relation $\redrel$ (Definition \ref{snartsommerigjen}), 
and recall the definition of the relations $\equiv_S$ and $\redrelstrict$ 
(Definition \ref{rodtihjornet}). Furthermore,
recall that the equivalence relation $\equiv_S$ induces a degree structure on the representations, where a degree  simply is
an $\equiv_S$-equivalence class, see page \pageref{naaskaljegsnartgaaaatrene}.
We will use the standard terminology of degree theory and say that a degree $\dga$ lies 
{\em (strictly) below} a degree $\dgb$
if we have ($R_1 \redrelstrict R_2$) $R_1 \redrel R_2$  
whenever  $R_1\in \dga$ and $R_2 \in \dgb$, moreover, we say that
 $\dgb$ lies {\em (strictly) above} $\dga$
if we have ($R_1 \redrelstrict R_2$) $R_1 \redrel R_2$ 
whenever  $R_1\in \dga$ and $R_2 \in \dgb$. Beware that Figure
\ref{fig:all_encompassing} shows an upside-down picture of the world, that is,
if a degree $\dga$ lies below a degree $\dgb$, then
$\dga$ is depicted above $\dgb$ in the figure.

Figure \ref{fig:all_encompassing} shows that the degree of the representation by left
best approximations is incomparable to, that is, lies neither above nor below, 
the degree of the representation by right
best approximations. That this indeed is the case can
be established  by a growth argument (see page \pageref{growthpage}). 
Let us see how such an  argument works.

 In Section \ref{skitteskebest}
we saw that the representation by dual Baire sequences is subrecursively equivalent
to the representation by left best approximations and that the representation by standard
Baire sequences is subrecursively equivalent
to the representation by right best approximations.
In order to make our growth
argument transparent, we will consider dual and standard Baire sequences
in place of left and right best approximations. Let $\xcb$ and $d\xcb$ denote 
the representations 
by standard and dual Baire sequences, respectively. We will argue
that $\xcb \not\redrel d\xcb$ and $d\xcb \not\redrel \xcb$.

Let $s$ be any time bound. We will prove that there exists an irrational $\alpha$
in the interval $(0,1)$ such that
\begin{align} \label{enrodhvit}
    \mbox{the standard Baire sequence of $\alpha$ is not computable in time $O(s)$}
\end{align}
but still
\begin{align} \label{torodhvit}
    \mbox{the dual Baire sequence of $\alpha$ is computable in time $\poly{2^n}$.}
\end{align}
where $n$ is the length of the input.

Consider  a
very fast increasing function $B:\nat \longrightarrow \nat$  
with a simple graph. Specifically,
we assume  that the graph of $B$,
that is the relation $B(x)=y$, 
is decidable in time $\poly{\max(x,y)}$, but still,
 $B$ increases too fast to be computable in time $O(s)$. Such a $B$
will always exist, and for convenience,  we will also choose $B$ so that
 $2^x\leq B(x)$ and  $B(x)\leq B(x+1)$.
Now, $B$ is the standard Baire sequence of some irrational number $\alpha$,
and since  an irrational number only has one  standard Baire sequence, the standard Baire sequence of $\alpha$ is not computable in time $O(s)$. Thus,  (\ref{enrodhvit}) holds.
It remains to prove that (\ref{torodhvit}) holds.

Let $a_n = B(0) + (\sum_{i=1}^n B(i) + 1)$. Let $A(x)=1$ if  $x=a_n$ for some $n$; otherwise, let $A(x)=0$. This defines a function $A$. We will prove that $A$ is the dual Baire sequence
of $\alpha$, but first we will argue that $A(x)$ is computable in time $\poly{2^{\bitlen{x}}}$:
Observe that $n$ will be smaller than $x$ whenever $x=a_n$ holds.
Thus, we can check if there exists $n$ such that $x=a_n$ by checking 
\begin{multline} \label{brillefint}
   (\exists n, y_0, \ldots , y_n < x) [ \ B(0)=y_0 \; \wedge \; \ldots 
   \; \wedge \; B(n)=y_n \\ \; \wedge \;  
   x = y_0 + (\sum_{i=1}^n y_i + 1) \ ]
\end{multline}
By assumption we can decide in 
time $\poly{\max(z_1,z_2)}$ if the relation $B(z_1)=z_2$ holds.
Thus we can check in time  $\poly{x}$ if (\ref{brillefint}) holds.
  This shows that $A(x)$ is computable in time $\poly{x}$, and
  thus in time $\poly{2^{\bitlen{x}}}$, since
$A(x)$ equals $1$,
if (\ref{brillefint}) holds, and $0$ if (\ref{brillefint})  does not hold.

For any natural number $n$, we define the strings  $\sigma_n$ and  $\tau_n$ by
$$\sigma_n = 0^{A(0)}1 0^{A(1)}1 \ldots 0^{A(a_n-1)}1 0^{A(a_n)}  \mbox{ and } 
\tau_n = 1^{B(0)}0 1^{B(1)}0 \ldots 1^{B(n)}0\; . $$
 We  prove by induction on $n$ that $\sigma_n  = \tau_n$. 
Let $n=0$. We have $a_0 = B(0)$, and thus, by the definition of $A$, we have  
$$
\sigma_0 \; = \; 0^{A(0)}1 0^{A(1)}1 \ldots 0^{A(a_0-1)}1 0^{A(a_0)} \; = \;  1^{a_0}0 \; = \;  1^{B(0)}0 \; = \; \tau_0\; .
$$
Let $n>0$. By the definition of $a_n$, we have $a_n = a_{n-1}+ B(n) + 1$, and thus $B(n)= a_n - (a_{n-1} + 1)$. Furthermore, we have
\begin{multline*}
\sigma_n \; \stackrel{\mbox{\scriptsize (1)}}{=} \;   \sigma_{n-1}10^{A(a_{n-1}+1)}10^{A(a_{n-1}+2)} \ldots 10^{A(a_n-1)}1 0^{A(a_n)}
 \; \stackrel{\mbox{\scriptsize (2)}}{=} \; \\
 \sigma_{n-1}1^{ a_n - (a_{n-1} + 1)} 0  \; \stackrel{\mbox{\scriptsize (3)}}{=} \; 
 \sigma_{n-1}1^{B(n)} 0 \; \stackrel{\mbox{\scriptsize (4)}}{=} \; \tau_{n-1}1^{B(n)} 0 
\; \stackrel{\mbox{\scriptsize (5)}}{=} \; \tau_{n}
\end{multline*}
where (1) holds by the definition of $\sigma_n$; (2) holds by the definition of $A$; (3) holds by the definition of $a_n$;
(4) holds by the induction hypothesis; and (5) holds by the  definition of $\tau_n$. This proves 
that $\sigma_n=\tau_n$ for any $n$, and  by
 Lemma \ref{lysnisse}, we have 
$$J^{a_n}_A \; = \; \fareypairtree(\sigma_n)  \; = \;  
\fareypairtree(\tau_n) \; = \;    I^n_B$$ 
for any $n$. By the definition of standard and dual Baire sequences (Definition \ref{juleto}),
it follows that
$A$ is the dual Baire sequence of $\alpha$. This completes our proof of (\ref{torodhvit}).

It follows from (\ref{enrodhvit}) and (\ref{torodhvit}) that we have
$$
\poly{2^n}_{d\xcb} \; \not\subseteq \; O(s)_{\xcb}
$$
for any time bound $s$, and thus we have $\xcb \not\redrel d\xcb$ by our  definition of 
of the ordering relation $\redrel$ (Definition \ref{snartsommerigjen}). A symmetric 
proof yields 
$d\xcb \not\redrel \xcb$.

Our growth argument shows that the degree of the representation 
by left best approximations (which is also the degree of 
$d\xcb$) 
is incomparable to the degree of the representation 
by right best approximations (which is also the degree of 
$\xcb$). 
We have seen that we can subrecursively compute the Dedekind cut  of $\alpha\in (0,1)$ if
we have access to a
 left, or to a right, best approximation of $\alpha$ (Lemma \ref{ineedthislabel}).
Thus our two degrees of representations by best approximations will both lie
 above the degree of the representation by Dedekind cuts,
and since the two degrees are incomparable, they
have to lie {\em strictly} above. 
We cannot subrecursively 
 convert the representation by Dedekind cuts into the representation by
 left, or the representation by right, best approximations.
In the next section we will see that both degrees lie below the degree of the representation
by continued fractions, and since the two degrees are incomparable, we can conclude that
they lie strictly below. We cannot subrecursively convert 
the representation by left, or the representation by right, best approximation  into the representation by continued fractions (but we will see in the next section 
that we indeed can  subrecursively  compute the continued fraction of $\alpha$ if  
we have access to {\em both a left and a right} best approximation of $\alpha$).
See Figure \ref{fig:all_encompassing}.

%% file: continued.tex
\section{Representations Equivalent to Continued Fractions}

\label{skittenskjecont}

\subsection{Continued fractions.}

We may assume some familiarity with continued fractions, but  we will state and explain some of
their  properties  below.
For more on continued fractions see Khintchine \cite{khi}
or Richards \cite{richards}. The latter is a very readable paper which carefully 
explains 
the relationship between continued fractions and Farey pairs.

Let  $a_0, a_1, a_2, \ldots $ be an infinite sequence of integers where
$a_1, a_2, a_3 \ldots $ are positive. 
The   {\em continued fraction} $[a_0; a_1, a_2, \ldots ]$ is defined by
 $$[ \ a_0; a_1, a_2,a_3 \ldots \ ]  \; = \;
a_0 + \frac{\displaystyle 1}{\displaystyle a_1 + \frac{\displaystyle 1}{\displaystyle a_2 + \frac{1}{a_3 + \ldots}}}
$$

We will work with continued fraction representations of 
irrational numbers between 0 and 1.
Every irrational number $\alpha$ in the interval  $(0,1)$ can be written uniquely of the form
$\alpha  = [  0; a_1, a_2,  \ldots ]$
where $a_1, a_2, a_3, \ldots $ are positive integers. Hence, the next definition makes sense.

\begin{definition}
Let $\alpha$ be an irrational number in the interval $(0,1)$, and let $\alpha  = [  0; a_1, a_2,  \ldots ]$. The {\em continued fraction} of $\alpha$ is the function $C:\nat^+ \longrightarrow \nat$
given by $C(i)= a_i$.
\qed
\end{definition}

We define $[a_0; a_1,\ldots , a_n ]$ by induction on $n$. If $n=0$, let $[a_0; \ ]= a_0$. If $n>0$, let 
$$[ \ a_0; a_1,\ldots , a_n \ ] \;\; = \;\; a_0 \;  + \; \frac{1}{    [ \ a_1; a_2,\ldots , a_n \ ]     } \; .$$
The rational number $[a_0; a_1,\ldots , a_n ]$ is known as the  $n$th \emph{convergent}
of the infinite continued fraction $[a_0; a_1, a_2, \ldots ]$.

 Let
\begin{align} \label{convdef}
\frac{p_0}{q_0} = \frac{0}{1} \; , \;\;\;   \frac{p_1}{q_1}=\frac{1}{a_1} \;\;\; \mbox{ and } \;\;\;
\frac{p_{k+2}}{q_{k+2}} = \frac{p_{k} + a_{k+2}p_{k+1}}{  q_{k} + a_{k+2}q_{k+1}    }\; .
\end{align}
It is well known that $p_n/q_n$ equals the $n$th convergent of 
$[  0; a_1, a_2, a_3, \ldots  ]$, that is, $p_n / q_n = [  0; a_1,  \ldots , a_n  ]$.
It is also well known that
$$\frac{p_0}{q_0} < 
\frac{p_2}{q_2} \;\; < \;\; \frac{p_{4}}{q_{4}} \;\; < \;\; 
\ldots \;\; < \;\; \alpha \;\; < \ldots \;\; < \;\; \frac{p_{5}}{q_{5}}
 < \;\; 
\frac{p_3}{q_3} \;\; < \;\;   \frac{p_{1}}{q_{1}}
$$
where $\alpha = [0; a_1, a_2, \ldots ]$, that is,  the even convergents approach the number
represented by the continued fraction from below whereas the odd convergents
approach it from above. Every convergent is a (left or right) best approximant,
but the converse  it not true, a best approximant will not necessarily be a convergent.
The next theorem relates the convergents of a continued fraction and the Farey pair tree $\fareypairtree$.

\begin{theorem}[Hurwitz \cite{hurwitz}] \label{juleferieslutt}
Let $[0; a_1, a_2, \ldots]$ be a continued fraction, and let the convergent $p_n/q_n$ be defined by (\ref{convdef}) above. We have
$$
\fareypairtree( \ 0^{a_1} 1^{a_2} \ldots 1^{a_{n-1}} 0^{a_n}  \ )\; = \; \left( \ \frac{p_{n-1}}{q_{n-1}} , \frac{p_{n}}{q_{n}} \ \right) \; = \; [0; a_1,\ldots , a_n ]
$$
when $n$ is odd, and we have
$$
\fareypairtree( \ 0^{a_1} 1^{a_2} \ldots 0^{a_{n-1}} 1^{a_n}  \ )\; = \; \left( \ \frac{p_{n}}{q_{n}} , \frac{p_{n-1}}{q_{n-1}} \ \right) \; = \; [0; a_1,\ldots , a_n ]
$$
when $n$ is even.
\end{theorem}

\subsection{Definitions.}
We will now define some representations which turn out to be subrecursively equivalent
  to the representation by continued fractions.

\begin{definition}
A function $T: [0,1] \cap \rational \longrightarrow \rational$  is  a {\em trace function} for  the irrational number $\alpha$ if we have
$\left\vert \alpha - r \right\vert >  \left\vert \alpha - T(r) \right\vert$
for every $r \in [0,1] \cap \rational$
\qed
\end{definition}

We will say that a trace function $T$ {\em moves $r$ to the right (left)} if $r<T(r)$ ($T(r)<r$). 
The easiest way to realize that a trace function indeed defines a unique real number,
is probably to observe that a trace function $T$ for $\alpha$ yields the Dedekind cut of $\alpha$: if $T$ moves $r$ the right, 
then we know that $r$ lies
below $\alpha$; if $T$ moves $r$ the left, then we know that $r$ lies
above $\alpha$. Obviously, $T$ cannot yield the Dedekind cut for any other number than $\alpha$.

Intuitively, a {\em contractor} is a function that moves two (rational) numbers closer to each other. We will see that also  contractors
can be used to represent irrational numbers.

\begin{definition} \label{ordspill}
A function $F:[0,1] \cap \rational \longrightarrow \rational$ is a {\em contractor} if we have $F(r) \neq r$ and 
$\vert F(r_1)- F(r_2) \vert < \vert r_1 - r_2 \vert$
for any rationals $r,r_1,r_2$ where $r_1\neq r_2$.
\qed
\end{definition}

\begin{lemma}\label{nuerlittoverseks}
Any contractor is a trace function for some irrational number.
\end{lemma}

\begin{proof} 
Let $F$ be a contractor.
 If $F$ moves $r$ to the right (left), then $F$ also moves any rational less (greater) than $r$ to the right (left); otherwise $F$ would
not be a contractor.
We define two sequences $r_0,r_1,r_2 \ldots$ and $s_0,s_1,s_2 \ldots$ of rationals.
Let $r_0 = 0$ and $s_0=1$. Let $r_{i+1} = (r_{i}+s_i)/2$ if $F$ moves $(r_{i}+s_i)/2$ to the right; otherwise, let
$r_{i+1} = r_{i}$.  Let $s_{i+1} = (r_{i}+s_i)/2$ if $F$ moves $(r_{i}+s_i)/2$ to the left; otherwise, let
$s_{i+1} = s_{i}$ (Definition \ref{ordspill} requires that a contractor moves any rational number). 
Obviously, we have 
$\lim_{i} r_i  =  \lim_{i} s_i$,
and obviously, this limit is an irrational number $\alpha$. It is easy to see that $F$ is a trace function for $\alpha$.   
\end{proof}

The previous lemma shows that the next definition makes sense.

\begin{definition}
A contractor $F$ is a {\em contractor for the irrational number $\alpha$} if $F$ is a trace function for $\alpha$.
\qed
\end{definition}

The representation by trace functions was introduced 
in Kristiansen \cite{ciejouren}, and the representation by
contractors was introduced in Kristiansen \cite{bairelars}.

\subsection{From complete best approximation to continued fractions.}

By Theorem \ref{juleferieslutt}, the continued fraction $[0; x_1, x_2, \ldots ]$ of
$\alpha \in (0,1)$ can be viewed as the infinite path $0^{x_1}1^{x_2}0^{x_3}1^{x_4}\ldots$ in the Farey pair tree $\fareypairtree$. By Lemma \ref{farey=approximants},
we can construct the path $0^{x_1}1^{x_2}0^{x_3}1^{x_4}\ldots$ if we have access to
the complete left and the complete right best approximations of $\alpha$.
This  insight yields an
algorithm for converting  complete best approximations 
into a continued fraction. The algorithm, which is given by pseudocode in Figure \ref{fig:pseudocode},  uses the complete left best approximation of $\alpha$, denoted 
$\{ a_i/b_i \}_{i\in \nat}$, and the complete right best approximation of $\alpha$, denoted 
$\{ c_i/d_i \}_{i\in \nat}$, as oracles.  The algorithm  outputs the $n$th
element $x_n$ of $\alpha$'s   continued fraction 
$[0; x_1, x_2, \ldots ]$ (the input is $n$). 
The comments embraced by $\texttt{(*} \ldots \texttt{*)}$  explain how the algorithm works,
and $m(a/b, c/d)$ denotes the mediant of the fractions $a/b$ and $c/d$, that is,
$m(a/b, c/d)= (a+c)/(b+d)$.

\begin{figure}
    \centering
    \renewcommand{\arraystretch}{1.6}
$$\begin{array}{ll}
&\ell:= 0;\; r:=0;\;  i:=1; \; x_i:= 1;\\
& \mbox{\texttt{(*} we have $a_\ell/b_\ell = 0/1$ and $c_r/d_r = 1/1$ \texttt{*)}}\\ 
\mbox{\sl more:} & \mbox{\texttt{(*} $\alpha$ lies in  $(a_\ell/b_\ell , c_r/d_r)$ \texttt{*)}}\\ 
&\texttt{WHILE}\;\;\; {\displaystyle m \left( \frac{a_\ell}{b_\ell} , \frac{c_r }{d_r} \right) = \frac{c_{r+1}}{d_{r+1}}}\;\;\; \texttt{DO}\;\;\; \\ &
 \;\;\; \;\;\;  \texttt{BEGIN} \;
 \mbox{\texttt{(*} $\alpha$ lies in  $(a_\ell/b_\ell , c_{r+1}/d_{r+1})$ and the path branches left \texttt{*)}} \\ &
   \;\;\; \;\;\;  r:= r+1; \; x_i:= x_i+1  \\ &
     \;\;\; \;\;\;  \texttt{END}\;\; \mbox{\texttt{(*} we have determined $x_i$ where $i$ is odd \texttt{*)}} \\ &
\texttt{IF}\;\;\; i=n \;\;\; \texttt{THEN}\;\;\; \langle\mbox{give output $x_i$ and halt}\rangle\\ &
\mbox{\texttt{(*} $\alpha$ lies in  $(a_{\ell+1}/b_{\ell+1} , c_{r}/d_{r})$ since the loop terminated \texttt{*)}}\\ &
  i:=i+1;\; x_i:= 1; \; \ell:= \ell + 1; 
  \mbox{\texttt{(*} $\alpha$ lies in  $(a_\ell/b_\ell , c_r/d_r)$ \texttt{*)}}\\ &
\texttt{WHILE}\;\;\; {\displaystyle  
m \left(\frac{a_\ell}{b_\ell},\frac{c_r}{d_r}\right) = \frac{a_{\ell+1}}{b_{\ell+1}}}\;\;\; 
\texttt{DO}\;\;\; \\ &
\;\;\; \;\;\;  \texttt{BEGIN} \;
 \mbox{\texttt{(*} $\alpha$ lies in  $(a_{\ell+1}/b_{\ell+1} , c_{r}/d_{r})$ and the path branches right \texttt{*)}} \\ &
   \;\;\; \;\;\;  \ell:= \ell +1; \; x_i:= x_i+1  \\ &
     \;\;\; \;\;\;  \texttt{END} \;\; \mbox{\texttt{(*} we have determined $x_i$ where $i$ is even \texttt{*)}} \\ &
\texttt{IF}\;\;\; i=n \;\;\; \texttt{THEN}\;\;\; \langle\mbox{give output $x_i$ and halt}\rangle\\ &
\mbox{\texttt{(*} $\alpha$ lies in  $(a_{\ell}/b_{\ell} , c_{r+1}/d_{r+1})$ since the loop terminated \texttt{*)}}\\ &
 i:=i+1;\; x_i:= 1; \; r:= r + 1; 
 \mbox{\texttt{(*} $\alpha$ lies in  $(a_\ell/b_\ell , c_r/d_r)$ \texttt{*)}}\\ &
 \texttt{GOTO}\;\;\; \mbox{\sl more}
\end{array}
$$
\renewcommand{\arraystretch}{1.0}

    \caption{ALGORITHM}
    \label{fig:pseudocode}
\end{figure}

The first while-loop in Figure \ref{fig:pseudocode}  counts consecutive zeros
found in the path $0^{x_1}1^{x_2}0^{x_3}1^{x_4}\ldots$ given by the best approximations,
whereas the second while-loop counts consecutive ones. 
It turns out that the actual counting is superfluous
and that these two loops can  be eliminated: We can directly
compute the number of times a while-loop will be executed from values available when the loop starts.
Let us consider the first while-loop:
$$
\texttt{WHILE}\;\;\; {\displaystyle m \left(\frac{a_\ell}{b_\ell},\frac{c_r}{d_r}\right)= \frac{c_{r+1}}{d_{r+1}}}\;\;\; \texttt{DO}\;\;\;
   \texttt{BEGIN} \;\;\; r:= r+1; \; x_i:= x_i+1 \;\;\; \texttt{END} 
$$
No variables except $r$ and $x_i$ are modified during the execution of this loop.
Let
\begin{itemize} 
\item  $S$ be the value of $r$ when the execution of the loop starts
\item  $T$ be the value of $r$ when the execution of the loop terminates
\item  $X$ be the value of $x_i$ when the execution of the loop terminates
(observe that $x_i$ is 1 when the execution  starts, and thus $X=1+T-S$).
\end{itemize} 
For  $r = S,\ldots, T-1$, we have 
$$
\frac{c_{r+1}}{d_{r+1}} \;\; = \; \;  m \left(\frac{a_\ell}{b_\ell},\frac{c_r}{d_r}\right)
\;\; = \; \;  \frac{a_\ell + c_{r}}{b_\ell + d_{r}}\; .
$$
Thus, when the loop terminates, we have
\begin{align} \label{badevekt}
 \frac{c_{T}}{d_{T}}
\;\; = \; \;  \frac{a_\ell (X-1) +  c_{S}}{b_\ell (X-1) + d_{S}}\; .
\end{align}
When the loop terminates, we also have 
$$
\frac{c_{T+1}}{d_{T+1}} \;\; \neq \; \; m \left(\frac{a_\ell}{b_\ell},\frac{c_T}{d_T}\right)
$$
but then,
as the mediant $m(a_\ell / b_\ell , c_T / d_T)$ is not the next fraction in the complete list of right best approximants, 
it will be fraction number $\ell+1$ in the complete list of left best approximants. That is, we have:
\begin{align}
 \label{passeroglinjal}
\frac{a_{\ell+1}}{b_{\ell+1}} 
\;\; = \;\; m \left(\frac{a_\ell}{b_\ell},\frac{c_T}{d_T}\right) \;\; = \; \;
  \frac{a_\ell+ c_T}{b_\ell + d_T}\; .
\end{align}
By (\ref{badevekt}) and (\ref{passeroglinjal}), we have
\begin{align}
 \label{xpasseroglinjal}
\frac{a_{\ell+1}}{b_{\ell+1}} 
\;\; = \; \;  \frac{a_\ell+ a_\ell (X-1) +  c_{S}}{b_\ell + b_\ell (X-1) + d_{S}}\; .
\end{align}
Now, (\ref{xpasseroglinjal}) yields the equation
$$
%ABA: commented  out equation which requires a_\ell > 0
%X = \frac{a_{\ell + 1}   - c_S}{a_\ell} \;\; \mbox{ and } \;\; 
X = \frac{b_{\ell + 1}   - d_S}{b_\ell}
$$
where $d_S$ is the value of $d_r$ when the execution of the loop starts and 
$X$ is the value of $x_i$ when the execution of the loop terminates.
Hence, the loop above can be replaced by 
$$
%x_i:= (a_{\ell + 1}   - c_r)/a_\ell; \; r:= r+(x_i-1)
x_i:= (b_{\ell + 1}   - d_r)/b_\ell; \; r:= r+(x_i-1)
$$
(recall that $x_i$ is 1 when the loop starts and
hence $r$ will be incremented exactly $x_i-1$ times before the loop terminates).
A symmetric argument shows that the loop
$$
\texttt{WHILE}\;\;\; {\displaystyle m \left(\frac{a_\ell}{b_\ell},\frac{c_r}{d_r}\right)= \frac{a_{\ell+1}}{b_{\ell+1}}}\;\;\; 
\texttt{DO}\;\;\; \texttt{BEGIN} \;\;\; \ell:= \ell+1; \; x_i:= x_i+1 \;\;\; \texttt{END} 
$$
can be replaced by the program
$$
x_i:= (d_{r + 1}  - b_\ell)/d_r; \; \ell:= \ell + (x_i -1) .
$$
Thus, we have the algorithm in Figure \ref{fig:morepseudocode}.

\begin{figure}
    \centering
$$\begin{array}{cll}
\scriptstyle 1&&\ell:= 0;\; r:=0;\;  i:=1; \\
%&a_0,b_0,c_0,d_0 := 0,1,1,1;\\
\scriptstyle 2&\mbox{\sl more:} & x_i:= (b_{\ell + 1}   - d_r)/b_\ell; \; r:= r+(x_i-1); \\ \scriptstyle 3&&
\texttt{IF}\;\;\; i=n \;\;\; \texttt{THEN}\;\;\; \langle\mbox{give output $x_i$ and halt}\rangle\\ \scriptstyle 4&&
  i:=i+1;\;  \ell:= \ell + 1;\\ \scriptstyle 5&&
  x_i:= (d_{r + 1}  - b_\ell)/d_r; \; \ell:= \ell + (x_i -1); \\ \scriptstyle 6&&
\texttt{IF}\;\;\; i=n \;\;\; \texttt{THEN}\;\;\; \langle\mbox{give output $x_i$ and halt}\rangle\\ \scriptstyle 7&&
 i:=i+1;\; r:= r + 1;\\ \scriptstyle 8&&
 \texttt{GOTO}\;\;\; \mbox{\sl more}
\end{array}
$$
\renewcommand{\arraystretch}{1.0}
    
    \caption{ALGORITHM}
    \label{fig:morepseudocode}
\end{figure}
\renewcommand{\arraystretch}{1.6}

\begin{lemma}\label{lem:completebest_continued}
Let $L: \mathbb{N} \longrightarrow \mathbb{Q}$ 
and $R: \mathbb{N} \longrightarrow \mathbb{Q}$ be the complete left and right, repstectively,  best approximations of an irrational $\alpha \in (0,1)$. Assume that both $L$ and $R$
are computable within a time-bound $s$, and let $f(n)= \lambda x. (s(1+x))^{(n)}(0)$.
There is a parameterized function-oracle Turing machine $M$ such that
\begin{itemize}
    \item $\Phi_\paramorac{M}{L,R}: \mathbb{N}^+ \longrightarrow \mathbb{N}$ is the continued
fraction  of $\alpha$
\item $M^{L,R}$ on input $n$ runs in time $\poly{f(n)}$ and
uses exactly $2n$ oracle calls, each of input size at most $f(n-1)$.
\end{itemize}
\end{lemma}

\begin{proof}
The Turing machine $M$ will iterate the 
 body of the loop implemented by the goto-statement in Figure \ref{fig:morepseudocode}  no more than $\lceil n/2 \rceil$ times. 
 To implement Line~2,
 it asks the oracle $L$ for the values of $a_\ell/b_\ell$ and $a_{\ell+1}/b_{\ell+1}$.
 Then $M$ can carry out the assignment as the value of $d_r$ is already known to $M$ and will be stored
 at the work tape (except for the very first time Line~2 is reached, then we have $r=0$ and
 $d_r=1$ by convention). 
 Similarly, to implement Line~5, it 
  asks the oracle $R$ for $c_{r}/d_{r}$ and $c_{r+1}/d_{r+1}$.
 The value of $b_\ell$ is already known to $M$. 
 %(The oracles will provide the values of 
 %the nominators even if $M$ does not need them.)
 Thus, $M$ consults each oracle four times each time the loop's body is iterated, and
 when $M$ outputs $x_n$ and halts, exactly $2n$ oracle calls have been made.
 
 The time required to perform the arithmetic operations is 
polynomial in the size of the (representation of) the involved values.  We do an induction to 
prove bounds on the  bit-lengths of $\ell,r$ throughout the algorithm.
Our induction claim is the conjunction of the following two statements:
(1) Whenever Line~6 is reached, we have
\begin{equation}
    \label{eq:ind:ell}
    \ell \; \le \;  d_{r+1} \;\;\; \text{ and } \;\;\; \bitlen{\ell} \;\le 
    \; s( 1 + f(i-1)) 
\end{equation}
and (2) whenever Line~3 is reached, we have
\begin{equation}
    \label{eq:ind:r}
    r \; \le \; b_{\ell+1} \;\;\; \text{ and }
    \;\;\; \bitlen{r} \; \le \; s( 1 + f(i-1))  \; .
\end{equation}

Note that these bounds, which we prove next, complete the justification of the lemma regarding the time complexity
and oracle input size.

Consider the first tour through the loop body. We reach Line~3 with
$$
 r = (b_1 - d_0)/b_0 - 1 = (b_1 -1)/1 - 1 = b_1 - 2 \le b_1 
$$
and we reach line~6 with 
$\ell = (d_1 - b_1)/d_0 = d_1 - b_1 \le d_1$.
It is easy to verify that both \eqref{eq:ind:ell} and \eqref{eq:ind:r} hold.

In general, suppose that we have reached Line~6 for the
$k$th time. Now $i=2k$.
Moving on towards Line~2, we increment $i$ and $r$.
The value of $x_i$ computed at Line~2 is bounded by
$$
 x_i  \; = \; (b_{\ell+1} - d_r)/b_\ell \; \le \; b_{\ell+1} - b_\ell ,
$$
where the last inequality uses the fact that $(x-1)/y \le x-y$ holds
for all integers where $x>y>0$.
We compare the new value of $r$ after Line~2 to the value $r'$ which it had when
we last visited Line~3. We have
$$
 r \; = \; r' + x_i \le b_{\ell'+1} + b_{\ell+1} - b_\ell \; \le \; b_{\ell+1}
$$
where $\ell'$ is the value $\ell$ had when we last visited Line~3 (clearly $\ell' \le \ell-1$).

This proves that $r \le  b_{\ell+1}$, and hence the first conjunct of (\ref{eq:ind:r}) holds
when Line~3 is reached.
In order to verify that the second conjunct of (\ref{eq:ind:r}) also holds,
observe that $b_{\ell+1}$ is obtained by an oracle query, and hence we have
$\bitlen{b_{\ell+1}} \leq s(\bitlen{\ell+1})$. By our induction hypothesis \eqref{eq:ind:ell},
we have $\bitlen{\ell} \leq s(1 + f(i-2))$ when Line~3 is reached.
Thus, as $r \le  b_{\ell+1}$, we have
\begin{multline*}
 \bitlen{r } \; \leq  \;  \bitlen{b_{\ell+1}} \; \leq \; s(\bitlen{\ell+1}) \; \leq \;
 s(1 +\bitlen{\ell})
  \; \leq \; \\  s(1 + s(1 + f(i-2))) \;  = \; s(1 + f(i-1)) 
\end{multline*}
(the final equality holds by the definition of $f$).
 A symmetric argument justifies \eqref{eq:ind:ell} in the inductive case.
\end{proof}

\subsection{From continued fractions to complete best approximations.}

Our algorithm for converting a continued fraction into a complete left (or right)
best approximation is pretty straightforward:
Let $ [0; x_1, x_2, \cdots ]$ be the continued fraction of $\alpha$.
By Theorem \ref{juleferieslutt}, $0^{x_1}1^{x_2}0^{x_3} \cdots$
is the unique path of $\alpha$ in the Farey pair tree $\fareypairtree$.
By Lemma \ref{farey=approximants}, complete  left and right approximations of $\alpha$
can be read off, in order, from  the Farey pairs along the path
$0^{x_1}1^{x_2}0^{x_3} \cdots$.  Every time a 1 occurs  in the path (branching right),
a new left best approximant will show up; every time a 0 occurs in the path (branching left),
a new right best approximant will show up.

\begin{lemma}\label{lem:cont_to_complete}
Let $C: \nat^+ \longrightarrow \nat$ be the continued fraction of the
irrational number $\alpha\in (0,1)$.
There is a parameterized function-oracle Turing machine $M$ such that 
\begin{itemize}
    \item  $\Phi_\paramorac{M}{C}: \nat \longrightarrow \mathbb{Q}$
is the complete left best approximation  of $\alpha$
\item $M^C$ on input $n$ runs in time  $\poly{\sum_{j=1}^{2n} C(j)}$ and at most $2n$ oracle
calls, each of input size at most $3 + \log n$ and output size at most $\max_{j=1}^{2n}\{\log C(j)\}$.
\end{itemize}
\end{lemma}

\begin{proof}
Let 
 $\underline{\sigma}$ denote the path $0^{C(1)}1^{C(2)}0^{C(3)} \cdots$.
 The $n$th element of the complete best left approximation is  obtained by 
finding the index $i_n$ such that the $n$th occurrence of $1$ in $\underline{\sigma}$
 occurs at 
index $i_n$. Thereafter $M$  computes the Farey pair at position $\underline{\sigma}_1 \cdots \underline{\sigma}_{i_n}$ and simply returns that pair's left component.

Observe that as $C(1),C(2),\ldots \geq 1$, we have $i_n \leq \sum_{j=1}^{2n} C(j)$.
Hence, $M$  needs to perform at most $2n$ oracle calls of input size
at most $1 + \lfloor \log 2n \rfloor \leq 3 + \log n$ 
(as the input needs to represent numbers of size at most $2n$)
and output size at most $\max_{j=1}^{2n}\{\log C(j)\}$. After having found $i_n$,
$M$ determines  the Farey pair at the position $\underline{\sigma}_1 \cdots \underline{\sigma}_{i_n}$. This can be accomplished in time $\poly{i_n} \leq
\poly{\sum_{j=1}^{2n} C(j)}$ by Proposition \ref{Hurwitzburdeslankesig}. 
Now it is easy to see that
the total time needed to compute the value of the oracle calls, reading off the results,
and computing the relevant Farey pair is bounded above by $\poly{\sum_{j=1}^{2n} C(j)}$.
\end{proof}

\subsection{Conversion between contractors, trace functions and complete best approximations.}

\begin{lemma}\label{krymper}
(i) Let  $\{r_i\}_{i\in \nat}$ be a complete left best approximation  of $\alpha$.
For any $i\in\nat$, we have 
$$
r_{i+1} - r_i \; > \; r_{i+2} - r_{i+1}\; .
$$
(ii) Let  $\{r_i\}_{i\in \nat}$ be a complete right best approximation  of $\alpha$.
For any $i\in\nat$, we have 
$$
r_{i} - r_{i+1} \; > \; r_{i+1} - r_{i+2}\; .
$$
\end{lemma}

 \newcommand{\mma} {\mathbf{a}}
 \newcommand{\mmb} {\mathbf{b}}

\begin{proof}
We prove (ii); the proof of (i) is symmetric.  By Lemma \ref{farey=approximants}, $r_i$ will be the right endpoint 
of some interval in the tree $\fareypairtree$ (we can w.l.o.g.~assume that $r_i> 0/1$). 
Thus, we have $\sigma\in \{0,1\}^*$ and fractions $a/b$ and $c/d$ such that 
$\fareypairtree(\sigma 0)= (a/b,c/d)$ and  $r_i = c/d$.

By Proposition \ref{slankemeg} and Lemma \ref{farey=approximants}, there exists $k$ such that
\begin{align} \label{fillijonka}
\fareypairtree(\sigma 0 1^k 0)  =  \left( \ \frac{a+ kc}{b+ kd}\; , \;  \frac{a+ kc + c}{b+ kd + d}  \ \right)
 \;\; \mbox{ and }  \;\;  r_{i+1} = \frac{a+ kc + c}{b+ kd + d}\; .
\end{align}
 Let $\mma = a+ kc$, let $\mmb = b+ kd$. We can now rewrite (\ref{fillijonka}) as
\begin{align} \label{jonka}
\fareypairtree(\sigma 0 1^k 0)  =  \left( \ \frac{\mma}{\mmb}\; , \;  \frac{\mma + c}{\mmb + d}  \ \right)
  \;\; \mbox{ and }  \;\; r_{i+1} = \frac{\mma + c}{\mmb + d}\; .
\end{align}
By (\ref{jonka}) and the definition of $\fareypairtree$, we have 
$$ \fareypairtree(\sigma 0 1^k ) = \left( \ \frac{\mma}{\mmb}  \; , \; \frac{c}{d} \ \right) 
\;\; \mbox{ and } \;\;
\fareypairtree(\sigma 0 1^k 1 ) = \left( \ \frac{\mma+c}{\mmb+d} \; , \;  \frac{c}{d} \ \right)\; .$$ 
This shows that $( (\mma+c) / (\mmb+d), c/d)$ is an interval in $\fareypairtree$.
Thus, by Theorem \ref{slankemeg}, we have 
\begin{align} \label{neinei}
 \frac{c}{d} \; - \;  \frac{\mma + c}{\mmb + d} \;\; = \;\; \frac{1}{d(\mmb + d)} \; .
\end{align}
By Proposition \ref{slankemeg} and Proposition \ref{prop:alleskalslankesig},
there exists $m$ such that $r_{i+2}$ is the right endpoint of the interval 
$\fareypairtree(\sigma 1^k 0 1^m 0)$.
We can assume that $m=0$ since $m=0$
yields the maximal distance between $r_{i+1}$ and $r_{i+2}$.
 Thus, by the definition of $\fareypairtree$, we have 
$$\fareypairtree( \sigma 1^k 00)  = \left( \ \frac{\mma}{\mmb} \; , \; 
 \frac{2\mma + c}{2\mmb + d} \ \right) 
  \;\; \mbox{ and }  \;\; r_{i+2} = \left( \ \frac{2\mma + c}{2\mmb + d} \ \right) \; .$$
 Moreover, again by  the definition of $\fareypairtree$, we have
$$\fareypairtree( \sigma 1^k 0) \; = \;  \left( \ \frac{\mma}{\mmb} \; , \; \frac{\mma + c}{\mmb + d} \ \right) \;\;
\mbox{ and } \;\;
 \fareypairtree( \sigma 1^k 01) \; = \; \left( \ \frac{2\mma + c}{2\mmb +d} \; , \; \frac{\mma + c}{\mmb + d} \ \right)\; .$$
This shows that $((2\mma + c) / (2\mmb +d), (\mma + c)/(\mmb + d)$ is and interval in  $ \fareypairtree$, and thus, by Theorem \ref{slankemeg}, we have 
\begin{align} \label{jaja}
 \frac{\mma + c}{\mmb + d} \; - \;  \frac{2\mma + c}{2\mmb + d} \;\; = \;\; \frac{1}{(\mmb + d)(2\mmb + d)} \; .
\end{align}
Now we can conclude our proof of (ii) with
\begin{multline*}
r_i - r_{i+1} \;\; = \;\; \frac{c}{d} \; - \; \frac{\mma + c}{\mmb + d}  \;\; \stackrel{\mbox{\scriptsize (\ref{neinei})}}{=} \;\;
 \frac{1}{d(\mmb + d)}  \;\; > \;\;  \frac{1}{(\mmb + d)(2\mmb + d)} \;\; \stackrel{\mbox{\scriptsize (\ref{jaja})}}{=} \;\; \\
 \frac{\mma + c}{\mmb + d} \; - \;  \frac{2\mma + c}{2\mmb + d}
 \;\; = \;\;  r_{i+1} -  r_{i+2}\; .
\end{multline*}
\end{proof}

It is necessary to assume in Lemma \ref{krymper} that the best approximations 
are complete. The lemma does not hold for best approximations in general: 
E.g., let $(1000)^{-1}> \alpha > 0$.
Then we have $1 > 1/2 > 1/3 > 1/99 > \ldots > \alpha$
%$$
%1 \; > \; \frac{1}{2} \; > \; \frac{1}{3} \; > \; \frac{1}{99} \; > \; \ldots \; > \; \alpha
%$$
where $1/2$, $1/3$ and $1/99$ all are right best approximants of $\alpha$, but it is false that
$1/2 - 1/3 > 1/3 - 1/99$.
A complete right best approximation  to $\alpha$ is of the form
$$
1 \; > \; \frac{1}{2} \; > \; \frac{1}{3} \; > \; \frac{1}{4} \; > \; \ldots \; > \; 
 \frac{1}{i} \; > \; \frac{1}{i+1} \; > \;   \ldots \; > \;  \frac{1}{1000} \; > \;  \ldots \; > \;  
\alpha \; .
$$

Let  $\{\hat{r}_i\}_{i\in \nat}$ be a complete left best approximation  of $\alpha$, and let
  $\{\check{r}_i\}_{i\in \nat}$ be a complete right best approximation  of $\alpha$.
Furthermore, let 
$$\hat{k}_i= (\hat{r}_{i+2}- \hat{r}_{i+1})/(\hat{r}_{i+1}- \hat{r}_{i}) \;\;\;\mbox{and} \;\;\;
\check{k}_i= (\check{r}_{i+1}- \check{r}_{i+2})/(\check{r}_{i}- \check{r}_{i+1})\; .
$$ 
By the previous lemma we have $\hat{k}_i< 1$ and $\check{k}_i< 1$ (for all $i$), and thus we can define
a contractor $F$ for $\alpha$ by
\begin{align} \label{contractordef}
  F(x)  =  \begin{cases}
  \hat{r}_{i+1} + \hat{k}_i(x - \hat{r}_i)
& \mbox{if  $\hat{r}_{i}\leq x < \hat{r}_{i+1}$ }   \\
 \check{r}_{i+1} + \check{k}_i(x - \check{r}_{i})
& \mbox{if  $\check{r}_{i+1}< x \leq \check{r}_{i}$ .}   
\end{cases}   
\end{align}
It order to verify  that $F$ is indeed is a contractor for $\alpha$, we will prove that we have
\begin{align}\label{hjemmeioslo}
|F(x)- F(y)| \; < \; | x - y| 
\end{align}
whenever $x\neq y$.
The proof splits into several cases. We can w.l.o.g.~assume that $x<y$.

\paragraph{Case (i)} There exist $i,j$ such that $\hat{r}_i\leq x < \hat{r}_{i+1}$ and $\check{r}_{j+1}< y \leq \check{r}_j$.
Then $F$ moves $x$ to the right and $y$ to the left, and hence (\ref{hjemmeioslo}) holds.

\paragraph{Case (ii)} There exists $i$ such that $\hat{r}_i\leq x < y <\hat{r}_{i+1}$.
Then we have 
$|F(x) - F(y)|  =  \hat{k}_i|x-y|$, and (\ref{hjemmeioslo}) holds since $\hat{k}_i<1$.

\paragraph{Case (iii)} There exist $i,j$ (where $i< j$) such that $\hat{r}_i \leq x < \hat{r}_{i+1}$ and 
$\hat{r}_j \leq y < \hat{r}_{j+1}$. 
We write
    $x = \hat{r}_{i+1} - b(\hat r_{i+1} - \hat r_i)$ and
    $y = \hat{r}_j + a(\hat r_{j+1} - \hat r_j)$
where $0\le a,b\le 1$. Then we have
$$y - x = a(\hat r_{j+1} - \hat r_j) + (\hat r_j - \hat r_{i+1}) + b(\hat r_{i+1} - \hat r_i) \;.
$$
Furthermore, we have
\begin{multline*}
    F(x) = \hat r_{i+1} + \hat k_i  (x - \hat{r}_i) 
    = \hat r_{i+2} - (\hat r_{i+2} -    
\hat r_{i+1}) + \hat k_i  (x - \hat{r}_i) \\
    = \hat r_{i+2} - \hat k_i (\hat r_{i+1} -    
\hat r_{i}) + \hat k_i  (x - \hat{r}_i) 
    = \hat r_{i+2} - \hat k_i (\hat r_{i+1} - x) \\
    = \hat r_{i+2} - \hat k_i  b(\hat{r}_{i+1} - \hat{r}_i) 
\end{multline*}
and
    $F(y) = \hat r_{j+1} + \hat k_j  a(\hat{r}_{i+1} - \hat{r}_i)$, and hence  
 $$F(y)-F(x) = \hat k_j a(\hat{r}_{j+1} - \hat{r}_j) + (\hat{r}_{j+1} - \hat{r}_{i+2}) + \hat k_ib(\hat{r}_{i+1} - \hat{r}_i)\; .$$
Now it is easy to verify that $y - x > F(y)-F(x)$.

\paragraph{Case (iv)}  There exists $i$ such that 
$\check{r}_{i+1}< x < y \leq \check{r}_{i}$.
This case is symmetric to (ii). Use that $\check{k}_i<1$.

\paragraph{Case (v)}   There exist $i,j$ (where $i<j$) such that 
$\check{r}_{i+1} < x \leq \check{r}_{i}$ and 
$\check{r}_{j+1} < y \leq \check{r}_{j}$.
This case is symmetric to (iii). Use Lemma \ref{krymper} (ii) in place of Lemma \ref{krymper} (i). This completes the
proof of (\ref{hjemmeioslo}).

\begin{lemma}\label{slutterlittoverseks}
Let $L: \mathbb{N} \longrightarrow \mathbb{Q}$ 
and $R: \mathbb{N} \longrightarrow \mathbb{Q}$ be the complete left and right, repstectively,  best approximation of an irrational $\alpha \in (0,1)$. Let  $L(i) = a_i/b_i$ and $R(i) = c_i/d_i$.
There is a parameterized function-oracle Turing machine 
$M$ such that 
\begin{itemize}
    \item $\Phi_\paramorac{M}{L,R} : [0,1] \cap \rational \longrightarrow \rational$ 
is a contractor for $\alpha$. 
\item $M^{L,R}$ on input $n/m$  runs
in time $O(\log m \log^2\max\{b_m,d_m\})$ and uses at most 
$2\log m$ oracle, each  of input size at most $\log m$.
\end{itemize}
\end{lemma}

\begin{proof}
Define $M$ to be the Turing machine that, when given oracle access
to $L,R$ computes the contractor $F$ given by (\ref{contractordef}).
On input $n/m \in \mathbb{Q}$,
the fact that $L$ and $R$ yield best approximations implies 
that there is an $i$ with
$1 < i < m$ such that
\begin{multline*}
    L(i) \; = \; \frac{a_i}{b_i} \; \leq \; \frac{n}{m} \; < \; 
    \frac{a_{i+1}}{b_{i+1}} \; = \; L(i+1) \\ 
\;\;\; \mbox{ or } \;\;\; R(i+1) \; = \; \frac{c_{i+1}}{d_{i+1}}
\; < \; \frac{n}{m} 
\; \leq  \; \frac{c_i}{d_i}  \; = \; R(i)\; .
\end{multline*}

As $L$ and $R$
are strictly increasing, resp.~decreasing,
$M$ can find $i$ as above by binary search,
requiring at most $\log m$ steps (hence $\log m$ queries to each oracle),
and each step requires 4 comparisons of rational
numbers that can be performed in time
$O(m_0^2)$ where $m_0$ is the length of the largest
binary representation of the integer components of the rational numbers (because $n/m < a/b$
if{f} $nb < am$, and schoolbook multiplication can be done in quadratic time in the size of the representation). The largest integer occurring
in the comparison above is bounded above by
$$
\max\{n,m,a_i,b_i,c_i,d_i\}_{i=1}^{m} \leq \max\{b_m,d_m\}
$$
where the inequality follows as all fractions are bounded above by $1$ and $\{a_i/b_i\}_{i \in \mathbb{N}}$
and $\{c_i/d_i\}_{i \in \mathbb{N}}$
are best approximations.
Hence, the total time needed to compute
$i$ is $O(\log m \log^2\max\{b_m,d_m \})$
using $2\log m$ oracle calls, each of size at most
$\log  m$.
Once $i$ has been found, $M$ first computes $F(n/m)$ 
as a fraction (not necessarily in lowest terms) using
a constant number of additions, subtractions and multiplications of numbers originally representable
by at most $\log\max\{b_m,d_m\}$ bits,
hence using at most
$O(\log^2\max\{b_m,d_m\})$ operations
and resulting in a number representable using
$O(\log\max\{b_m,d_m\})$ bits. Obtaining a reduced fraction can be done by first computing the gcd of the numerator and the denominator,
and then performing the 2 requisite integer divisions, for a total of $O(\log^2\max\{b_m,d_m\})$ further operations. 
\end{proof}

Before we give our algorithm for converting trace functions (and thus also contractors)
into best approximations, we will make a couple of observations. The first observation is trivial:
If  a trace function for $\alpha$ moves
a rational number  to the right (left), then the rational number lies below (above) $\alpha$. Hence, if we have access to
a trace function for $\alpha$, we can easily compute the Dedekind cut  $D^\alpha$ of $\alpha$. The next observation
is slightly more sophisticated: Let $T_0(r)= (r + T(r))/2$ where $T$ is a trace function for $\alpha$.  Then, also $T_0$
will be a trace function for $\alpha$, moreover, we have 
\begin{align}\label{nybygg}
    \mbox{$T_0(r)< \alpha$ if $r< \alpha$}
\end{align}
 and  
 \begin{align} \label{nybygger}
    \mbox{$T_0(r)> \alpha$ if $r> \alpha$.}
\end{align}
 In order to see that (\ref{nybygg}) holds, assume that $r<\alpha$ and
 $T_0(r) = (r + T(r))/2 \geq \alpha$. Then, we have 
$T(r) - \alpha \geq  \alpha - r > 0$, and thus $ \vert T(r) - \alpha \vert = T(r) - \alpha \ge \alpha - r  = \vert \alpha - r \vert$, contradicting that $T$ is a trace function. A symmetric argument shows that (\ref{nybygger}) holds.

 Our algorithm for converting a trace function $T$
for $\alpha$ into the complete left best approximation 
$\{ a_i/b_i \}_{i\in \nat}$
of $\alpha$ uses the Dedekind cut
$D^\alpha$ and the trace function $T_0$. When $n=0$, the algorithm simply lets
$a_n/b_n=0/1$. When $n>0$, the algorithm performs the following steps.
\begin{itemize}
    \item Step 1: Recursively, compute the value $a'/b'$ of $T_0(a_{n-1}/b_{n-1})$.
    \item Step 2: Using $D^\alpha$, search for the least natural number $b''\leq b'$ such that for some $a''$,
    $a_{n-1}/b_{n-1}< a''/b'' < \alpha$.    
    \item Step 3: Using $D^\alpha$, find the greatest $a''< b''$ such that $a''/b'' < \alpha$.
    \item Step 4: Set $a_n/b_n$ to $a''/b''$.
\end{itemize}
 Such $a''$ and $b''$ will for sure exist as, if no $b'' < b'$ satisfies the requirement in Step 2, then
 $b'$ itself satisfies it.
 It is easy to see that $a''/b''$ will be smallest left best approximation to $\alpha$ that is strictly greater
 than $a_{n-1}/b_{n-1}$.
 
 A trace function can of course be converted into a right best approximation
 by a symmetric algorithm.

\begin{lemma}\label{lem:trace_to_best}
Let $T : [0,1] \cap \mathbb{Q} \longrightarrow \mathbb{Q}$ be a trace function
 for an irrational number $\alpha \in (0,1)$.
Assume $T(x)$ is computable in time $s(\bitlen{x})$, and let 
$f(n)= \lambda x .(2s(x))^{(n)}(2)$.
 There is a parameterized function-oracle Turing 
 machine $M$ such
 \begin{itemize}
     \item $\Phi_\paramorac{M}{T}: \nat \longrightarrow \mathbb{Q}$
 is the complete left best approximation  of $\alpha$
 \item $M^T$ on input $n$ runs in time 
$O(2^{ f(n) } + s( f(n) ))$ and uses at most $2^{ f(n) }$
 oracle calls, each  of size at most $f(n)$.
 \end{itemize}
\end{lemma}

\begin{proof}
Observe that computing $T_0(r)$ can be performed with a single oracle call
to $T(r)$ followed by three arithmetical operations, hence in time
polynomial in the size of the representations of $r$ and $T(r)$; the latter is bounded by $s(\bitlen{r})$.
We first bound the size of $b'$: by assumption, $a'/b'$ is computable in $s(\bitlen{a_{n-1}/b_{n-1}})$ time, so this is also a bound on its 
size, and the size of $b''$ is by definition at most that of $b'$; the same goes for $a''$. A straightforward induction on $n$ will show that $\bitlen{a_n/b_n} \leq f(n)$.  

The search for $b''$ in Step~2 requires
at most $b'$ calls to $T$ (in order to compute the Dedekind cut of $\alpha$), 
each with arguments consisting
of a rational with representations  size at most $\bitlen{a'/b'}$.
We use a rough bound on $b'$, namely $2^{\bitlen{a'/b'}}$, to estimate the number of calls.

Regarding execution time: arithmetic operations, tests etc.~add an overhead polynomial in the number of oracle calls, and we should also take into
account that the result of the call may be bigger than its argument and another application of $s$ is necessary to cover this cost. Hence we
arrive at the expression $O(2^{ f(n) } + s(f(n)) )$.
The complexity of the search for $a''$ in Step~3 is bounded in precisely the same way.
This leads to the conclusions stated in the lemma.
\end{proof}

The brute-force search in the proof of Lemma \ref{lem:trace_to_best} can also be performed by using the Stern-Brocot tree, but we have been unable to derive better bounds for this approach.

%The algorithm above can be extended to compute not only the best approximant $a_n/b_n$, but also a
%node $\displaystyle{\left(\frac{a_n}{b_n}, \frac{c}{d} \right)}$ in the tree such that 
%$\alpha\in \displaystyle{\left(\frac{a_n}{b_n}, \frac{c}{d} \right)}$. Then we will also have 
%$a'/b' \in \displaystyle{\left(\frac{a_n}{b_n}, \frac{c}{d} \right)}$, and we can find the number we are looking for,
%that is $a_{n+1}/b_{n+1}$, by using the Dedekind cut $D^\alpha$ and computing mediants until we find a node  
%$\displaystyle{\left(\frac{a_{n+1}}{b_{n+1}}, \frac{c'}{d'} \right)}$ such that 
%$\alpha \in \displaystyle{\left(\frac{a_{n+1}}{b_{n+1}}, \frac{c'}{d'} \right)}$.

%In many cases this tree-search will find $a_{n+1}/b_{n+1}$ considerably faster than a brute force search, but in certain
%cases it will not be much more efficient than a brute force search (as far I can see without thinking too much).

\subsection{Summary.}
We can now give a summary of our results on the complexity of the conversions
presented in this sections
in the same style as
we have given  summaries   of corresponding results presented in earlier sections.

\begin{theorem} \label{millehossille}
Consider the representations by (1)  continued fractions,
(2) complete left best approximation together with  complete right best approximations, (3)  trace functions and (4)  contractors, and let $R_1$ and $R_2$ be
any two of these four representations. Then, for  an arbitrary time-bound $t$,
there exists a time-bound $s$ primitive recursive in $t$ such that
$O(t)_{R_2} \subseteq O(s)_{R_1}$. 
\end{theorem}

The reader familiar with the Grzegorcyk hierarchy can easily check that the time-bound 
$s$ in Theorem \ref{millehossille},
for any $i\geq 3$,
indeed will be in the Grzegorcyk class $\mathcal{E}_{i+1}$ if the time-bound $t$
in the Grzegorcyk class $\mathcal{E}_{i}$.